\documentclass[microtype]{gtpart}     
\usepackage{graphicx}  
\usepackage[all]{xy}

\usepackage{float}
\usepackage{tikz}
\usetikzlibrary{
  knots,
  hobby,
  decorations,
  decorations.pathreplacing,
  shapes.geometric,
  calc,
  intersections
}
\usepackage{caption} 
\usepackage{subcaption}
\usepackage{comment}

\title{Diagrammatic representations of 3-periodic entanglements}

\author{Toky Andriamanalina}
\givenname{Toky}
\surname{Andriamanalina}
\address{University of Potsdam, Institute for Mathematics, Karl-Liebknecht-Str. 24-25, 14476 Potsdam-Golm, Germany}

\author{Myfanwy E. Evans}
\givenname{Myfanwy}
\surname{Evans}
\address{University of Potsdam, Institute for Mathematics, Karl-Liebknecht-Str. 24-25, 14476 Potsdam-Golm, Germany}
\email{evans@uni-potsdam.de}

\author{Sonia Mahmoudi}
\givenname{Sonia}
\surname{Mahmoudi}
\address{Advanced Institute for Materials Research, Tohoku University, 2-1-1 Katahira, Aoba-ku, Sendai 980-8577, Japan;  RIKEN iTHEMS, 2-1 Hirosawa, Wako, Saitama 351-0198, Japan}


\arxivreference{2401.14254}

\theoremstyle{plain}
\newtheorem{theorem}{Theorem}[section]

\newtheorem{proposition}[theorem]{Proposition}
\newtheorem{lemma}[theorem]{Lemma}

\theoremstyle{definition}
\newtheorem{definition}[theorem]{Definition}
\newtheorem{example}[theorem]{Example}
\newtheorem{remark}[theorem]{Remark}

\def\qed{\hfill $\Box$}

\begin{document}

\begin{abstract}    
Diagrams enable the use of various algebraic and geometric tools for analysing and classifying knots. In this paper we introduce a new diagrammatic representation of triply periodic entangled structures (\textit{TP tangles}), which are embeddings of simple curves in $\mathbb{R}^3$ that are invariant under translations along three non-coplanar axes. As such, these entanglements can be seen as preimages of links embedded in the 3-torus $\mathbb{T}^3 = \mathbb{S}^1 \times \mathbb{S}^1 \times \mathbb{S}^1$ in its universal cover $\mathbb{R}^3$, where two non-isotopic links in $\mathbb{T}^3$ may possess the same TP tangle preimage. We consider the equivalence of TP tangles in $\mathbb{R}^3$ through the use of diagrams representing links in $\mathbb{T}^3$. These diagrams require additional moves beyond the classical Reidemeister moves, which we define and show that they preserve ambient isotopies of links in $\mathbb{T}^3$. The final definition of a \textit{tridiagram} of a link in $\mathbb{T}^3$ allows us to then consider additional notions of equivalence relating non-isotopic links in $\mathbb{T}^3$ that possess the same TP tangle preimage.

\end{abstract}

\maketitle


\section{Introduction}\label{sec:1}
A classical knot is an embedding of the circle $\mathbb{S}^1$ in $\mathbb{R}^3$. A knot diagram is a planar projection of a given knot with crossing information. It is an efficient planar representation that helps in understanding the complexity of a knot (see, for example, Adams \cite{Adams.book}). A diagram allows one to translate the equivalence of knots given by ambient isotopies of the surrounding space $\mathbb{R}^3$, to simple moves on a plane called the \textit{Reidemeister moves} (Figure \ref{fig:usual_Reid_moves}) (see Burde et al. \cite{burdezies_chap1}, Murasugi \cite{Murasugi1996_chap4}, Ozsváth et al. \cite{grid_homology_appendix_B}). Knot diagrams enable the description of some numerical knot invariants, such as the crossing number or the unknotting number \cite{Adams.book}, and some algebraic knot invariants such as the Jones or Alexander polynomials (Lickorish \cite{raymondlickorish}), among others, which are used to classify knots and links in space.

\begin{figure}[htbp]
    \centering
    \begin{minipage}[b]{0.4\textwidth}
        \centering
        \begin{tikzpicture}
    \begin{knot}[
        consider self intersections=true,
        flip crossing/.list={2,4},
        line width=5pt, 
        clip width=1.25pt, 
        only when rendering/.style={
        }
    ]
        \strand (2,0) .. controls +(0,1.0) and +(54:1.0) .. (144:2) .. controls +(54:-1.0) and +(18:-1.0) .. (-72:2) .. controls +(18:1.0) and +(162:-1.0) .. (72:2) .. controls +(162:1.0) and +(126:1.0) .. (-144:2) .. controls +(126:-1.0) and +(0,-1.0) .. (2,0);
    \end{knot}
\end{tikzpicture}
    \end{minipage}
   \hspace{0.5cm}
    \begin{minipage}[b]{0.4\textwidth}
        \centering
        \includegraphics[width=0.7\textwidth]{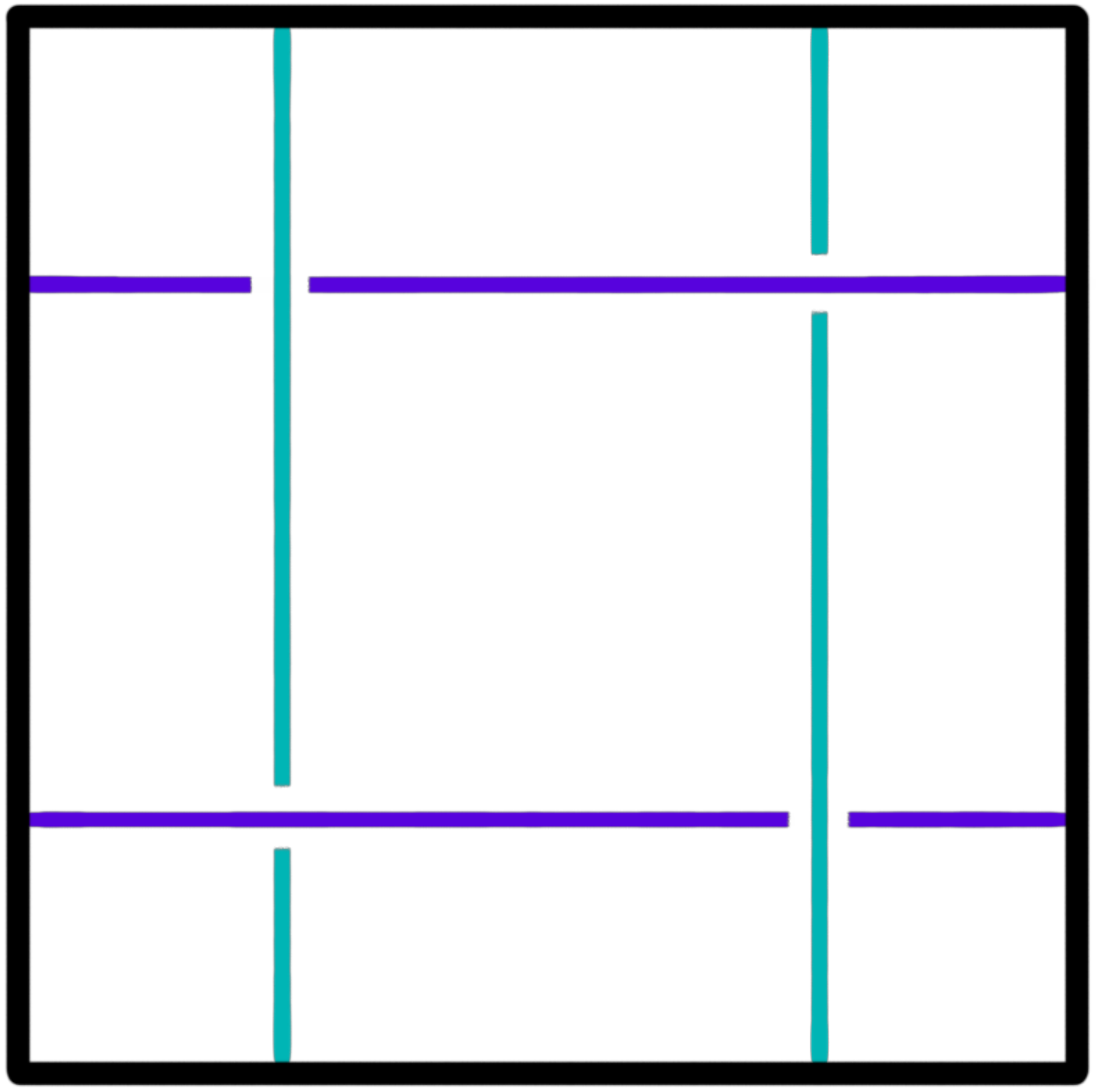}
        \vspace{0.275cm}
    \end{minipage}
    \caption{Diagrams of knotted structures: On the left, a diagram of the cinquefoil knot, obtained from a planar projection of an embedding of the knot. On the right, a motif diagram of a doubly periodic alternating weave, obtained from a projection of one of its motifs onto a 2-torus, represented as a square with identified edges.}
    \label{fig:diagrams_knots_2pw}
\end{figure}

Although knots or links are compact structures, periodic knottings or linkings can also be considered (Fukuda et al. \cite{fukuda2023construction}, Morton and Grishanov \cite{mortongrishanov2009}). In the case of doubly periodic structures, called \textit{2-structures} by Grishanov et al. in \cite{grishanovmeshkov2007} or \textit{DP tangles} by Diamantis et al. in \cite{diamantis2023equivalence}, these objects can be thought of as preimages in the universal cover $\mathbb{R}^2 \times [0,1]$ of links in the thickened torus $\mathbb{T}^2 \times [0,1]$.
As an infinite (doubly periodic) diagram of a DP tangle would contain infinitely many crossings, one may prefer to study diagrams of their motifs, that are compact domains generating the whole DP tangle under translations, since they have less complexity. These motif diagrams are defined as link diagrams on a 2-torus $\mathbb{T}^2$, as illustrated in Figure \ref{fig:diagrams_knots_2pw} right (see Fukuda et al. \cite{fukudakotanimahmoudi}, or Grishanov et al. \cite{grishanovmeshkov2009,grishanovmeshkov2009_2,grishanovmeshkov2007}). 
As there are multiple choices of motifs for a DP tangle, some complementary notions of equivalence have been added to account for these different motifs as detailed in \cite{diamantis2023equivalence}, such as {\it torus twists} \cite{grishanovmeshkov2007} and \textit{scale equivalence} (Mahmoudi \cite{mahmoudi2023tait}).

\begin{figure}[ht]
\centerline{\includegraphics[width=7.25cm]{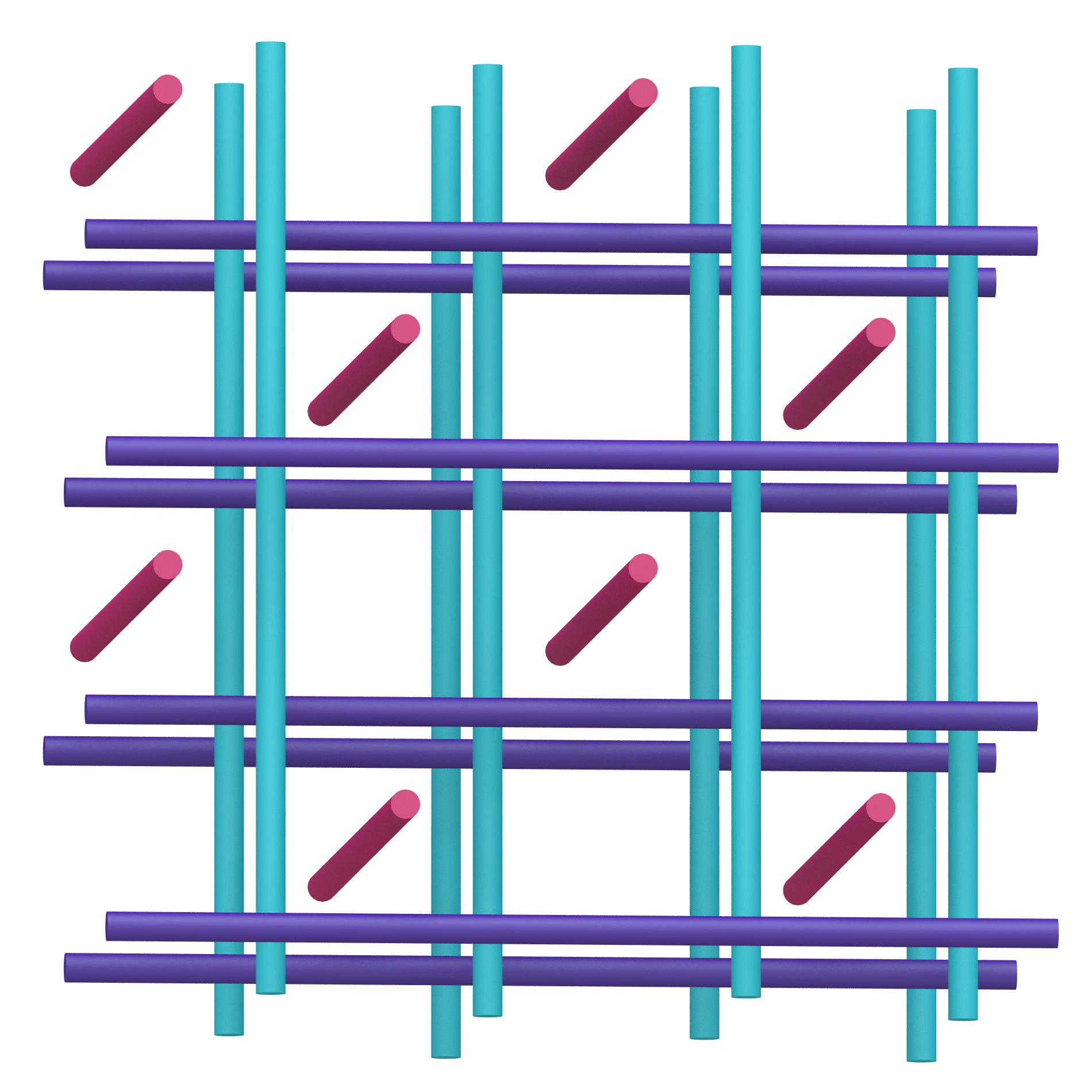}
\includegraphics[width=7.25cm]{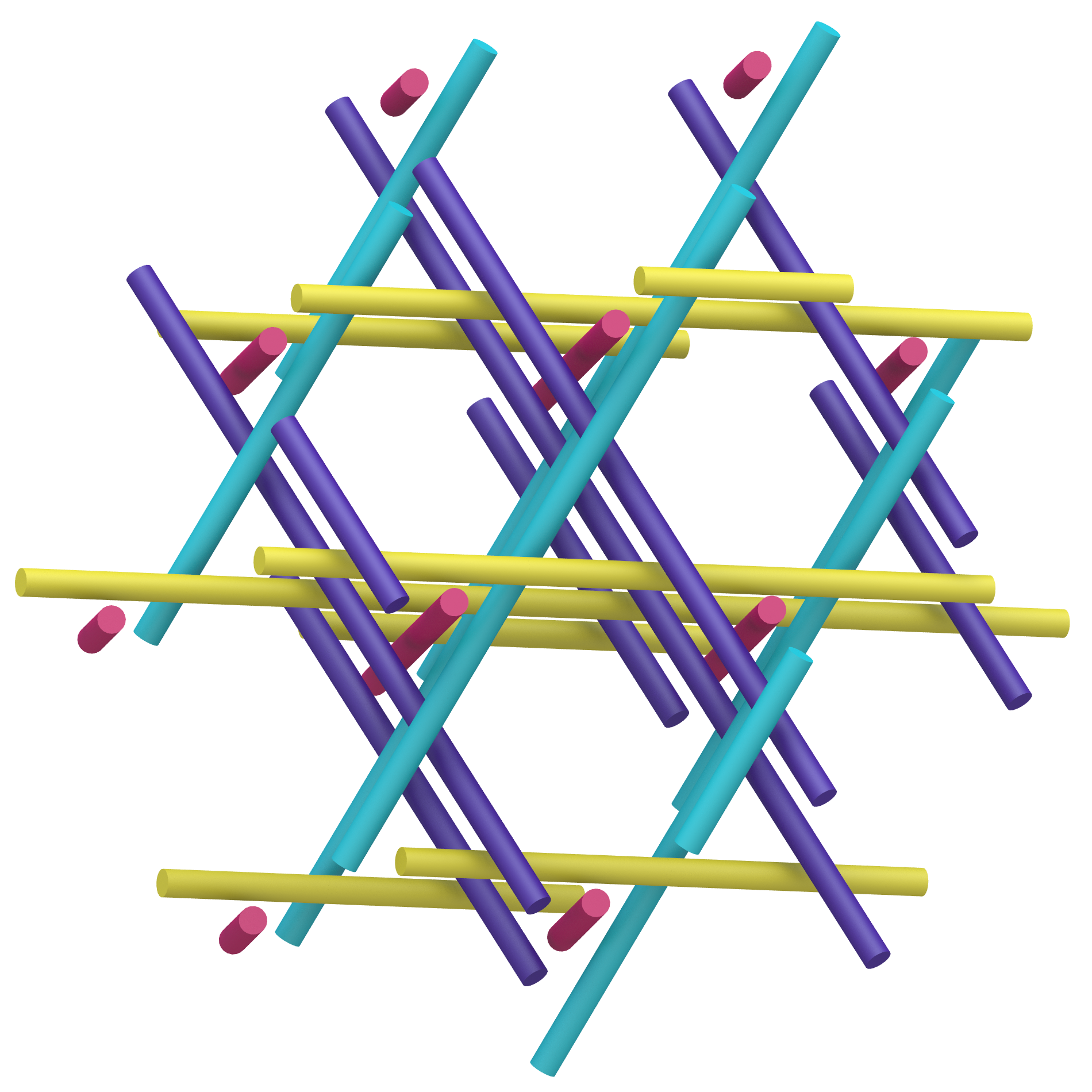}}
\vspace*{8pt}
\caption{Simple examples of 3-periodic tangles are given by periodic packings of straight cylinders. These two structures are well known in structural chemistry, referred to as the $\Pi^{+}$ cylinder packing on the left and the $\Sigma^{+}$ cylinder packing on the right \cite{okeeffe2001}. The cylinders are coloured by those that are parallel.}
\label{fig:TPT}
\end{figure}

There is interest in characterising entangled structures that are periodic in three directions of space (see, for example, Evans et al. \cite{Evans:eo5020}, and also Hui and Purcell \cite{hui2024geometry}), that we call \textit{3-periodic tangles} or {\it TP tangles}. TP tangles are found in many different systems in the natural sciences. For example, the $\Pi^{+}$ and $\Sigma^{+}$ cylinder packings (Figure \ref{fig:TPT}), which are packings of straight lines in $\mathbb{R}^3$ (see O'Keeffe et al. \cite{okeeffe2001}), are well known structures from structural chemistry which can be described as 3-periodic collections of curves embedded in $\mathbb{R}^3$. The $\Sigma^{+}$ structure has been used to model the mesoscale structure of mammalian skin cells (corneocytes) (Evans and Hyde \cite{evans2011}, Evans and Roth \cite{evans2014}), where the curves are helical rather than straight. Another setting where TP tangles are found is in DNA origami crystals (Seeman \cite{SEEMAN1982237}, Lu et al. \cite{luvecchioni2021}), where the DNA strands tangle around each other to form a complicated 3-periodic structure. In polymeric systems, topological features of the entangled polymer chains influence their physical properties (Likhtman and Ponmurugan \cite{likhtman2014}), and thus some knot invariants have been extended into measures for periodic systems, such as the \textit{periodic linking number} (Panagiotou \cite{PANAGIOTOU2015533}) or the \textit{Jones polynomial} (Barkataki and Panagiotou \cite{Barkataki_2024}) to characterise the entanglement, where the invariants encompass both geometry and topology of the filaments.

In this paper we describe a new diagrammatic representation of 3-periodic tangles, extending notions of equivalence to the 3-periodic setting. As two \textit{non-isotopic} links in $\mathbb{T}^3$ may possess the same 3-periodic preimage in the universal cover, studying isotopy in a 3-torus, such as in the paper of Lambropoulou and Rourke \cite{LAMBROPOULOU199795}, is insufficient for 3-periodic tangles. In Section \ref{sec:2} we present the definition of 3-periodic tangles and the equivalence relation between them. The latter is composed by three notions of equivalence, which are the isotopies in the 3-torus, as well as the \textit{torus twists} and the equivalence via covering maps relating non-isotopic links in the 3-torus that are quotients of the same TP tangle.  Section \ref{sec:3} states the definitions of diagrams and \textit{tridiagrams}, which are ordered 3-tuples of diagrams projected along three non-coplanar vectors, and their respective equivalences.  To the usual Reidemeister moves, we add six new moves that we call $R_4$, $R_5$, $R_6$, $R_7$, $R_8$ and $R_9$ moves. 
In particular, the $R_4$ and $R_5$ moves are similar to the $\Omega_4$ and $\Omega_5$ moves found in the paper of Mroczkowski and Dabkowski \cite{MROCZKOWSKI20091831} in the setting of links in $F\times \mathbb{S}^1$ where $F$ is an orientable surface. The $R_6$, $R_7$ and $R_8$ moves are similar to the moves of mixed link diagrams \cite{LAMBROPOULOU199795}, and the moves of DP tangle diagrams \cite{diamantis2023equivalence}. We end this Section by stating and proving a \textit{generalised Reidemeister theorem} for 3-periodic tangles using \textit{singularity theory}. Similar diagrams were hinted by Carrega in \cite{Carrega2017} and defined by Vuong independently from our work in \cite{vuong2023fundamental} with a different approach, limited to links in $\mathbb{T}^3$.
Finally, in Section \ref{sec:4} we demonstrate the practicality of the diagrams defined in the previous sections by defining the \textit{crossing number} of TP tangles that we illustrate with some examples.

\section{3-periodic tangles and their equivalence}\label{sec:2}
In this section, we start by introducing the notion of triply periodic entangled structures, which we refer to as \textit{3-periodic tangles} or \textit{TP tangles}. We then define the equivalence of 3-periodic tangles in $\mathbb{R}^3$ by investigating different equivalence relations of their quotients under lattices of periodicity, which can be regarded as links embedded in the 3-torus $\mathbb{T}^3 = \mathbb{S}^1 \times \mathbb{S}^1 \times \mathbb{S}^1$. We do this by considering the different choices of lattices and bases of those lattices, leading to equivalence relations between links in $\mathbb{T}^3$ via torus twists or covering maps. This serves as a basis to state and prove our generalised Reidemeister theorem for TP tangles in the next section.

\subsection{Definition of 3-periodic tangles} \label{subsec:2.1}

To define a 3-periodic tangle embedded in $\mathbb{R}^3$, we start by introducing different types of possible curve components. 

\begin{definition}
    We call a \textit{thread} an embedding of $\mathbb{R}$ in $\mathbb{R}^3$ having its two endpoints at infinity, and a \textit{loop} an embedding of $\mathbb{S}^1$ in $\mathbb{R}^3$.
\end{definition}

\begin{remark}
    Note that curves that do not have both of their endpoints at infinity, like finite strands such as an embedding given by the arctangent map, are not included in the definition of a thread. 
\end{remark}

Threads and loops form the components of the following object.

\begin{definition}
    An \textit{entanglement} is a disjoint union, possibly infinite, of threads and loops in $\mathbb{R}^3$.
\end{definition}

A 3-periodic tangle can be defined as a particular class of entanglement as follows.

\begin{definition}
    A \textit{3-periodic tangle} or \textit{TP tangle} $K$ is an entanglement that is invariant under some 3-dimensional point lattice $\Lambda$ of type $\mathbb{Z}.v_1 \oplus \mathbb{Z}.v_2 \oplus \mathbb{Z}.v_3 + v_4$, where the $v_i$'s are vectors of space, and $\lbrace v_1,v_2,v_3 \rbrace$ is a basis respecting the usual orientation of $\mathbb{R}^3$. When $K$ is equipped with such a lattice, we denote the structure by $\left\lbrace K,\Lambda \right\rbrace $.
\end{definition}

\begin{remark}
Strictly speaking, a lattice $\Lambda$ of a TP tangle, that we call lattice of periodicity for simplicity, is an affine collection of points. However, depending on the context, we will regard $\Lambda$ as either the affine object, or the $\mathbb{Z}$-module $\mathbb{Z}.v_1 \oplus \mathbb{Z}.v_2 \oplus \mathbb{Z}.v_3$, that is, the set of linear combinations of $v_1,v_2,v_3$ with scalars in $\mathbb{Z}$, to facilitate the proofs involving the lattices of a TP tangle.
\end{remark}

For 3-periodic tangles, we prefer working within their quotients under a lattice of periodicity, which are compact, and thus, have less complexity than the infinite structure.

\begin{definition}\label{def:unit_cell}
    Let $\left\lbrace K,\mathbb{Z}.v_1 \oplus \mathbb{Z}.v_2 \oplus \mathbb{Z}.v_3 + v_4 \right\rbrace $ be a 3-periodic tangle. By invariance under translation, $K$ maps to the quotient space $\mathbb{R}^3/ \left( \mathbb{Z}.v_1 \oplus \mathbb{Z}.v_2 \oplus \mathbb{Z}.v_3 + v_4 \right)$, represented by a parallelepiped delimited by the vectors $v_1,v_2,v_3$, where the opposite faces are identified. We call a \textit{unit cell} of $K$ such a quotient structure. It can be described as a union of disjointly embedded closed curves in the 3-torus, that is, a link embedded in $\mathbb{T}^3$. In particular, the vectors $v_1,v_2,v_3$ are identified with the generating circles of the $3$-torus that we denote by $l_1,l_2,l_3$.
\end{definition}

\begin{remark}\label{rmk:unit_cell_and_links_in_t3}
    In Definition \ref{def:unit_cell}, our understanding of the term `unit cell' is different from the usual meaning, where it is a fundamental region for the action of the three non-coplanar translations. We understand a unit cell as being a link embedded in the flat $3$-torus, represented by a parallelepiped with identified faces.
\end{remark}

\begin{remark}    
    Whenever the parallelepiped is not of interest, we regard a unit cell as a link in the $3$-torus that we embed in $\mathbb{C}^3$. We choose to embed $\mathbb{T}^3$ in $\mathbb{C}^3$, instead of any other manifold, to facilitate the description of any transformation that we apply to the links embedded in the $3$-torus. A consequence of this is that the two unit cells of Figure \ref{fig:pi_star_motifs_comparison} are exactly the same when regarded as links in the $3$-torus.

    \begin{figure}[hbtp]
        \centering
        \includegraphics[width=0.3\textwidth]{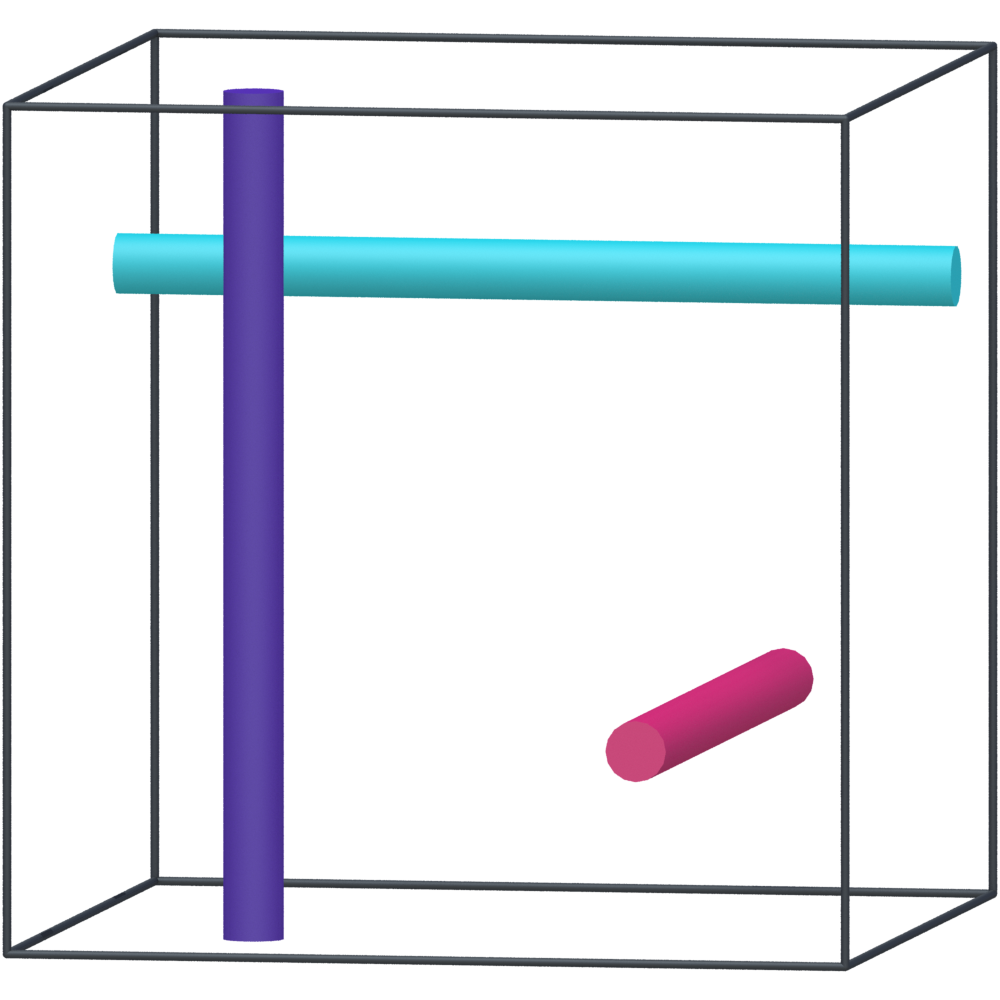}
        \hspace{0.5cm}
        \includegraphics[width=0.6125\textwidth]{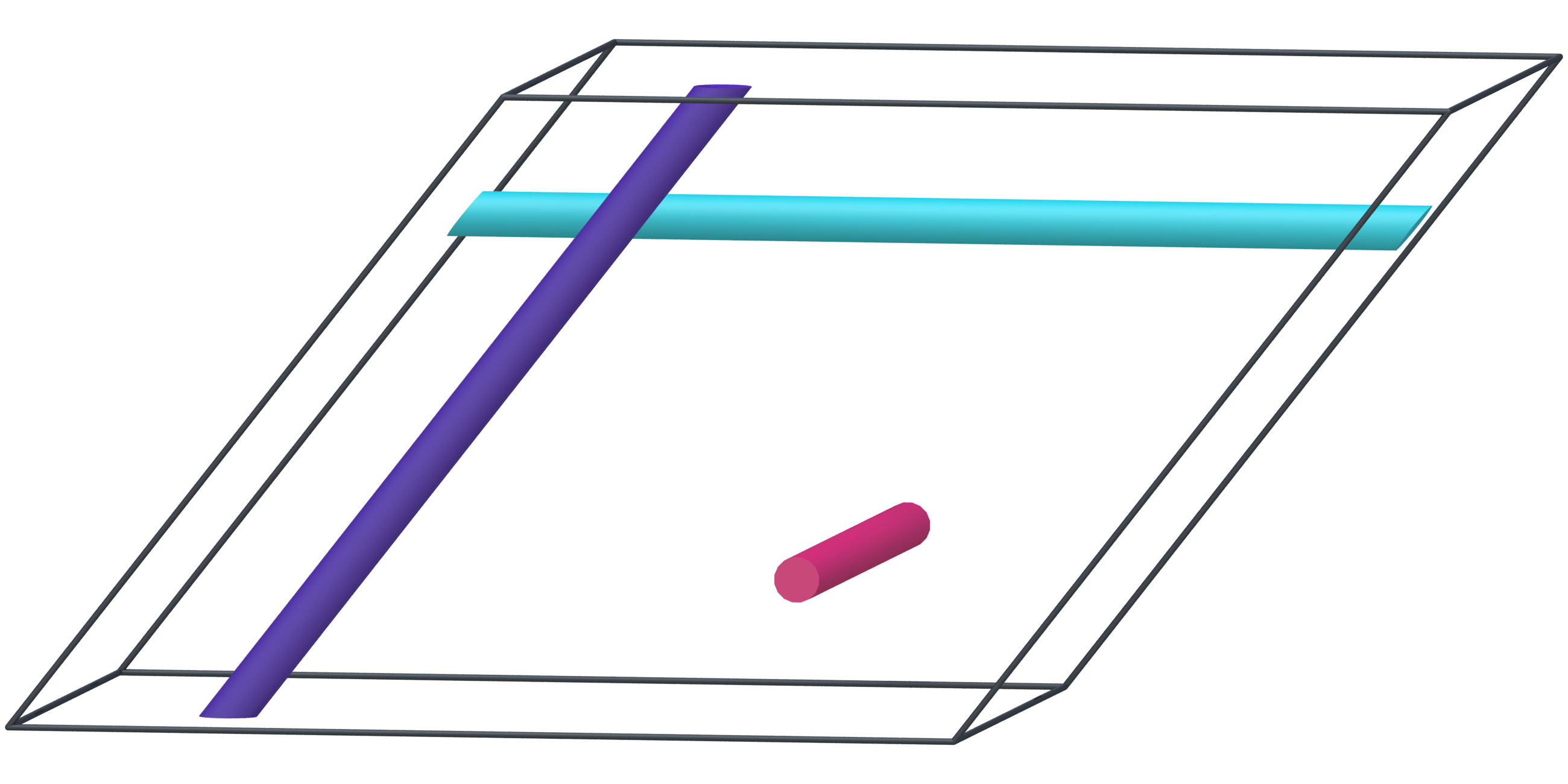}
        \caption{Two unit cells of two TP tangles whose parallelepiped differ by a shearing and a uniform scaling of space. They are associated to the same link embedded in the $3$-torus.}
        \label{fig:pi_star_motifs_comparison}
    \end{figure}

    Conversely, one can assign any parallelepiped to represent a link in the $3$-torus, by choosing any three vectors of space to represent the generating circles of the 3-torus.\\

    These considerations are important for the definition of different notions of equivalence that we will present later in the paper.
\end{remark}

\begin{remark}
    In a unit cell, the preimages in the universal cover $\mathbb{R}^3$ of a curve component homotopic to a point, are loops, and the preimages of all other curves are threads \cite{Farb_chap1}.
\end{remark}

\subsection{Equivalence of 3-periodic tangles}\label{subsec:equi_TP_tangles}

In general, one would like to tell whether two structures can continuously be deformed from one to the other. Equivalence of knotted objects is usually defined with the notion of ambient isotopy, and one can do the same for 3-periodic tangles. However, for a comprehensive diagrammatic description, we prefer giving another definition of equivalence, using unit cells of 3-periodic tangles, that we prove to be a particular case of ambient isotopies of TP tangles.\\

Some orientation-preserving affine transformations of the space $\mathbb{R}^3$, such as changes of basis, rotations, translations, uniform scalings, or shear deformations, preserve the isotopy class of TP tangles. These transformations modify the parallelepiped delimiting a unit cell. This motivates the following definition.

\begin{definition}\label{def:a-equi}
    Two 3-periodic tangles $K$ and $K'$ are \textit{$a$-equivalent} if they are connected by a finite sequence of orientation-preserving affine transformations such as orientation-preserving changes of basis of space, rotations, translations, uniform scalings or shear deformations. The equivalence is denoted by $K \sim_a K'$.
\end{definition}

\begin{figure}[ht]
\centering
\includegraphics[width= 10cm]{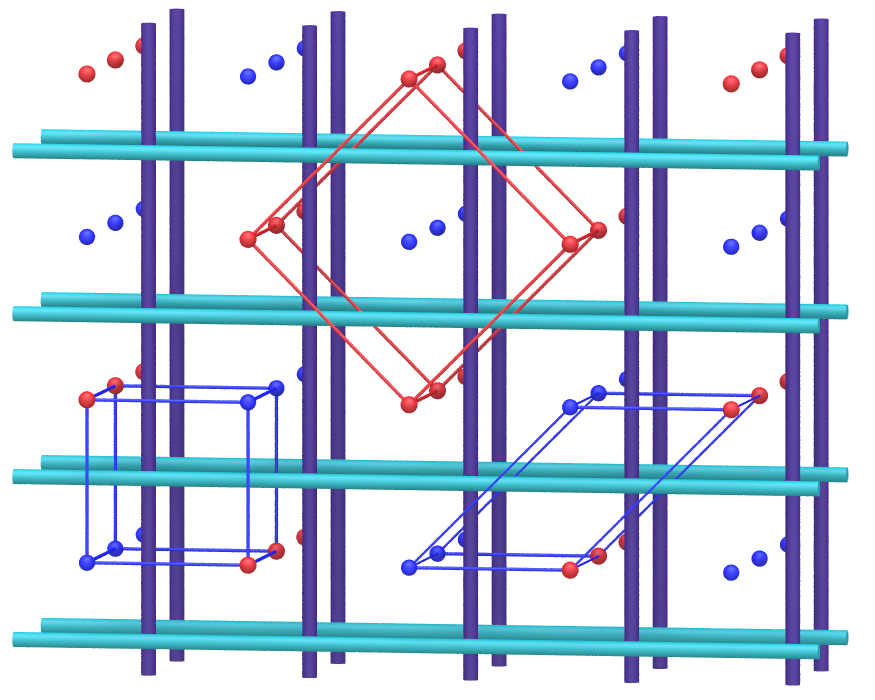}
\caption{Two different lattices (represented by red points and blue points, the red ones belonging to the blue lattice as well) that preserve the periodicity of the TP tangle. Different unit cells can be taken out of those lattices, like the blue cube and the blue parallelepiped both belonging to the blue lattice or the red parallelepiped which is a unit cell of volume 2 with respect to the blue lattice. Those three unit cells are non-isotopic links in the 3-torus and highlight the difference between such links and 3-periodic tangles.}
\label{fig:diff_unit_cells}
\end{figure}

For a given TP tangle, there are multiple choices of unit cells that represent it, as shown in Figure \ref{fig:diff_unit_cells}. Those unit cells are non-isotopic links in the 3-torus, yet they represent the same 3-periodic tangle, and thus, should be somehow related.\\

From a given lattice, there are many possible unit cells that can be obtained. Consider a 3-periodic tangle $K$ equipped with a lattice $\Lambda$. Assume, without loss of generality, that $\Lambda$ is the $\mathbb{Z}$-module $\mathbb{Z}^3$. Two different bases of $\mathbb{Z}^3$ give rise to two different unit cells $\Gamma$ and $\Gamma'$. Suppose, therefore, that $\left\lbrace v_1,v_2,v_3 \right\rbrace$ is a basis of $\mathbb{Z}^3$ that respects the usual orientation of $\mathbb{R}^3$. Suppose also that $M$ is the matrix whose columns are the coefficients of $\left\lbrace v_1,v_2,v_3 \right\rbrace$ in the canonical basis. It follows that $M$ is an invertible matrix in the set of integer matrices, with positive determinant, and therefore, belongs to $\mathrm{SL}(3,\mathbb{Z})$. Such a matrix induces an orientation-preserving linear homeomorphism of $\mathbb{R}^3$ that is equivariant with respect to the deck transformation group $\mathbb{Z}^3$, and thus descends to a homeomorphism $\phi_M$ of the torus $\mathbb{T}^3 = \mathbb{R}^3/\mathbb{Z}^3$ \cite{Farb_chap2}. This homeomorphism $\phi_M$ connects the links in the $3$-torus associated to the two unit cells $\Gamma$ and $\Gamma'$.

More precisely, if $f_M$ is the linear mapping associated to $M$, then we have
\[
\begin{array}{cccc}
    f_M : & \mathbb{R}^3 & \longrightarrow & \mathbb{R}^3 \\
     & (x,y,z) & \longmapsto & \left(f_{M,1}(x,y,z),f_{M,2}(x,y,z),f_{M,3}(x,y,z)\right)\\
     & \downarrow & & \downarrow \\
    \phi_M : & \mathbb{S}^1 \times \mathbb{S}^1 \times \mathbb{S}^1 & \longrightarrow & \mathbb{S}^1 \times \mathbb{S}^1 \times \mathbb{S}^1\\
    & \left(e^{2\mathrm{i}\pi x},e^{2\mathrm{i}\pi y},e^{2\mathrm{i}\pi z} \right) & \longmapsto & \left(e^{2\mathrm{i}\pi f_{M,1}(x,y,z)},e^{2\mathrm{i}\pi f_{M,2}(x,y,z)},e^{2\mathrm{i}\pi f_{M,3}(x,y,z)}\right) \\
\end{array}.
\]

Furthermore, it is possible to give a more accurate description of $\phi_M$. Indeed, recall that $\mathrm{SL}(3,\mathbb{Z})$ is generated by the following six shear matrices.
    \[
   M_1 = \left[\begin{array}{ccc}
   1  & 1 & 0  \\
   0 & 1 & 0 \\
   0 & 0 & 1
\end{array} \right], \quad
 M_2 = \left[\begin{array}{ccc}
   1  & 0 & 0  \\
   1 & 1 & 0 \\
   0 & 0 & 1
\end{array} \right], \quad
 M_3 = \left[\begin{array}{ccc}
   1  & 0 & 1  \\
   0 & 1 & 0 \\
   0 & 0 & 1
\end{array} \right],\]
\[
M_4 = \left[\begin{array}{ccc}
   1  & 0 & 0  \\
   0 & 1 & 0 \\
   1 & 0 & 1
\end{array} \right], \quad
M_5 = \left[\begin{array}{ccc}
   1  & 0 & 0  \\
   0 & 1 & 1 \\
   0 & 0 & 1
\end{array} \right], \quad
M_6 = \left[\begin{array}{ccc}
   1  & 0 & 0  \\
   0 & 1 & 0 \\
   0 & 1 & 1
\end{array} \right].
\]

We have the following property.

\begin{lemma} \label{lem:M_words_shear_matrices}
    If $L,N \in \mathrm{SL}(3,\mathbb{Z})$, then $\phi_{LN} = \phi_{L}\circ \phi_N$.
\end{lemma}

\begin{proof}
    If $f_{LN}$ is the linear mapping associated to $LN$, then 
    \[ f_{LN}(x,y,z) = LN \left[ \begin{array}{c}
         x\\
         y\\
         z
    \end{array} \right] = f_L \circ f_N (x,y,z).
    \]
    By descending to the quotient, the lemma follows.
\end{proof}

Lemma \ref{lem:M_words_shear_matrices} ensures that $\phi_M$ is a finite composition of the $\phi_{M_j}$'s, $j= 1,\dots,6$, and that one can restrict oneself to the homeomorphisms of the $3$-torus associated to the six shear matrices $M_j$'s listed above.

The action of one of the $\phi_{M_j}$'s on a link in the $3$-torus can actually be seen as `twisting' the $3$-torus along one of its generating circles. We give here the case of the matrix $M_1$.

\begin{example}
Consider the case of $M_1$, and the TP tangle associated to the unit cells displayed in Figure \ref{fig:torus_twists_M1}.

\begin{figure}[hbtp]
\centering
\includegraphics[width= 8cm]{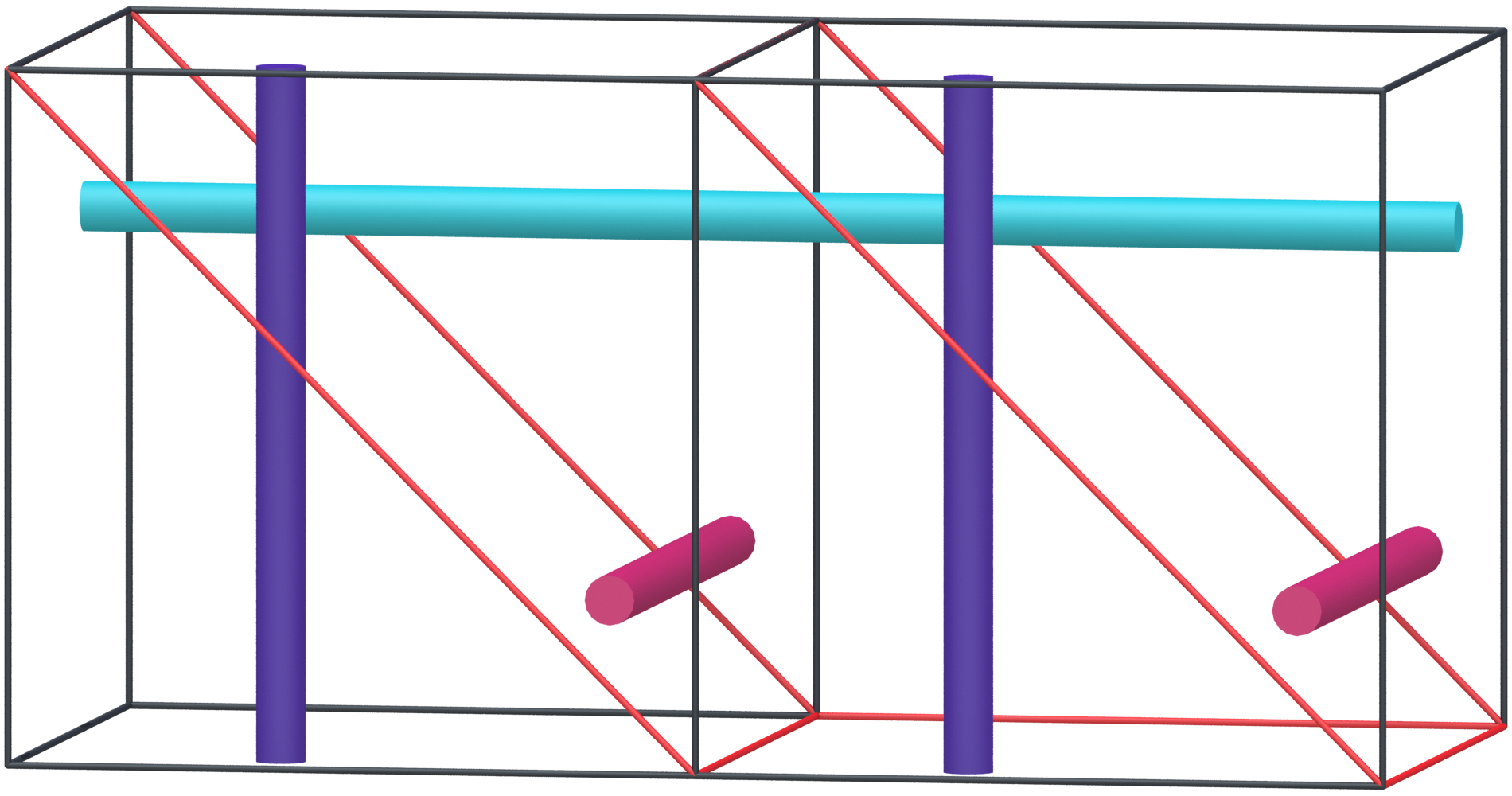}
\caption{Two unit cells that differ by the action of the matrix $M_1$.}
\label{fig:torus_twists_M1}
\end{figure}

The action of $M_1$ in $\mathbb{R}^3$ yields the unit cell delimited by the red parallelepiped of Figure \ref{fig:torus_twists_M1}. As links in the $3$-torus, this unit cell and the unit cell delimited by the unit cube are thus connected by $\phi_{M_1}$, whose formula states as follows:
\[
\begin{array}{cccc}
    \phi_{M_1} : & \mathbb{S}^1 \times \mathbb{S}^1 \times \mathbb{S}^1 & \longrightarrow & \mathbb{S}^1 \times \mathbb{S}^1 \times \mathbb{S}^1\\
    & \left(e^{2\mathrm{i}\pi x},e^{2\mathrm{i}\pi y},e^{2\mathrm{i}\pi z}\right) & \longmapsto & \left(e^{2\mathrm{i}\pi (x+y)},e^{2\mathrm{i}\pi y},e^{2\mathrm{i}\pi z}\right) \\
\end{array}.
\]

The action of $\phi_{M_1}$ can be interpreted as follows. Glue four faces of the cube delimiting $\mathbb{T}^3$, to get a thickened 2-torus with the inner and outer surfaces being identified as in Figure \ref{fig:thick_torus_twist_b}.
 The map $\phi_{M_1}$ acts on the thickened torus with identified surfaces by twisting it along the longitude as shown in Figure \ref{fig:thick_torus_twist_c}. By gluing four faces of the unit cell of Figure \ref{fig:thick_torus_twist_d}, one gets the twisted torus of Figure \ref{fig:thick_torus_twist_c}. This means that the links associated to the two unit cells of Figures \ref{fig:thick_torus_twist_a} and \ref{fig:thick_torus_twist_d} are connected by a twist of the torus. Furthermore, the unit cell of Figure \ref{fig:thick_torus_twist_d} and the red unit cell of Figure \ref{fig:torus_twists_M1} are actually associated to the same link in $\mathbb{T}^3$.
 
\begin{figure}[ht]
    \centering
    \begin{subfigure}[b]{4cm}
        \includegraphics[width=\textwidth]{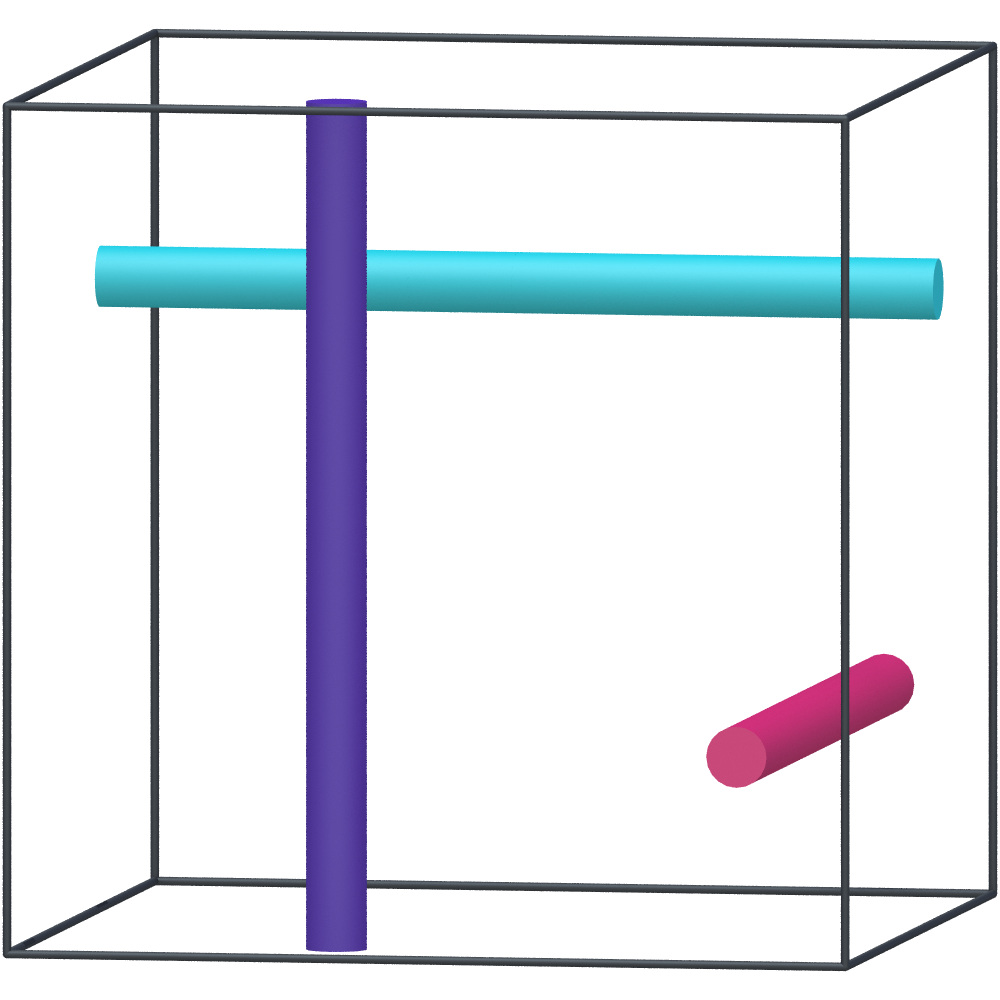}
        \caption{Unit cell of the TP tangle of Figure \ref{fig:torus_twists_M1} associated to the unit cube.}
        \label{fig:thick_torus_twist_a}
    \end{subfigure}
    \hspace{1cm}
    \begin{subfigure}[b]{5cm}
        \includegraphics[width=\textwidth]{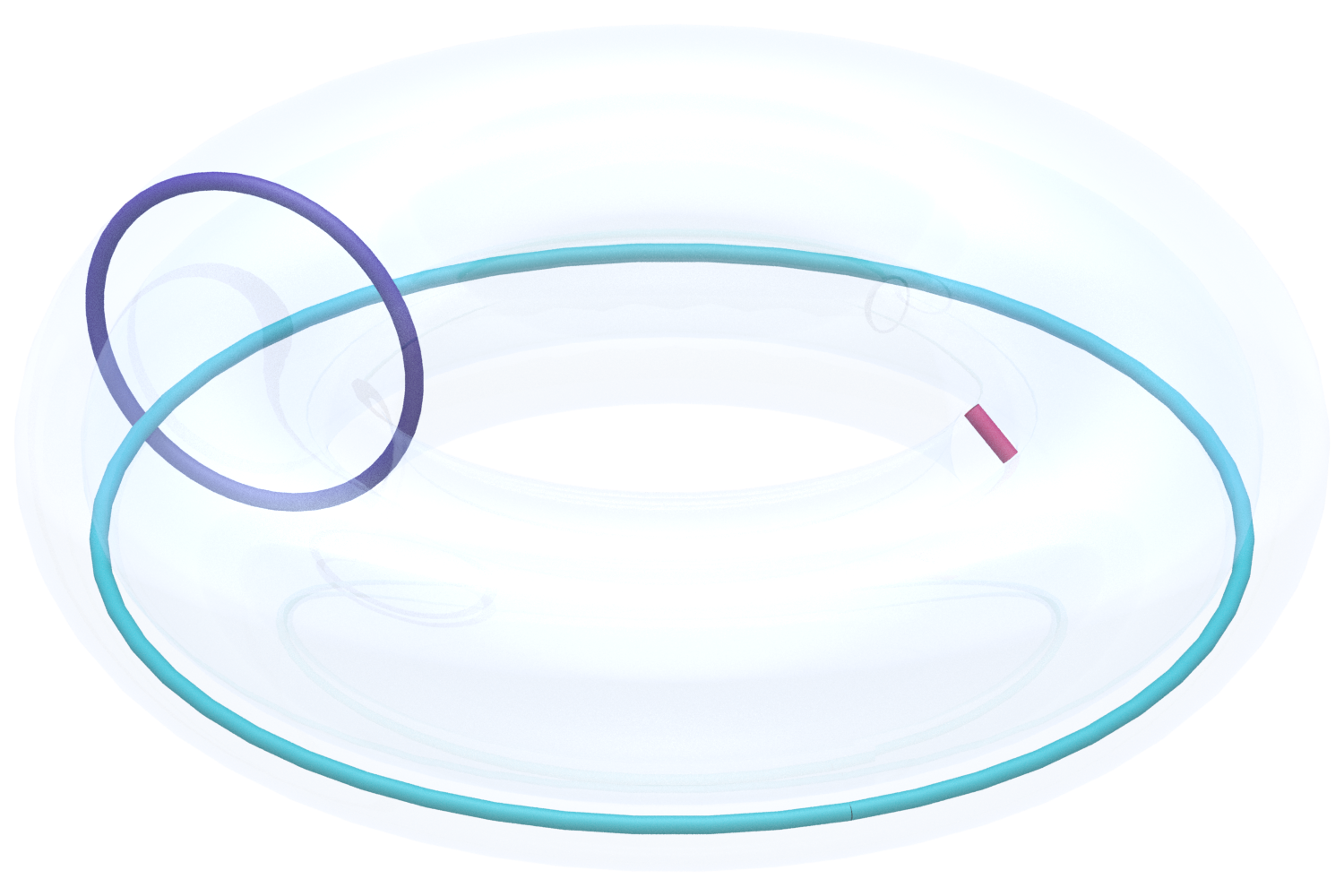}
        \caption{Same link as that of (a) where the $3$-torus is represented as a thickened 2-torus with identified inner and outer surfaces.}
        \label{fig:thick_torus_twist_b}
    \end{subfigure}
    \vskip\baselineskip
    \begin{subfigure}[b]{5cm}
        \includegraphics[width=\textwidth]{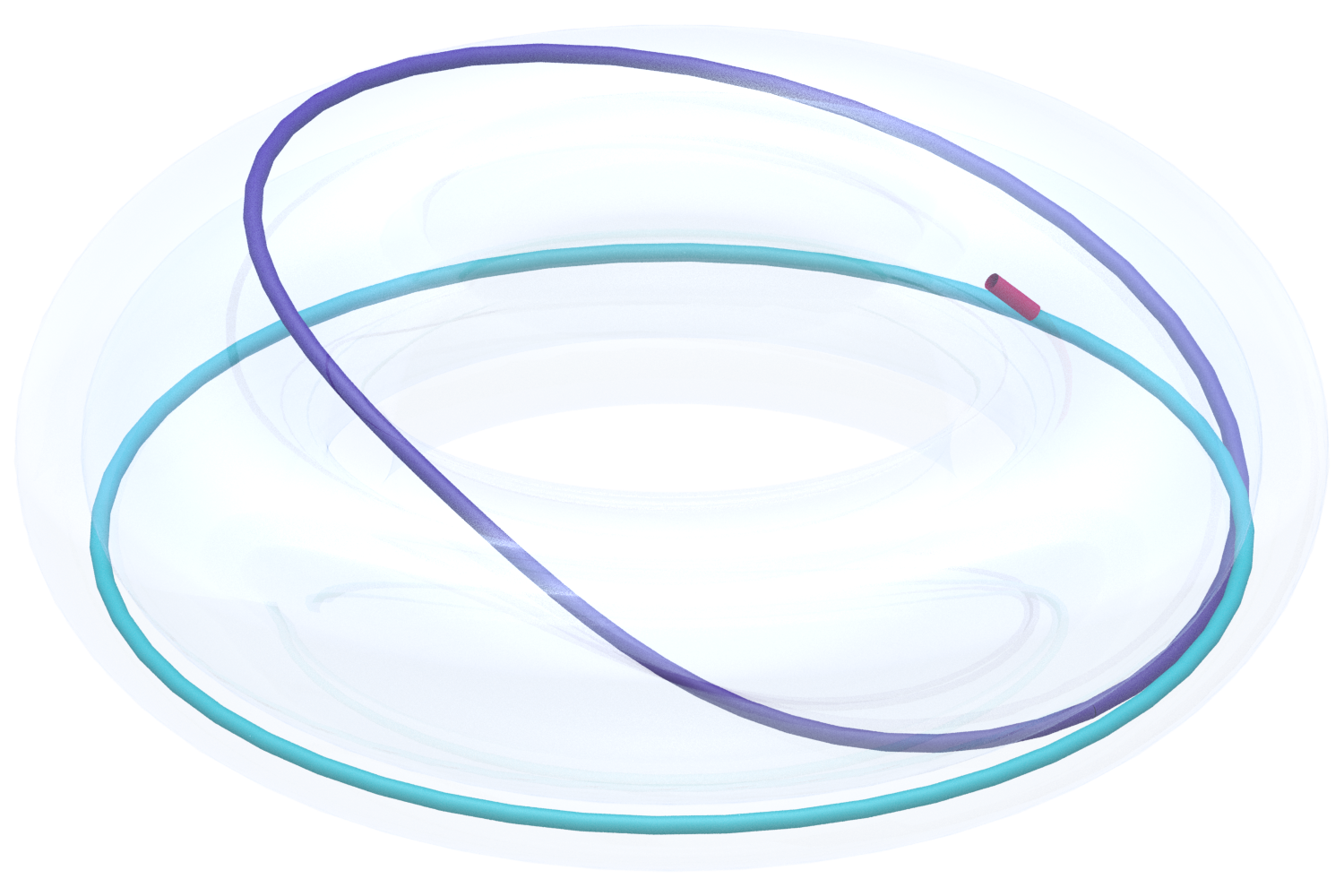}
        \caption{A link in the thickened 2-torus with identified surfaces connected to the link of (b) by the action of $\phi_{M_1}$.}
        \label{fig:thick_torus_twist_c}
    \end{subfigure}
    \hspace{1cm}
    \begin{subfigure}[b]{4cm}
        \includegraphics[width=\textwidth]{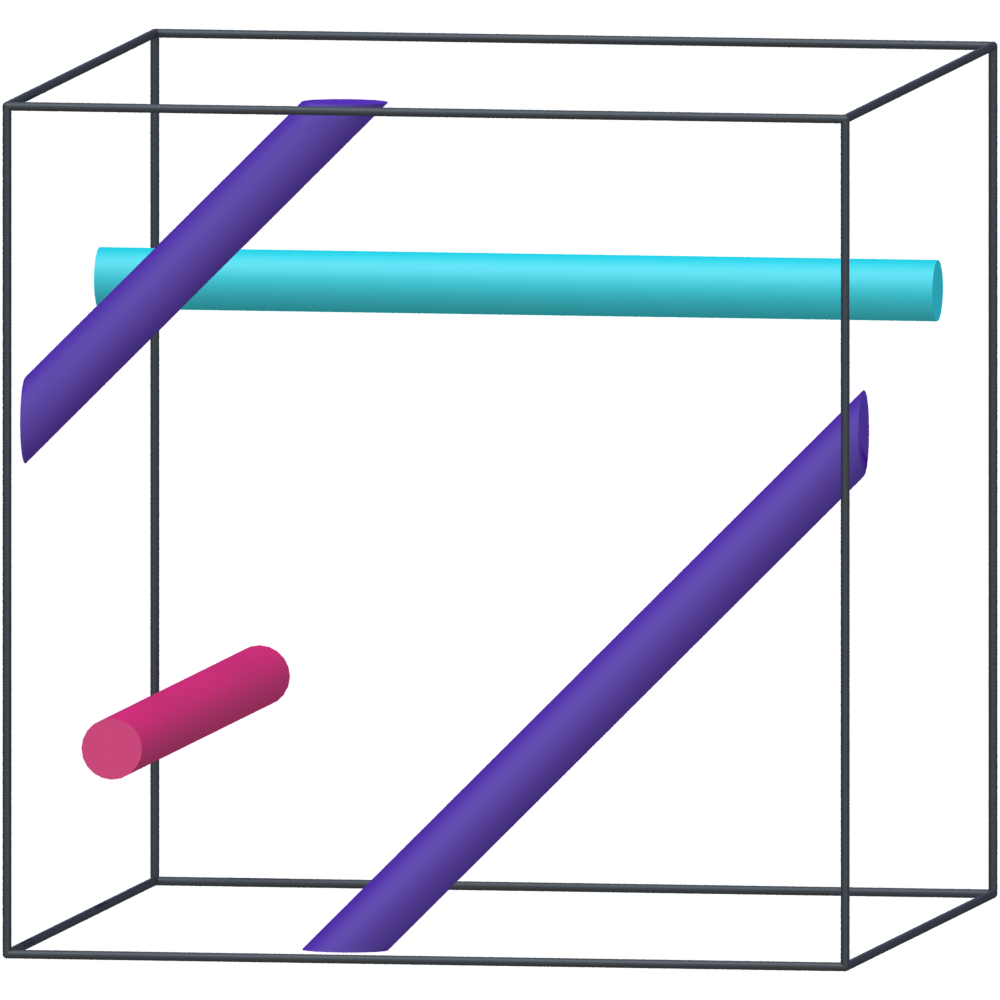}
        \caption{A unit cell whose associated link is the same as that of (c).}
        \label{fig:thick_torus_twist_d}
    \end{subfigure}
    \caption{Visualisation of the action of $\phi_{M_1}$ on a unit cell of a simple TP tangle: The link  in (c) is obtained from that of (b) by twisting the torus along the longitude.}
    \label{fig:thick_torus_twist}
\end{figure}
\end{example}

The actions of the other $\phi_{M_j}$'s are similar to that of $\phi_{M_1}$, with different gluing orders of the faces and different twists.

Using this as motivation, we give the following definition.
\begin{definition}\label{def:3_torus_twists}
    A \textit{$3$-torus twist}, or \textit{torus twist} when the context leaves no ambiguity, is one of the following six homeomorphisms:
    \[
\begin{array}{cccc}
    \phi_1 : & \mathbb{S}^1 \times \mathbb{S}^1 \times \mathbb{S}^1 & \longrightarrow & \mathbb{S}^1 \times \mathbb{S}^1 \times \mathbb{S}^1\\
    & \left(z_1,z_2,z_3\right) & \longmapsto & \left(z_1z_2,z_2,z_3\right) \\
\end{array},
\]
\[
\begin{array}{cccc}
    \phi_2 : & \mathbb{S}^1 \times \mathbb{S}^1 \times \mathbb{S}^1 & \longrightarrow & \mathbb{S}^1 \times \mathbb{S}^1 \times \mathbb{S}^1\\
    & \left(z_1,z_2,z_3\right) & \longmapsto & \left(z_1,z_1z_2,z_3\right) \\
\end{array},
\]
\[
\begin{array}{cccc}
    \phi_3 : & \mathbb{S}^1 \times \mathbb{S}^1 \times \mathbb{S}^1 & \longrightarrow & \mathbb{S}^1 \times \mathbb{S}^1 \times \mathbb{S}^1\\
    & \left(z_1,z_2,z_3\right) & \longmapsto & \left(z_1z_3,z_2,z_3\right) \\
\end{array},
\]
\[
\begin{array}{cccc}
    \phi_4 : & \mathbb{S}^1 \times \mathbb{S}^1 \times \mathbb{S}^1 & \longrightarrow & \mathbb{S}^1 \times \mathbb{S}^1 \times \mathbb{S}^1\\
    & \left(z_1,z_2,z_3\right) & \longmapsto & \left(z_1,z_2,z_1z_3\right) \\
\end{array},
\]
\[
\begin{array}{cccc}
    \phi_5 : & \mathbb{S}^1 \times \mathbb{S}^1 \times \mathbb{S}^1 & \longrightarrow & \mathbb{S}^1 \times \mathbb{S}^1 \times \mathbb{S}^1\\
    & \left(z_1,z_2,z_3\right) & \longmapsto & \left(z_1,z_2z_3,z_3\right) \\
\end{array},
\]
\[
\begin{array}{cccc}
    \phi_6 : & \mathbb{S}^1 \times \mathbb{S}^1 \times \mathbb{S}^1 & \longrightarrow & \mathbb{S}^1 \times \mathbb{S}^1 \times \mathbb{S}^1\\
    & \left(z_1,z_2,z_3\right) & \longmapsto & \left(z_1,z_2,z_2z_3\right) \\
\end{array}.
\]
\end{definition}

\begin{remark}
    The $\phi_j$'s of Definition \ref{def:3_torus_twists} correspond to the $\phi_{M_j}$'s of the shear matrices.
\end{remark}

We now give a notion of equivalence of $3$-torus links that are connected by torus twists.

\begin{definition}\label{def:3_torus_equivalence}
    Two links in the $3$-torus, $\Gamma_1$ and $\Gamma_2$, are \textit{$T^3$-equivalent} if they are connected by a finite composition of $3$-torus twists (and their inverses as well as the identity map).
\end{definition}

\begin{remark}
Definition \ref{def:3_torus_equivalence} generalises the notion of equivalence of unit cells of DP (Doubly Periodic) tangles under Dehn twists for a fixed lattice, called {\it Dehn equivalence} in \cite{diamantis2023equivalence}, first explained in \cite{grishanovmeshkov2007}.
\end{remark}

Given a TP tangle $K$ and a lattice associated to $K$, a translation of the points of the lattice generates another lattice of periodicity of $K$. Suppose that the points of the two lattices differ by a translation by $(x_0,y_0,z_0)$. The map

\[
\begin{array}{cccc}
     \mathbb{R}^3 & \longrightarrow & \mathbb{R}^3 \\
     (x,y,z) & \longmapsto & (x+x_0,y+y_0,z+z_0)\\
     \downarrow & & \downarrow \\
     \mathbb{S}^1 \times \mathbb{S}^1 \times \mathbb{S}^1 & \longrightarrow & \mathbb{S}^1 \times \mathbb{S}^1 \times \mathbb{S}^1\\
    \left(e^{2\mathrm{i}\pi x},e^{2 \mathrm{i} \pi y},e^{2 \mathrm{i} \pi z} \right) & \longmapsto & \left(e^{2 \mathrm{i}\pi (x+x_0)},e^{2 \mathrm{i} \pi (y+y_0)},e^{2 \mathrm{i} \pi (z+z_0)}\right) \\
\end{array}
\]
connects the unit cells of the two lattices.

This motivates the following definition.

\begin{definition}\label{def:shift_equi}
    Two links in the $3$-torus are \textit{$S$-equivalent} if they are connected by a map of type
    \[
\begin{array}{cccc}
     \mathbb{S}^1 \times \mathbb{S}^1 \times \mathbb{S}^1 & \longrightarrow & \mathbb{S}^1 \times \mathbb{S}^1 \times \mathbb{S}^1\\
    \left(z_1,z_2,z_3 \right) & \longmapsto & \left(z_1 e^{2 \mathrm{i}\pi x_0},z_2 e^{2 \mathrm{i} \pi y_0},z_3e^{2 \mathrm{i} \pi z_0}\right) \\
\end{array}, (x_0,y_0,z_0) \in \mathbb{R}^3.
\]
\end{definition}

There can be many lattices, whose associated $\mathbb{Z}$-modules are not the same, that preserve the periodicity of a TP tangle. We are going to show that some unit cells from those lattices are covers of one another. Several lemmas and definitions, that we present now, will be necessary to prove this statement.

\begin{lemma}\label{lem:sublattice}
    Consider two lattices $\Lambda$ and $\Lambda'$ of periodicity of a given TP tangle. Up to translation, one can assume that they have the same origin, and thus, can be regarded as $\mathbb{Z}$-modules. If $\mathcal{B} = \lbrace u_1,u_2,u_3 \rbrace$ and $ \mathcal{B}' = \lbrace u^{\prime}_1,u^{\prime}_2,u^{\prime}_3 \rbrace$ are two respective bases of the lattices, then the $\mathbb{Z}$-module $\Lambda_0$ generated by $\mathcal{B} \cup \mathcal{B}'$ also preserves the periodicity of the TP tangle.
\end{lemma}

\begin{proof}
    Consider the element $p_0 = \lambda_1u_1 + \lambda_2u_2 + \lambda_3u_3 + \lambda_1'u_1' + \lambda_2'u_2' + \lambda_3'u_3'$ of $\Lambda_0$. The point $\lambda_1'u_1' + \lambda_2'u_2' + \lambda_3'u_3'$ is an element of $\Lambda'$. Regarded as an affine collection of points, one can change the origin of $\Lambda'$ by doing a translation by $ - \left(\lambda_1'u_1' + \lambda_2'u_2' + \lambda_3'u_3'\right)$. With the new basis, the point $p_0$ becomes $\lambda_1u_1 + \lambda_2u_2 + \lambda_3u_3$. This latter is an element of $\Lambda$ which also preserves the periodicity of the TP tangle. This means that the point $\lambda_1u_1 + \lambda_2u_2 + \lambda_3u_3 + \lambda_1'u_1' + \lambda_2'u_2' + \lambda_3'u_3'$ is also a point of a lattice of periodicity of the TP tangle. This implies that $\Lambda_0$ is a lattice of periodicity of the TP tangle.
\end{proof}

\begin{definition}
    Consider a 3-periodic tangle $K$. Regard its lattices as $\mathbb{Z}$-modules. For a given lattice $\Lambda_1$ associated to $K$, a lattice $\Lambda_2$ is called a \textit{refinement} of $\Lambda_1$ if $\Lambda_1$ is a subgroup of $\Lambda_2$, and if $\Lambda_2$ also preserves the periodicity of $K$.
\end{definition}

    Lemma \ref{lem:sublattice} ensures that, in general, one can restrict oneself to studying submodules of a lattice whenever one considers different lattices. Therefore, one wishes to have a description of submodules of a lattice. To obtain such a description, we recall the following definitions and properties regarding submodules of an integer lattice.

\begin{proposition}\label{prop:hermite_normal_form}
    For any $A \in \mathrm{GL}(n,\mathbb{R}) \cap \mathrm{M}_n(\mathbb{Z})$, where $\mathrm{M}_n(\mathbb{Z})$ denotes the set of $n\times n$ matrices with integer coefficients, there is a unique $H \in \mathrm{GL}(n,\mathbb{R}) \cap \mathrm{M}_n(\mathbb{Z})$ such that $H = AU$ for some $U \in \mathrm{GL}(n,\mathbb{Z})$, and $H = \left[h_{ij}\right]$ satisfies the following conditions:
    \begin{itemize}
        \item $h_{ij} = 0$ for $i < j$, that is, $H$ is lower triangular,
        \item $0 \leqslant h_{ij} < h_{ii}$ for $j < i$.
    \end{itemize}
\end{proposition}

The matrix $H$ is called the \textit{Hermite normal form} of $A$. This is a generalisation over $\mathbb{Z}$ of the reduced echelon form of a matrix. A proof of Proposition \ref{prop:hermite_normal_form} is given in \cite{Cohen1993_chap2}. It allows us to characterise a basis of a submodule of a lattice.

\begin{lemma}\label{lem:submodule_basis}
    Consider a submodule $\Lambda$ of $\mathbb{Z}^3$ that is itself isomorphic to $\mathbb{Z}^3$ (full-dimensional). There exists a basis of $\Lambda$ of type $\left\lbrace h_{11}e_1 + h_{21}e_2 + h_{31}e_3,\right.$ $\left. h_{22}e_2 +h_{32}e_3,\right.$ $\left. h_{33}e_3 \right\rbrace$, where the $e_i$'s form the canonical basis of $\mathbb{Z}^3$, and $0 \leqslant h_{ij} < h_{ii}$ for $j < i$.
\end{lemma}

\begin{proof}
    As $\Lambda$ is a full-dimensional submodule of $\mathbb{Z}^3$, it possesses a basis of type $\left\lbrace a_1,a_2,a_3 \right\rbrace$ that respects the usual orientation of $\mathbb{R}^3$. Denote by $A$ the matrix whose columns are the coefficients of $\left\lbrace a_1,a_2,a_3 \right\rbrace$ in the canonical basis, that is, with abuse of language, $[a_1,a_2,a_3] = [e_1,e_2,e_3]A$.
    It follows that the matrix $A$ belongs to $\mathrm{GL}(3,\mathbb{R}) \cap \mathrm{M}_3(\mathbb{Z})$.

    Proposition \ref{prop:hermite_normal_form} states that there is a unique $H$ whose coefficients satisfy the conditions of the lemma, and some $U$ such that $H = AU$. 
 Notice in particular that, since $H$ and $A$ are of positive determinant, $U$ belongs to $\mathrm{SL}(3,\mathbb{Z})$. The equality $H = AU$ implies that $[e_1,e_2,e_3]AU = [e_1,e_2,e_3]H$ constitutes a basis of $\Lambda$.
\end{proof}

Suppose now that for a TP tangle $K$, we have $\Lambda_2$ a refinement of $\Lambda_1$. Without loss of generality, assume that $\Lambda_2$ is $\mathbb{Z}^3$. Since $\Lambda_1$ is a full dimensional submodule of $\Lambda_2$, Lemma \ref{lem:submodule_basis} states that $\Lambda_1$ is of type $\mathbb{Z}. (h_{11}e_1 + h_{21}e_2 + h_{31}e_3) \oplus \mathbb{Z}. (h_{22}e_2 +h_{32}e_3) \oplus \mathbb{Z}. h_{33}e_3$. As $\mathbb{Z}^3$ is the fundamental group of $\mathbb{T}^3$, and as $\mathbb{R}^3$ is simply connected, by the connection between subgroups of the fundamental group and path-connected covering spaces of $\mathbb{R}^3/\mathbb{Z}^3$, called Galois connection \cite{pierre_guillot}, there is a covering map
\[
\begin{array}{ccccc}
   \mathbb{R}^3/\Lambda_1 \cong & \mathbb{S}^1 \times \mathbb{S}^1 \times \mathbb{S}^1 & \longrightarrow & \mathbb{S}^1 \times \mathbb{S}^1 \times \mathbb{S}^1 & \cong \mathbb{R}^3/\mathbb{Z}^3 \\
   &   (z_1,z_2,z_3) & \longmapsto & \left({z_1}^{h_{11}}{z_2}^{h_{21}}{z_3}^{h_{31}},{z_2}^{h_{22}}{z_3}^{h_{32}},{z_3}^{h_{33}}\right) &
\end{array},
\]
where $z_i \in \mathbb{C}, \|z_i\| = 1$.

\noindent This is a covering map with $h_{11}h_{22}h_{33}$ sheets for each point of $\mathbb{T}^3$.

\begin{example} \label{eg:multi_unit_cell}
Suppose $\Lambda_1 = \mathbb{Z}. 2e_1 \oplus \mathbb{Z} .e_2 \oplus \mathbb{Z}. e_3$. Consider the TP tangle having a unit cell displayed in Figure \ref{fig:thick_torus_twist_a}. In this case, the covering map is
\[ \begin{array}{ccccc}
   \mathbb{R}^3/ (\mathbb{Z}. 2e_1 \oplus \mathbb{Z} .e_2 \oplus \mathbb{Z}. e_3) \cong & \mathbb{S}^1 \times \mathbb{S}^1 \times \mathbb{S}^1 & \longrightarrow & \mathbb{S}^1 \times \mathbb{S}^1 \times \mathbb{S}^1 & \cong \mathbb{R}^3/\mathbb{Z}^3  \\
     & (z_1,z_2,z_3) & \longmapsto & \left({z_1}^2,z_2,z_3\right)
\end{array},\]
and the unit cell corresponding to the lattice $\mathbb{Z}. 2e_1 \oplus \mathbb{Z} .e_2 \oplus \mathbb{Z}. e_3$ is the double of the one of Figure \ref{fig:thick_torus_twist_a} along the $x$-axis as seen in Figure \ref{fig:pi_star_double_unit}.

\begin{figure}[hbtp]
\centering
\includegraphics[width= 8cm]{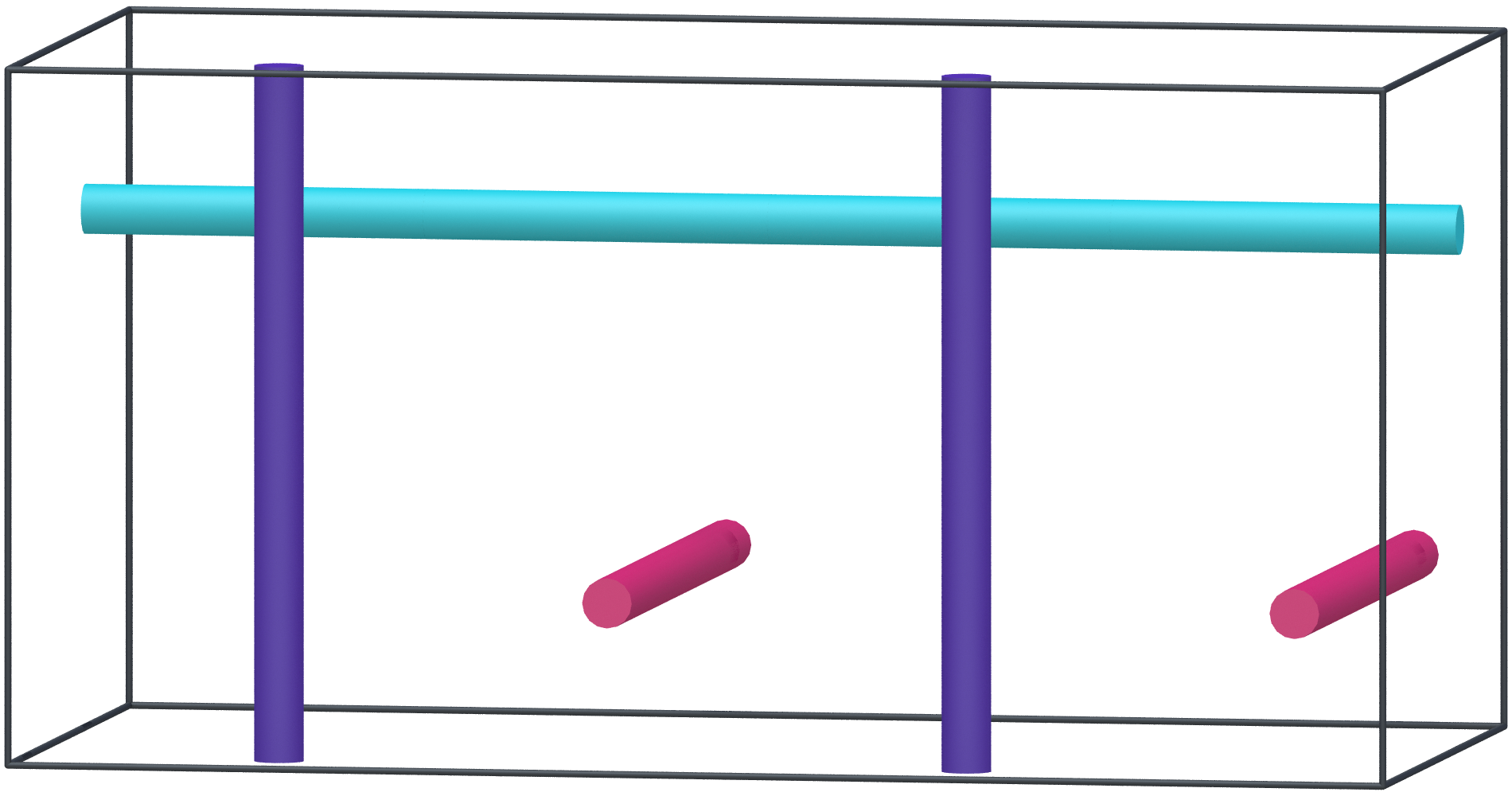}
\caption{A unit cell of the simple packing of cylinders of Figure \ref{fig:thick_torus_twist} corresponding to the lattice $\mathbb{Z}. 2e_1 \oplus \mathbb{Z} .e_2 \oplus \mathbb{Z} .e_3$, that is of volume twice that of the unit cell in Figure \ref{fig:thick_torus_twist_a}.}
\label{fig:pi_star_double_unit}
\end{figure}
\end{example}

In the following, we state a notion of equivalence for links in the $3$-torus that are connected by covering maps.

\begin{definition}\label{def:cover_3}
    Two links in the $3$-torus, $\Gamma$ and $\Gamma'$, are \textit{$CM^3$-equivalent} if there is a sequence of links in the $3$-torus $(\Gamma_k)_{k=0,\dots,m}$ where $\Gamma_0 = \Gamma$, $\Gamma_m = \Gamma'$, and for all $k=1,\dots,m$, either $\Gamma_{k-1}$ is a cover $\Gamma_{k}$ or $\Gamma_{k}$ is a cover $\Gamma_{k-1}$, where the covering map is of type 
    \[
\begin{array}{ccc}
    \mathbb{S}^1 \times \mathbb{S}^1 \times \mathbb{S}^1 & \longrightarrow & \mathbb{S}^1 \times \mathbb{S}^1 \times \mathbb{S}^1 \\
      (z_1,z_2,z_3) & \longmapsto & \left({z_1}^{h_{11}}{z_2}^{h_{21}}{z_3}^{h_{31}},{z_2}^{h_{22}}{z_3}^{h_{32}},{z_3}^{h_{33}}\right)
\end{array}, z_i \in \mathbb{C}, \|z_i\| = 1,
\]
with $0 \leqslant h_{ij} < h_{ii}$ for $j < i$.
\end{definition}

\begin{proposition}\label{prop:unit_cells_are_equi}
    Any two unit cells $\Gamma_1$ and $\Gamma_2$ of a given embedding of a TP tangle $K$, regarded as links in the $3$-torus, are connected by a finite sequence of links in $\mathbb{T}^3$ that are $T^3$-equivalent, $S$-equivalent or $CM^3$-equivalent.
\end{proposition}

\begin{proof}
    Suppose that $\Gamma_1$ and $\Gamma_2$ are respectively associated to two lattices $\Lambda_1$ and $\Lambda_2$, where $\Lambda_1$ and $\Lambda_2$ are not necessarily different. Up to a translation, and thus an $S$-equivalence of the links in the $3$-torus, one can assume that $(0,0,0)$ belongs to $\Lambda_1$ and $\Lambda_2$. Lemma \ref{lem:sublattice} ensures the existence of a refinement $\Lambda_0$,  common to both $\Lambda_1$ and $\Lambda_2$. This implies the existence of three unit cells $\Gamma_0'$, $\Gamma_1'$ and $\Gamma_2'$, respectively associated to $\Lambda_0$, $\Lambda_1$ and $\Lambda_2$, that are $CM^3$-equivalent. Immediately, we have $\Gamma_1$ and $\Gamma_2$ respectively $T^3$-equivalent to $\Gamma_1'$ and $\Gamma_2'$, as they differ only by the different bases delimiting the parallelepipeds.
\end{proof}

We are now in a position to state a definition of equivalence of TP tangles.

\begin{definition} \label{def:u_equivalence}
Let $\lbrace K,\Lambda \rbrace$ and $\lbrace K',\Lambda' \rbrace$ be two TP tangles with given lattices of periodicity $\Lambda$ and $\Lambda'$. They are \textit{$U$-equivalent} if there are two TP tangles $K_a$ and $K_a'$ with given lattices of periodicity $\Lambda_a$ and $\Lambda_a'$ that satisfy the following conditions:
\begin{itemize}
    \item[-] $K_a \sim_a K$ and $K_a' \sim_a K'$,
    \item[-] there is a finite sequence of links in $\mathbb{T}^3$, $(\Gamma_k)_{k=0,\dots,n}$ such that:
    \begin{itemize}
        \item[.] $\Gamma_0$ is a unit cell of $\lbrace K_a,\Lambda_a \rbrace$,
        \item[.] $\Gamma_n$ is a unit cell of $\lbrace K_a',\Lambda_a' \rbrace$,
        \item[.] $\forall k=1,\dots,n$, $\Gamma_{k-1}$ and $\Gamma_{k}$ are either ambient isotopic, or $T^3$-equivalent as given by Definition \ref{def:3_torus_equivalence}, or $S$-equivalent as per Definition \ref{def:shift_equi}, or $CM^3$-equivalent as given by Definition \ref{def:cover_3}.
    \end{itemize}
\end{itemize}
\end{definition}

\begin{proposition}\label{prop:u_equi_is_equi_relation}
    The $U$-equivalence of Definition \ref{def:u_equivalence} is an equivalence relation.
\end{proposition}

\begin{proof}
    Reflexivity and symmetry are immediate.
    
    For transitivity, suppose that the TP tangles $\lbrace K,\Lambda \rbrace$ and $\left\lbrace K^{(1)},\Lambda^{(1)} \right\rbrace$, as well as $\left\lbrace K^{(1)},\Lambda^{(1)} \right\rbrace$ and $\left\lbrace K^{(2)},\Lambda^{(2)} \right\rbrace$ are $U$-equivalent. We have:
    \begin{itemize}
    \item[-] $K_a \sim_a K$ and $K_{a,1}^{(1)} \sim_a K^{(1)}$,
    \item[-] a sequence of links in $\mathbb{T}^3$, $(\Gamma_{0 \rightarrow 1, k})_{k=0,\dots,n_1}$ connecting $K_a$ and $K_{a,1}^{(1)}$, and satisfying the conditions of the definition of $U$-equivalence,
    \item[-] $K_{a,2}^{(1)} \sim_a K^{(1)}$ and $K_{a}^{(2)} \sim_a K^{(2)}$,
    \item[-] a sequence of links in $\mathbb{T}^3$, $(\Gamma_{1 \rightarrow 2, k})_{k=0,\dots,n_2}$ connecting $K_{a,2}^{(1)}$ and $K_{a}^{(2)}$, and satisfying the conditions of the definition of $U$-equivalence.
\end{itemize}
     Denote by $\varphi_1$, $\varphi_2$ and $\varphi_3$ the composite maps of the affine transformations that respectively connect $K_{a,1}^{(1)}$ to $K^{(1)}$,  $K^{(1)}$ to $K_{a,2}^{(1)}$, and $K_{a}^{(2)}$ to $ K^{(2)}$. We have $\varphi_1^{-1} \circ \varphi_2^{-1}\left( K_{a,2}^{(1)} \right) = K_{a,1}^{(1)}$. This fact and Remark \ref{rmk:unit_cell_and_links_in_t3}, imply that $\Gamma_{1 \rightarrow 2, 0}$ can be regarded as a link that is associated to a unit cell of $K_{a,1}^{(1)}$. By Proposition \ref{prop:unit_cells_are_equi}, $\Gamma_{0 \rightarrow 1, n_1}$ and $\Gamma_{1 \rightarrow 2, 0}$ are connected by a finite sequence of links in $\mathbb{T}^3$, that we denote by $(\Gamma_{(1,1) \rightarrow (1,2), k})_{k=0,\dots,m}$, that are $T^3$-equivalent, $S$-equivalent or $CM^3$-equivalent. We obtain a sequence
     \[\left(\Gamma_{0 \rightarrow 1, 0}, \dots, \Gamma_{0 \rightarrow 1, n_1} =\Gamma_{(1,1) \rightarrow (1,2), 0},\dots, \Gamma_{(1,1) \rightarrow (1,2), m} = \Gamma_{1 \rightarrow 2, 0}, \dots, \Gamma_{1 \rightarrow 2, n_2} \right)\]
     of links in $\mathbb{T}^3$ that are ambient isotopic, $T^3$-equivalent, $S$-equivalent or $CM^3$-equivalent. This sequence connects $K_{a,1}^{(1)}$ to a TP tangle $K_{a}^{(3)}$ satisfying $\varphi_3 \circ \varphi_2 \circ \varphi_1 \left(K_{a}^{(3)} \right) = K^{(2)}$, which completes the proof.
\end{proof}

\begin{remark}
     In Definition \ref{def:u_equivalence} on $U$-equivalence, the $a$-equivalence, encompassing geometric transformations, is considered only at the beginning and at the end of the transformations that connect two TP tangles. One may wonder what would happen if some affine transformation $\zeta$ were applied at some intermediary stage. For example, if one had $K$ transformed into $K'$ in the following way: First, we have $K_a= \xi_1(K)$, with $\xi_1$ an affine transformation. Then, there is the sequence $(\Gamma_{0\rightarrow 1,k})_{k=0,\dots,n_1}$ of links in the $3$-torus, where $\Gamma_{0\rightarrow 1,0}$ is a unit cell of $K_a$. This is followed by $K_{a,2}= \zeta(K_{a,1})$, where $K_{a,1}$ is the TP tangle associated to $\Gamma_{0\rightarrow 1,n_1}$. Then, there is another sequence of links in the $3$-torus $(\Gamma_{1\rightarrow 2,k})_{k=0,\dots,n_2}$, with $\Gamma_{1\rightarrow 2,0}$ a unit cell of $K_{a,2}$. Finally, $K' = \xi_2(K_{a,3})$, where $K_{a,3}$ is the TP tangle associated to $\Gamma_{1\rightarrow 2,n_2}$, and $\xi_2$ an affine transformation.     
     The proof of Proposition \ref{prop:u_equi_is_equi_relation} on transitivity actually shows that it is possible to not apply $\zeta$ on $K_{a,1}$, but only after the second sequence of links in the $3$-torus $(\Gamma_{1\rightarrow 2,k})_{k=0,\dots,n_2}$. Indeed, $\zeta$ can be regarded as $\varphi_2 \circ \varphi_1$ of the proof of Proposition \ref{prop:u_equi_is_equi_relation}. By doing so, we come back to the definition of $U$-equivalence.
\end{remark}

\begin{proposition}\label{prop:implication_of_equivalences}
    Two $U$-equivalent TP tangles $K$ and $K'$ are ambient isotopic.
\end{proposition}
   
\begin{proof}
The transformations considered in an $a$-equivalence can all be regarded as ambient isotopies. What remains to prove is that the links in the $3$-torus $(\Gamma_k)_k$ of Definition \ref{def:u_equivalence}, correspond to unit cells of ambient isotopic TP tangles.

 Suppose, without loss of generality, that $\Gamma_{k-1}$ and $\Gamma_{k}$ are ambient isotopic. One can represent them as embeddings of curves in the unit cube with identified faces. By doing so, their preimages in the universal cover are naturally endowed with the lattice $\mathbb{Z}^3$ with its canonical basis, and are (periodically) ambient isotopic in $\mathbb{R}^3$.

Suppose $\Gamma_{k-1}$ and $\Gamma_{k}$ are $T^3$-equivalent. The torus twists can be regarded in $\mathbb{R}^3$ as changes of basis of space. Assuming, without loss of generality, that $\Gamma_{k-1}$ is represented by the unit cube with identified faces, there exists a parallelepiped whose delimiting vectors are related to the canonical basis of space by a matrix of $\mathrm{SL}(3,\mathbb{Z})$ that determines the torus twists connecting $\Gamma_{k-1}$ and $\Gamma_{k}$. By representing $\Gamma_{k}$ with this parallelepiped, the preimages in the universal cover of the two links in $\mathbb{T}^3$ actually are the exact same embedding of a TP tangle.

If $\Gamma_{k-1}$ and $\Gamma_{k}$ are $S$-equivalent, the link components can be connected by ambient isotopies of $\mathbb{T}^3$. Therefore, this case has already been treated previously.

Suppose $\Gamma_{k-1}$ and $\Gamma_{k}$ are $CM^3$-equivalent. By Definition \ref{def:cover_3}, there is a finite sequence $(\gamma_l)_l$ of links in the $3$-torus, where the ones are covers of the others, connecting $\Gamma_{k-1}$ and $\Gamma_{k}$. Without loss of generality, suppose that $\gamma_{l}$ is a cover of $\gamma_{l-1}$ where the covering map is $(z_1,z_2,z_3)  \longmapsto  ({z_1}^{h_{11}}{z_2}^{h_{21}}{z_3}^{h_{31}},{z_2}^{h_{22}}{z_3}^{h_{32}},{z_3}^{h_{33}})$. Assume also that $\gamma_{l-1}$ is represented by the unit cube. Then, by representing $\gamma_{l}$ by the parallelepiped delimited by the vectors $h_{11}e_1 + h_{21}e_2 + h_{31}e_3$, $h_{22}e_2 +h_{32}e_3$ and $h_{33}e_3$, the preimages in the universal cover $\mathbb{R}^3$, of $\gamma_{l-1}$ and $\gamma_{l}$ are actually the exact same embedding of a TP tangle. Continuing in a similar way for each element of the sequence $(\gamma_l)_l$, the preimages of $\Gamma_{k-1}$ and $\Gamma_{k}$ in $\mathbb{R}^3$ are the exact same TP tangle.
\end{proof}

\begin{remark}
    In the proof of Proposition \ref{prop:implication_of_equivalences}, we actually also demonstrated the converse of Proposition \ref{prop:unit_cells_are_equi}. Indeed, Proposition \ref{prop:unit_cells_are_equi} states that any two unit cells of the same TP tangle are connected by the $T^3$-equivalence, $S$-equivalence and $CM^3$-equivalence. The proof of Proposition \ref{prop:implication_of_equivalences} shows that any links in the $3$-torus connected by those three equivalences can be regarded as unit cells of one single TP tangle.
\end{remark}

Proposition \ref{prop:implication_of_equivalences} ensures that the $U$-equivalence is a particular case of ambient isotopy in $\mathbb{R}^3$. It would be interesting to see whether it is possible to approximate every ambient isotopy by the $U$-equivalence. However, that is beyond the scope of our interests in this paper.

From now on, we consider only the $U$-equivalence, which we simply call equivalence when there is no ambiguity.


\section{Reidemeister theorem for diagrams of 3-periodic tangles}\label{sec:3}

In this section, we formalise the diagrammatic description of 3-periodic tangles, and define the notions of equivalence for diagrams. We then prove a generalised Reidemeister theorem for 3-periodic tangles. These diagrams were suggested in \cite{Carrega2017} and defined in \cite{vuong2023fundamental} independently of this paper in the context of links in the 3-torus. Even though the results coincide, the methods used and the main goal remain different, as we are interested in the description of 3-periodic entanglements and not just links in the 3-torus, as well explained in Section \ref{sec:2}. Furthermore, in this paper we detail the proof of the generalised Reidemeister theorem.

\begin{remark}\label{rmk:tame_TP_tangle}
    Before explaining how to obtain a diagram, we would like to point out that in general not all 3-periodic tangles are worth our attention. One example that one wants to avoid is a periodic embedding of the wild knot invented by R.H. Fox \cite{burdezies_chap1}. Another example is an embedding of parallel lines along the $z$-axis passing through every point $(p,q,0) \in \mathbb{Q}.e_1 \oplus \mathbb{Q}.e_2$. The structure is 3-periodic, but within any of its unit cells, there are infinitely many components.

    Therefore, every TP tangle that we consider in the following will be simple enough so that every unit cell will be ambient isotopic to a polygonal link with finitely many components.
\end{remark}

\subsection{Diagrams}
Our aim is to encode the 3-dimensional topological information of TP tangles into 2-dimensional diagrams. In classical knot theory, one can obtain a diagram by projecting a knot onto a plane and by adding crossing information to the projection. In the following, we generalise this approach by making a projection of a unit cell of a TP tangle onto a 2-torus represented as a square with identified edges.

Projecting the 3-torus $ \mathbb{T}^3 \cong \mathbb{S}^1_1 \times \mathbb{S}^1_2 \times  \mathbb{S}^1_3$ represented as a cube with identified faces onto a 2-torus, say $\mathbb{T}_1^2 \cong \mathbb{S}^1_1 \times \mathbb{S}^1_2$, can be seen as flattening $\mathbb{S}^1_3$. It is important to keep track of the two identified faces along $\mathbb{S}^1_3$ as they are indistinct on the 2-torus. Just as with crossings for usual knots with an over-strand and an under-strand, a point intersecting those faces must also carry the information of which part is over and which is under. To capture this, we start with parametrising $ \mathbb{T}^3 \cong \mathbb{S}^1_1 \times \mathbb{S}^1_2 \times  \mathbb{S}^1_3$ as follows.

\begin{lemma}\label{lem:circle_is_interval_and_point}
    The $3$-torus is fully described by two smooth maps: $f_N: \mathbb{S}^1 \times \mathbb{S}^1 \times \left( \mathbb{S}^1 \setminus \lbrace N \rbrace \right)\longrightarrow {\mathbb{S}^1 \times \mathbb{S}^1 \times (0,1)}_N$, $f_S : \mathbb{S}^1 \times \mathbb{S}^1 \times \left( \mathbb{S}^1 \setminus \lbrace S \rbrace \right) \longrightarrow {\mathbb{S}^1 \times \mathbb{S}^1 \times (0,1)}_S$, and a smooth map transitioning between them, where $f_N$ and $f_S$ are obtained from the stereographic projections of the circle from the North and South poles.
\end{lemma}

\begin{proof}
    We recall that $\mathbb{S}^1$ is smoothly diffeomorphic to a union of an open interval and a point $(0,1) \cup \lbrace N \rbrace$. One possibility to get this is by considering an atlas given by two stereographic projections from the North and South poles, each of which transforms the punctured circle into the real line $\mathbb{R}$ \cite{Lee2011}. By composing one stereographic map with the smooth map $s: x \longmapsto \frac{1}{\pi} \arctan(x) + \frac{1}{2}$, one gets the open interval $(0,1)$. This gives rise to two smooth maps $f_N: \mathbb{S}^1 \times \mathbb{S}^1 \times \left( \mathbb{S}^1 \setminus \lbrace N \rbrace \right)\longrightarrow {\mathbb{S}^1 \times \mathbb{S}^1 \times (0,1)}_N$ and $f_S : \mathbb{S}^1 \times \mathbb{S}^1 \times \left( \mathbb{S}^1 \setminus \lbrace S \rbrace \right) \longrightarrow {\mathbb{S}^1 \times \mathbb{S}^1 \times (0,1)}_S$ and a smooth map transitioning between them.
\end{proof}

We represent ${\mathbb{S}^1_1 \times \mathbb{S}^1_2 \times (0,1)}_N$ as a cube with identified opposite faces, save the front and back ones, where we define the front and back faces of a 3-torus and of ${\mathbb{S}^1_1 \times \mathbb{S}^1_2 \times (0,1)}_N$ as follows.

\begin{definition}\label{def:front_back_faces}
Consider the 3-torus $ \mathbb{T}^3 \cong \mathbb{S}^1_1 \times \mathbb{S}^1_2 \times  \mathbb{S}^1_3$. The \textit{front face} of the torus is defined as the submanifold $\mathbb{S}^1_1 \times \mathbb{S}^1_2 \times \lbrace 0+\varepsilon \rbrace$ of ${\mathbb{S}^1_1 \times \mathbb{S}^1_2 \times (0,1)}_N$, and the \textit{back face} is defined as $\mathbb{S}^1_1 \times \mathbb{S}^1_2 \times \lbrace 1-\varepsilon \rbrace$, with $\varepsilon > 0$ arbitrarily small. The submanifold $\mathbb{S}^1_1 \times \mathbb{S}^1_2 \times \lbrace N\rbrace$ of ${\mathbb{S}^1_1 \times \mathbb{S}^1_2 \times (0,1)}_S$ is referred to as the {\it $N$-face}.
\end{definition}

Let $\Gamma$ be a link in $\mathbb{T}^3$. It follows that the embedding of $\Gamma$ in $\mathbb{T}^3$ has two descriptions via the two maps $f_N$ and $f_S$ of Lemma \ref{lem:circle_is_interval_and_point}. Note that as those two maps and the map transitioning between them are all smooth, they preserve the smoothness of the curves when $\Gamma$ is smoothly embedded.

\begin{figure}[hbtp]
\centering
\includegraphics[width= 11cm]{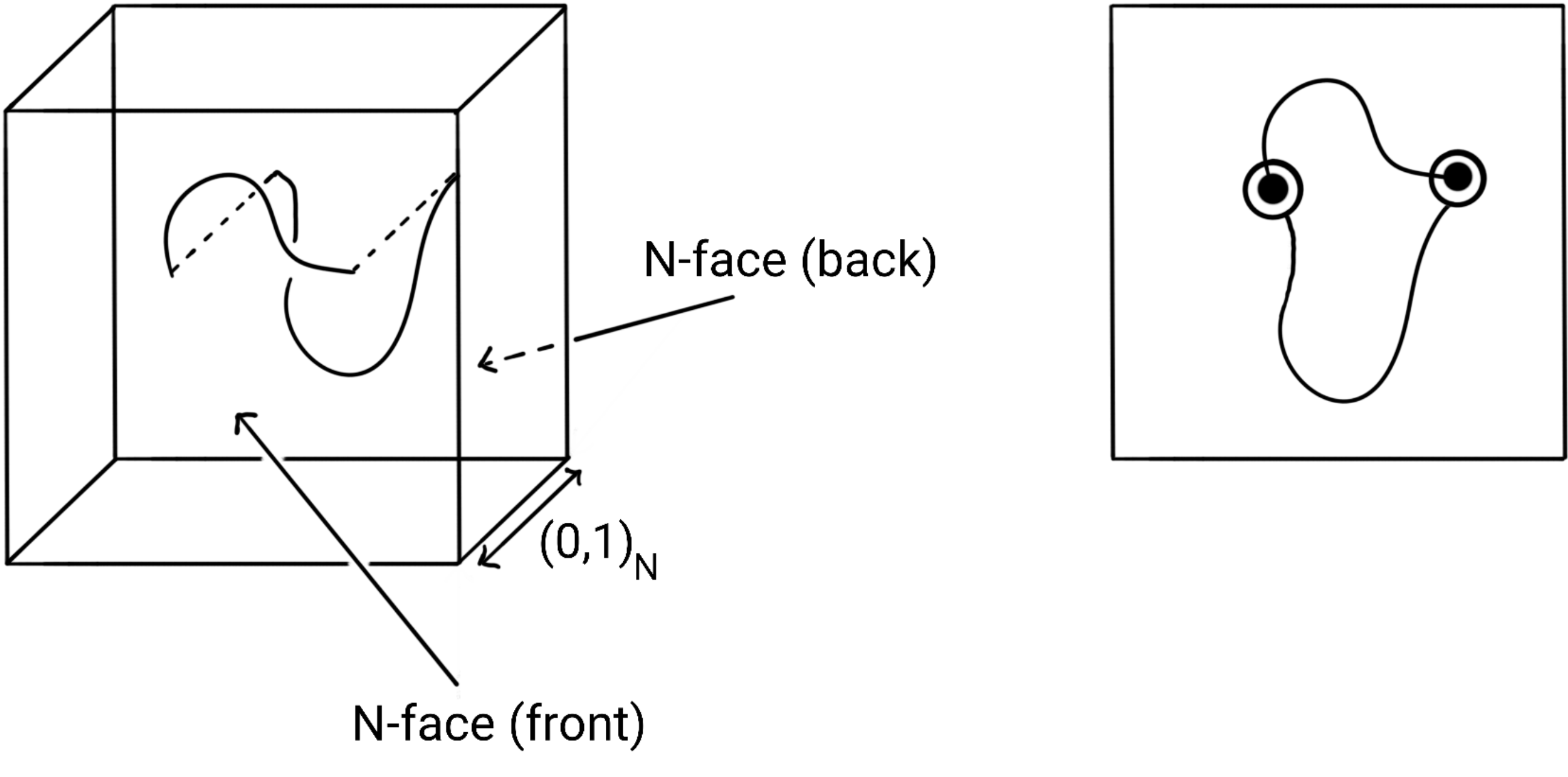}
\caption{Projecting a diagram of a 3-periodic tangle: on the left, a link with a single component in the 3-torus. One circle of the torus is represented as the union of the open interval (0,1) and the North point $N$. The identified front and back faces are called the $N$-face. On the right, the resulting diagram with a symbol representing the points belonging to the $N$-face.}
\label{fig:dia_proj}
\end{figure}

Consider the projection map  $\pi : \mathbb{T}^3 \longrightarrow \mathbb{T}_1^2 $ which is described by the maps
$ \pi_N : {\mathbb{S}^1_1 \times \mathbb{S}^1_2 \times (0,1)}_N \longrightarrow \mathbb{T}_1^2 $ via $f_N$ and $ \pi_S : {\mathbb{S}^1_1 \times \mathbb{S}^1_2 \times (0,1)}_S \longrightarrow \mathbb{T}_1^2 $ via $f_S$. Visually, we consider the information given by $\pi_N$ as in Figure \ref{fig:dia_proj}. The map $\pi_S$ is used only to complete the missing pieces of information on the $N$-face. By doing so, one can understand the `height' of a curve, especially regarding the front face and back face from Definition \ref{def:front_back_faces}. If one considers ${\mathbb{S}^1_1 \times \mathbb{S}^1_2 \times (0,1)}_N \longrightarrow (0,1)$, one can regard it as a map that determines the `height' of $\mathbb{S}^1_1 \times \mathbb{S}^1_2 \times \lbrace t \rbrace$ on $(0,1)$. Visualising the map $\mathbb{S}^1_1 \times \mathbb{S}^1_2 \times \mathbb{S}^1_3 \longrightarrow \mathbb{S}^1_3$ as a `height map' is difficult as $\mathbb{S}^1_3$ loops around. Given this, we can introduce new symbols to represent a point intersecting the $N$-face.

\begin{figure}[hbtp]
\centering
\includegraphics[width= 7cm]{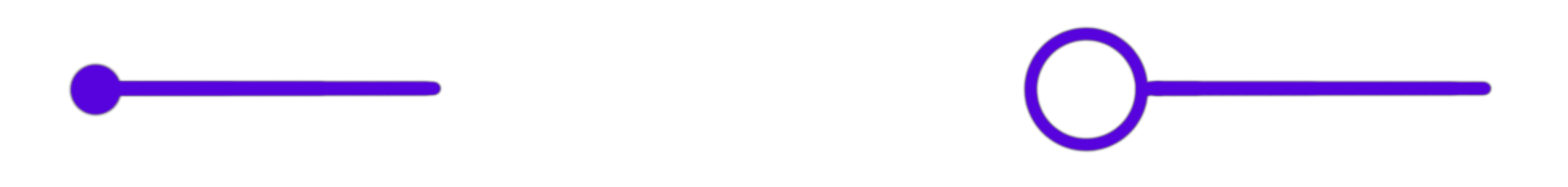}
\caption{New symbols: on the left, the $F$-point, representing the part of the curve that intersects the front face $\mathbb{S}^1_1 \times \mathbb{S}^1_2 \times \lbrace 0+\varepsilon \rbrace$. On the right, the $B$-point, representing the part of the curve that intersects the back face $\mathbb{S}^1_1 \times \mathbb{S}^1_2 \times \lbrace 1-\varepsilon \rbrace$.}
\label{fig:newsignsdiagram}
\end{figure}

\begin{definition}\label{definition_F_B_N_points}
    The thick dot in Figure \ref{fig:newsignsdiagram} left, called \textit{$F$-point}, depicts an intersection with the front face, and the circle in Figure \ref{fig:newsignsdiagram} right, called \textit{$B$-point}, depicts an intersection with the back face. Together, they form an \textit{$N$-point} which is the projection of a point intersecting the $N$-face.
\end{definition}

Diagrams of knotted structures are obtained from \textit{regular projections}, avoiding ambiguity in aspects like the nature of a crossing of a triple point (a point $P \in \pi(l)$ with $\pi^{-1}(P)$ containing three points) \cite{burdezies_chap1}. For similar reasons, we want to define a notion of regular projection for TP tangles. To do so, we can restrict ourselves to studying smoothly embedded curves since the set of such curves is dense in the set of embeddings of curves \cite{Hirsch1976_chap2}.

We make use of the genericness of transversality, meaning that the set of transverse maps with respect to a manifold is dense in the set of smooth maps. This is a consequence of Thom's transversality theorem which states the following. Suppose $f_t: X \longrightarrow Y$ is a family of smooth maps indexed by a parameter $t$ that ranges over some set $T$. Consider the map $F: X\times T \longrightarrow Y$ defined by $F(x,t) = f_t(x)$. We require that the family varies smoothly by assuming $T$ to be a manifold and $F$ to be smooth.

\begin{theorem} \label{thm:thom_transversality}
(Thom, \cite{thom_transversality}). Suppose that $F: X\times T \longrightarrow Y$ is a smooth map of manifolds, where only $X$ has boundary, and let $Z$ be any boundaryless submanifold of $Y$. If both $F$ and $\partial F$ are transversal to $Z$, then for almost every $t \in T$, both $f_t$ and $\partial f_t$ are transversal to $Z$.
\end{theorem}

\begin{figure}[hbtp]
\centering
\includegraphics[width= 10cm]{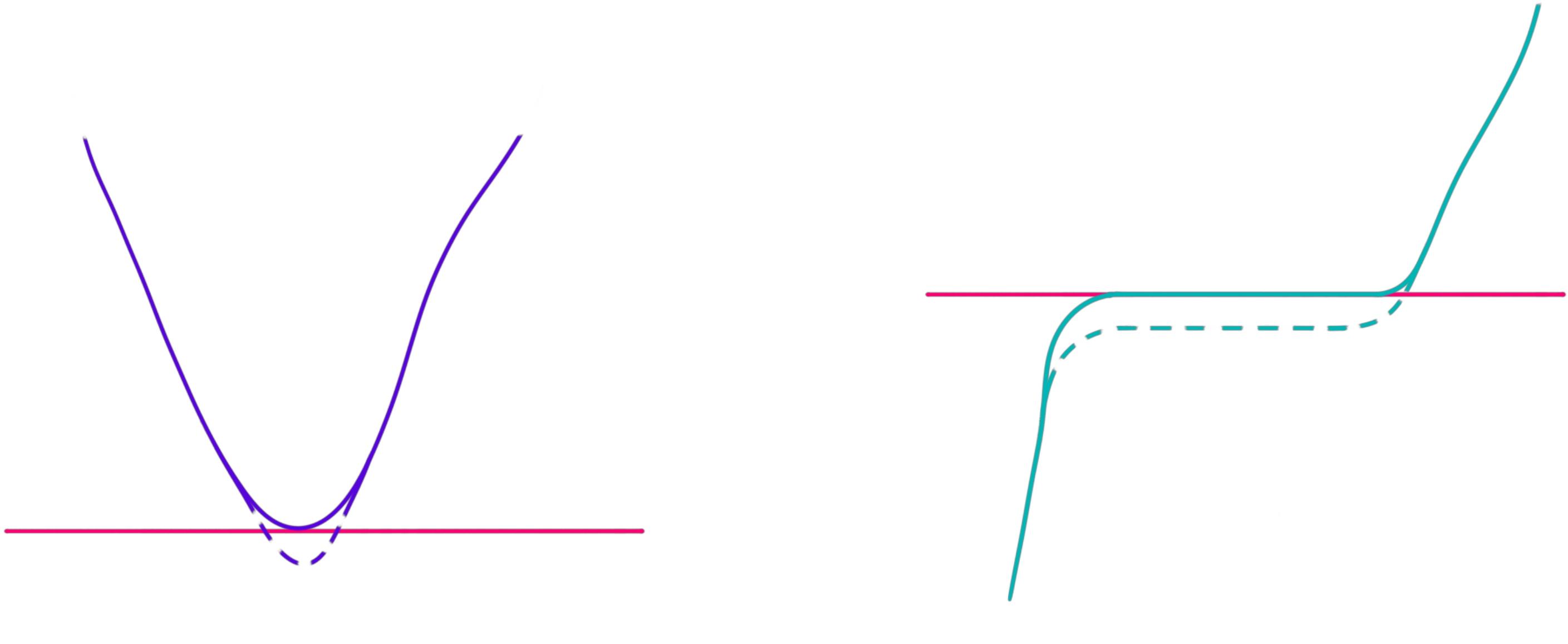}
\caption{Deforming curves to get transverse intersections.}
\label{fig:transversality}
\end{figure}

A proof of Theorem \ref{thm:thom_transversality} can be found in \cite{guillemin_differential_topology_chap2,Hirsch1976_chap3}. In our case the transversality theorem implies that any curve in $\mathbb{T}^2$ may be deformed by an arbitrarily small amount to have a transverse intersection with a second curve as in Figure \ref{fig:transversality}.

This allows us to state the following definition.

\begin{definition}\label{def:general_position}
\textbf{(Regular projection)}\\
Consider a link $\Gamma$ embedded in the $3$-torus, represented by a cube with identified faces, where the edges are considered as the generating circles $l_1$, $l_2$ and $l_3$. Consider a projection $\pi: \mathbb{T}^3 \longrightarrow \mathbb{T}^2$, where $\mathbb{T}^2$ is represented as a square with identified edges, where the edges are its generating circles $l_1$ and $l_2$. We say that the projection $\pi$ of $\Gamma$ is \textit{regular}, or that $\Gamma$ is in \textit{general position} with respect to the projection $\pi$, if $\pi$ satisfies the following eight requirements:
\begin{enumerate}
\item there are only finitely many multiple points $\lbrace P_i : 1 \leqslant i \leqslant n \rbrace $, and all multiple points are double points, that is $\char"0023 \pi^{-1}(P_i)=2 $,
\item all intersections, that is, double points, are transverse,
\item no point in the $N$-face is mapped onto a double point,
\item in $\mathbb{T}^3$, all intersections of $\Gamma$ with the $N$-face are transverse,
\item all intersections with $l_1$ and $l_2$ are transverse,
\item no double point of $\pi(\Gamma)$ lies on $l_1$ and $l_2$,
\item as $l_1$ and $l_2$ intersect once, $\pi(\Gamma)$ does not intersect $l_1$ and $l_2$ at their intersection point,
\item no point belonging to the $N$-face is projected onto $l_1$ and $l_2$.
\end{enumerate}
\end{definition}

\begin{definition}\label{def:def_diagram}
A \textit{unit-diagram} of a TP tangle, or simply \textit{diagram}, is a regular projection of one of its unit cells, with crossing information and $N$-points.
\end{definition}

\begin{remark}
Usually, the diagram of a given knotted object is obtained from a projection of that object. We do not give a similar definition for TP tangles in this paper.
\end{remark}

\begin{figure}[hbtp]
    \centering
    \begin{subfigure}[b]{\textwidth}
        \centering
        \includegraphics[width=0.8\textwidth]{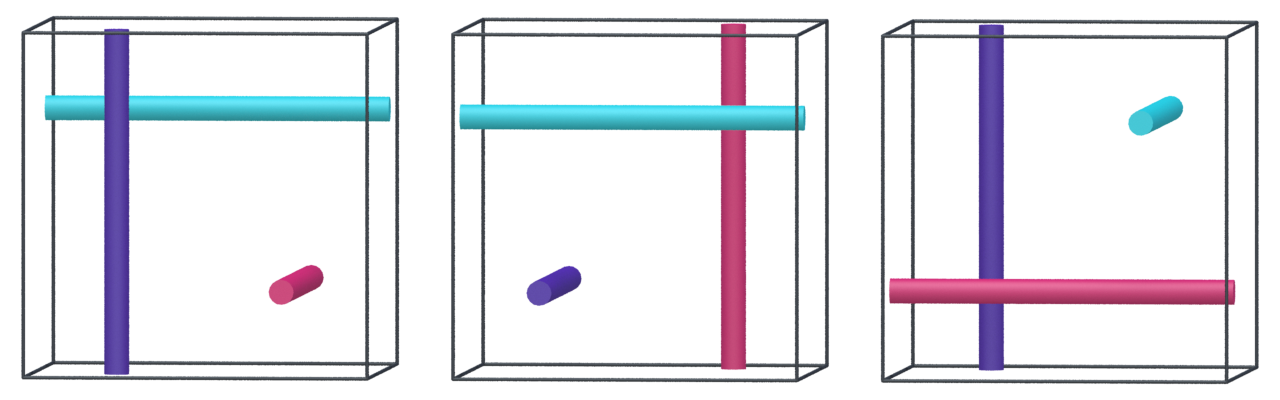}
        \caption{}
        \label{fig:pi_star_uc_xyz}
    \end{subfigure}
    \vskip\baselineskip
    \begin{subfigure}[b]{\textwidth}
        \centering
        \includegraphics[width=0.75\textwidth]{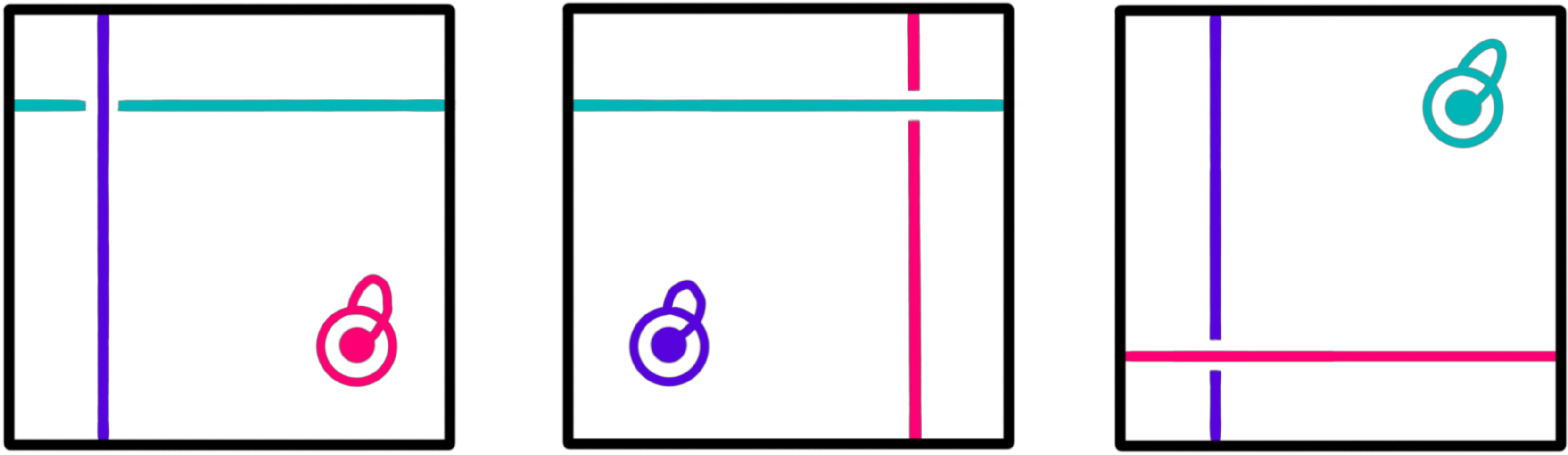}
        \caption{}
        \label{fig:pi_star_tridiagram}
    \end{subfigure}
    \caption{An example of a tridiagram of a TP tangle: (a) The unit cell of Figure \ref{fig:pi_star_motifs_comparison} left, viewed along three directions. (b) A tridiagram obtained from the unit cell of (a), where the diagrams are projected from the front, the top and the right faces.}
    \label{fig:pi_star_xyz_and_tridia}
\end{figure}

\begin{remark}\label{rmk:diagram_loss_of_geometry}
    Unless otherwise mentioned, a diagram of a TP tangle is always represented by a square with identified edges. By doing so, we actually lose some geometric information. For example, the two unit cells of Figure \ref{fig:pi_star_motifs_comparison} are both represented by the same diagrams that we display in Figure \ref{fig:pi_star_tridiagram}. This implies that the $a$-equivalence is not detected by diagrams.
\end{remark}

One single diagram is enough to fully characterise an embedding of a link in the $3$-torus. However, as will be seen later, the generalised Reidemeister theorem for TP tangles necessitates the consideration of three diagrams projected along three non-coplanar axes. This motivates the following definition.

\begin{definition}
    A \textit{tridiagram} of a 3-periodic tangle is an ordered $3$-tuple of diagrams obtained from three regular projections of a unit cell along three non-coplanar axes.
\end{definition}

An example of a tridiagram is given in Figure \ref{fig:pi_star_xyz_and_tridia}. Unless otherwise mentioned, every tridiagram in this paper is an ordered 3-tuple of diagrams projected respectively from the front, the top and the right faces of a unit cell.

\begin{lemma}\label{lem:dia_has_tridia}
    From any diagram $A$ of a TP tangle, one can obtain a tridiagram $\left\lbrace A,B,C\right\rbrace$.
\end{lemma} 

\begin{proof}
    One can embed $A$ as a link in $\mathbb{T}^3$ and obtain two other diagrams from that embedding by considering two regular projections along two other axes, the three axes forming a basis of space. From those three diagrams, one obtains a tridiagram $\lbrace A,B,C \rbrace$.
\end{proof}

\begin{remark}
    The associated tridiagram to a diagram $A$ need not be unique, as there are infinitely many ways of embedding $A$.
\end{remark}

\subsection{Notions of equivalence of diagrams}
\begin{figure}[ht]
    \centering
    \begin{subfigure}[b]{5.5cm}
        \includegraphics[width=\textwidth]{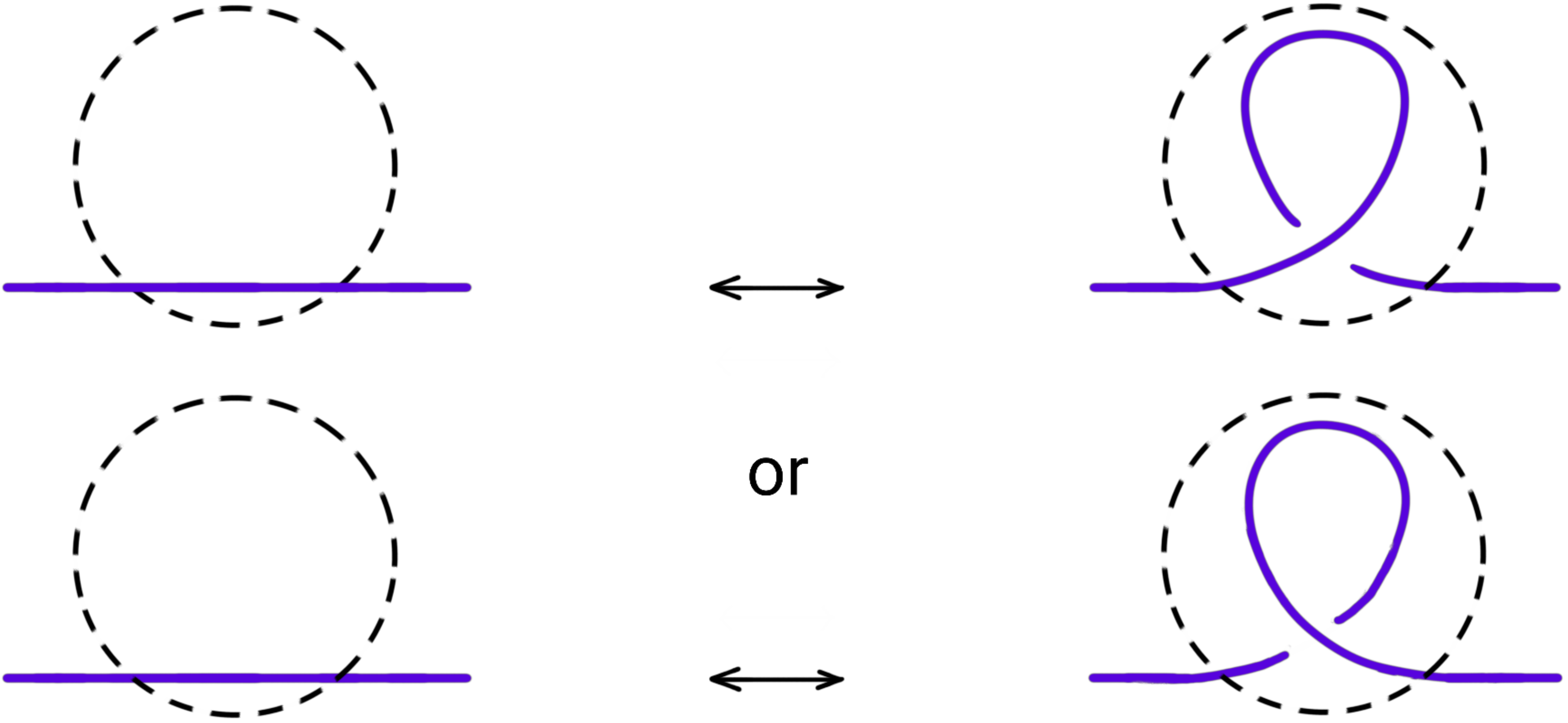}
        \caption{The $R_1$ move}
        \label{fig:Reid_move_I}
    \end{subfigure}
    \hspace{0.5cm}
    \begin{subfigure}[b]{6cm}
        \includegraphics[width=\textwidth]{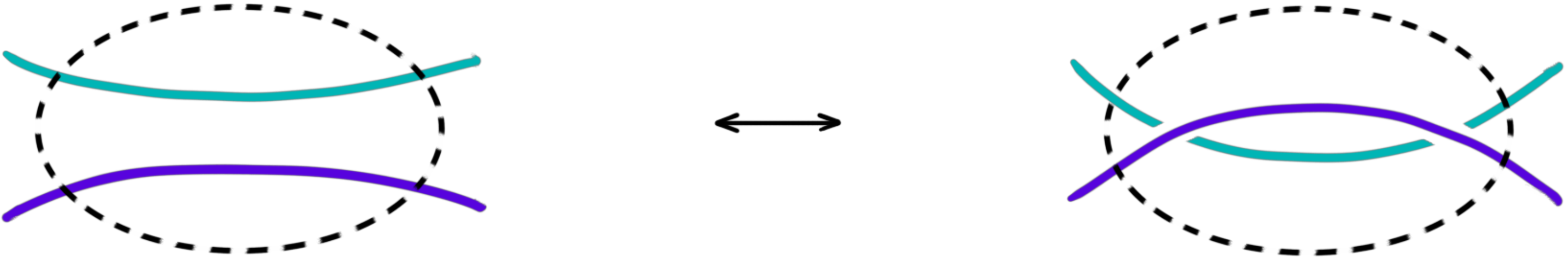}
        \caption{The $R_2$ move}
        \vspace{0.75cm}
        \label{fig:Reid_move_II}
    \end{subfigure}
    \vskip\baselineskip
    \begin{subfigure}[b]{5.5cm}
        \includegraphics[width=\textwidth]{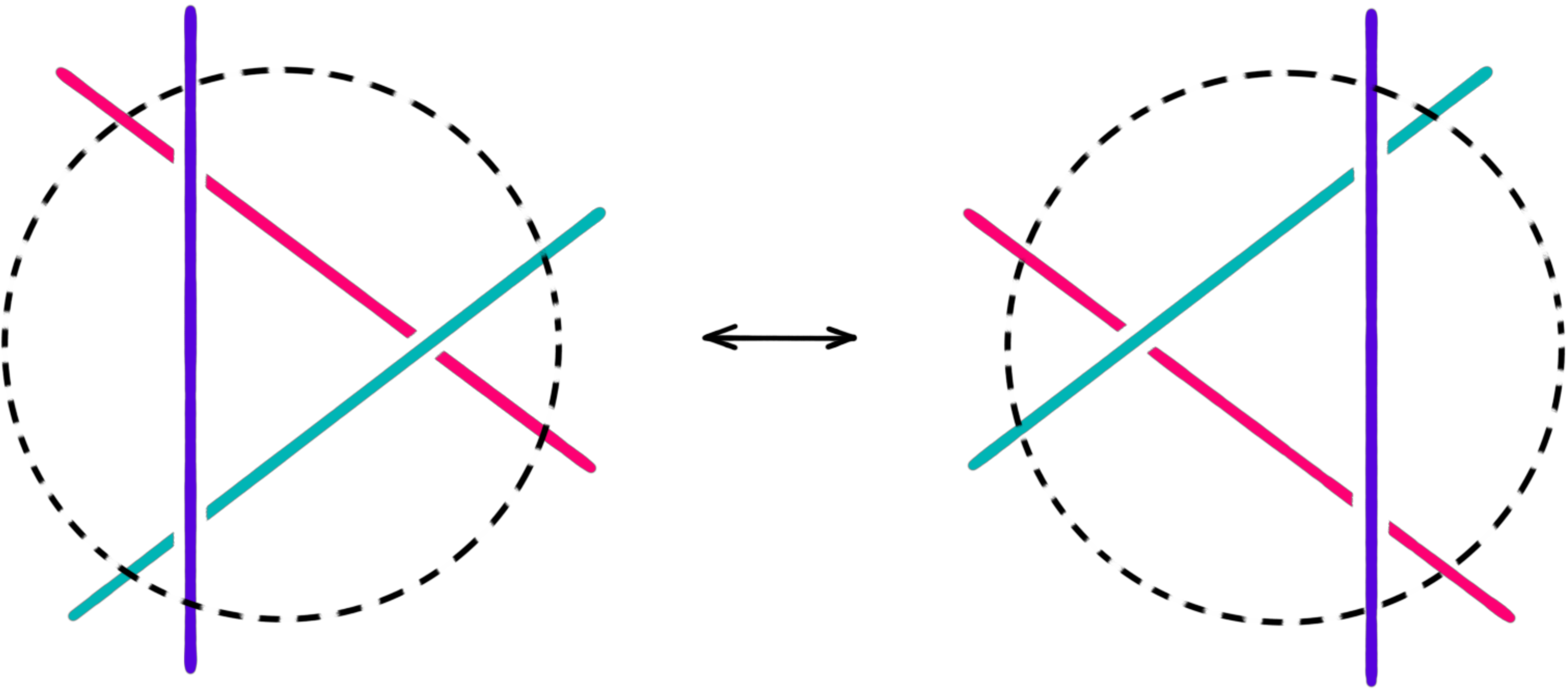}
        \caption{The $R_3$ move}
        \label{fig:Reid_move_III}
    \end{subfigure}
    \caption{The three Reidemeister moves.}
    \label{fig:usual_Reid_moves}
\end{figure}

Ambient isotopies of usual knots in $\mathbb{R}^3$ are translated at the diagrammatic level to planar isotopies and Reidemeister moves \cite{burdezies_chap1}. We recall that the first Reidemeister move ($R_1$) allows us to put in or take out a twist in the knot, which also adds or removes a crossing as shown in Figure \ref{fig:Reid_move_I}. The second Reidemeister move ($R_2$) adds or removes two crossings as displayed in Figure \ref{fig:Reid_move_II}. The third Reidemeister move ($R_3$) allows us to slide a strand of the knot from one side of a crossing to the other side of the crossing, as depicted in Figure \ref{fig:Reid_move_III}.

In the case of TP tangles, the diagrams are delimited by squares, and we also introduced new symbols (Figure \ref{fig:newsignsdiagram}) to describe the direction of periodicity along which the projection is made. Due to these new considerations, new moves arise to fully encode ambient isotopies in the $3$-torus.

\begin{figure}[hbtp]
    \centering
    \begin{subfigure}[b]{6cm}
    \centering
        \includegraphics[width=\textwidth]{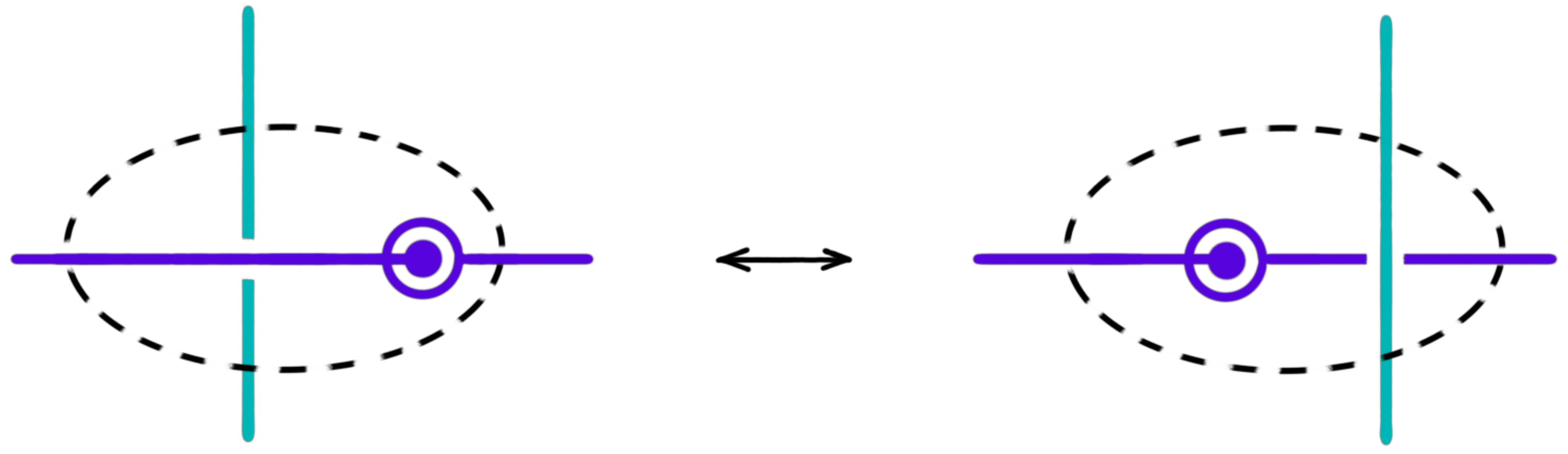}
        \caption{The $R_4$ move}
        \label{fig:Reid_move_IV}
        \vspace{0.5cm}
    \end{subfigure}
    \hspace{0.5cm}
    \begin{subfigure}[b]{6cm}
    \centering
        \includegraphics[width=\textwidth]{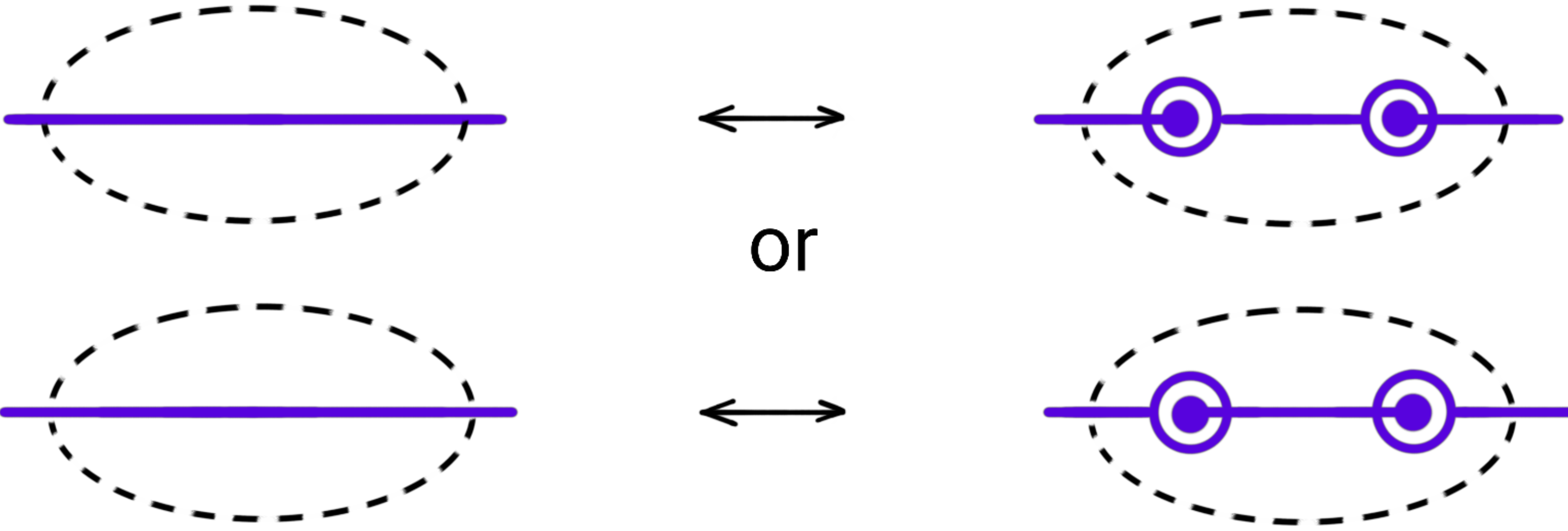}
        \caption{The $R_5$ move}
        \label{fig:Reid_move_V}
    \end{subfigure}
    \vskip\baselineskip
    \begin{subfigure}[b]{5.5cm}
    \centering
        \includegraphics[width=\textwidth]{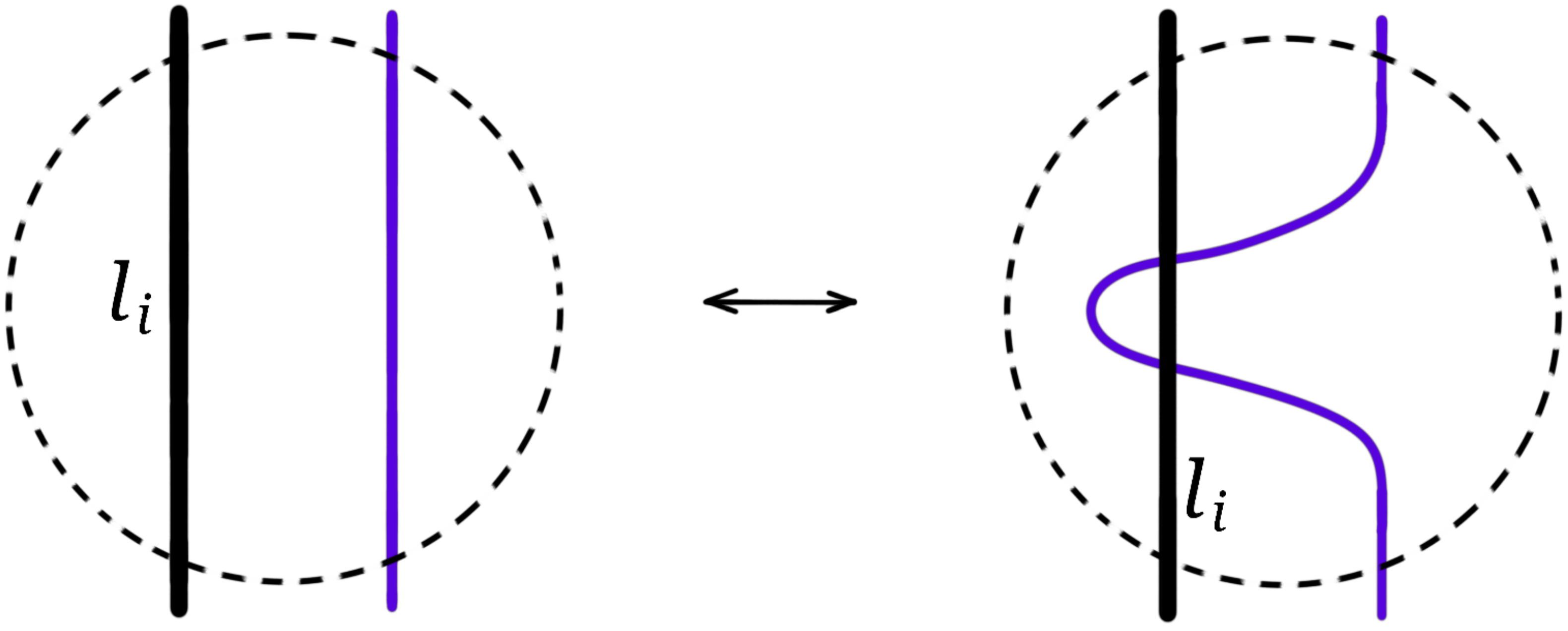}
        \caption{The $R_6$ move}
        \vspace{1.25cm}
        \label{fig:Reid_move_VI}
    \end{subfigure}
    \hspace{1cm}
    \begin{subfigure}[b]{5.5cm}
    \centering
        \includegraphics[width=\textwidth]{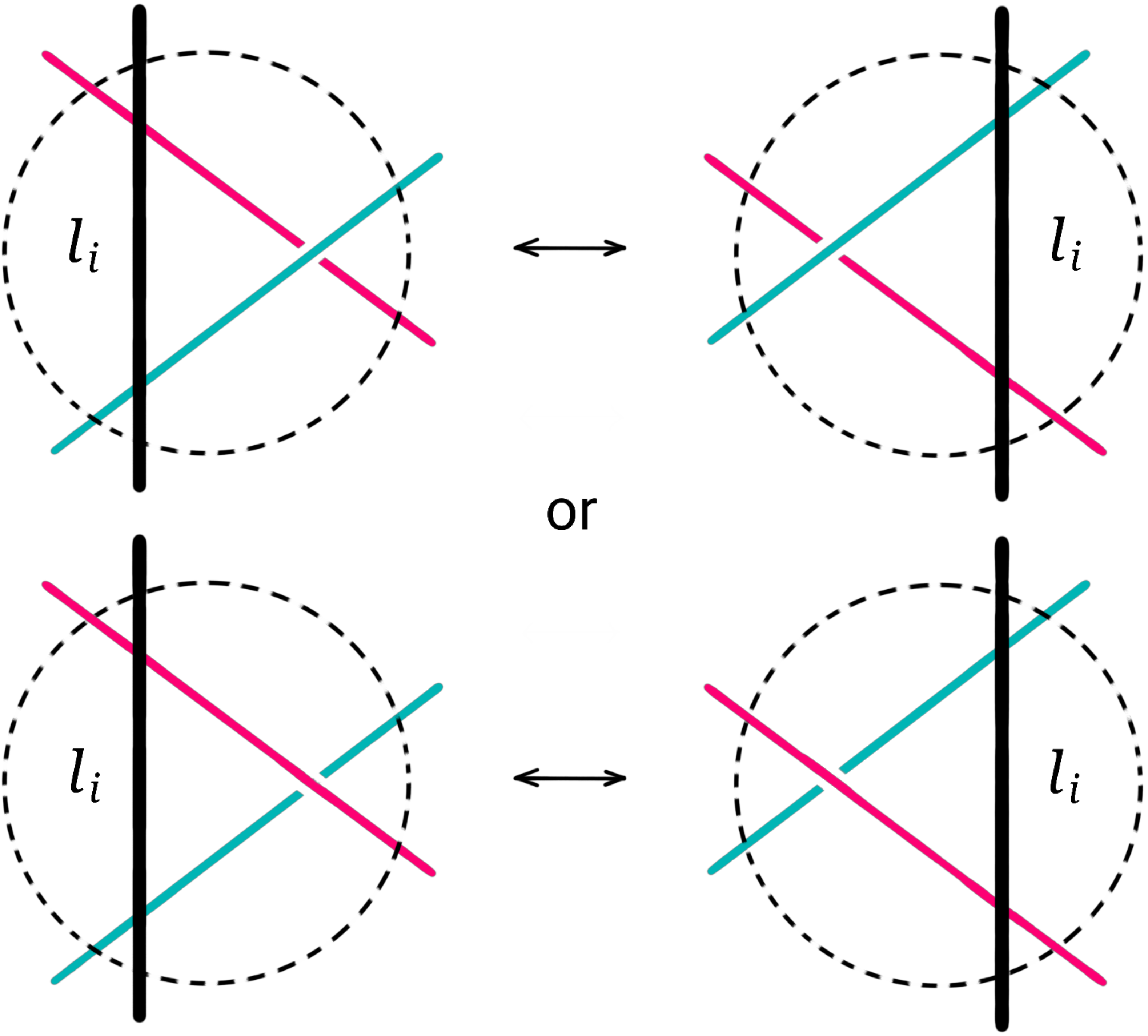}
        \caption{The $R_7$ move}
        \label{fig:Reid_move_VII}
    \end{subfigure}
    \vskip\baselineskip
    \begin{subfigure}[b]{5.5cm}
    \centering
        \includegraphics[width=\textwidth]{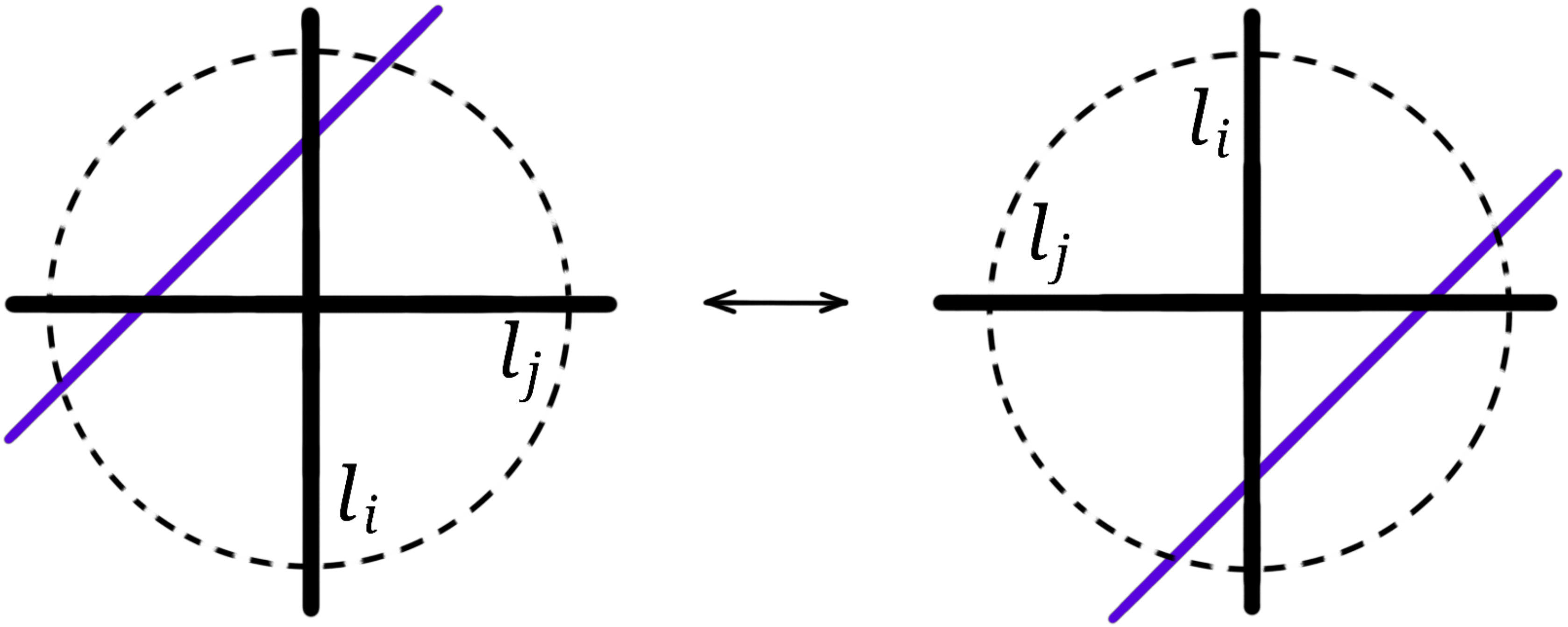}
        \caption{The $R_8$ move}
        \label{fig:Reid_move_VIII}
    \end{subfigure}
    \hspace{0.75cm}
    \begin{subfigure}[b]{6cm}
    \centering
        \includegraphics[width=\textwidth]{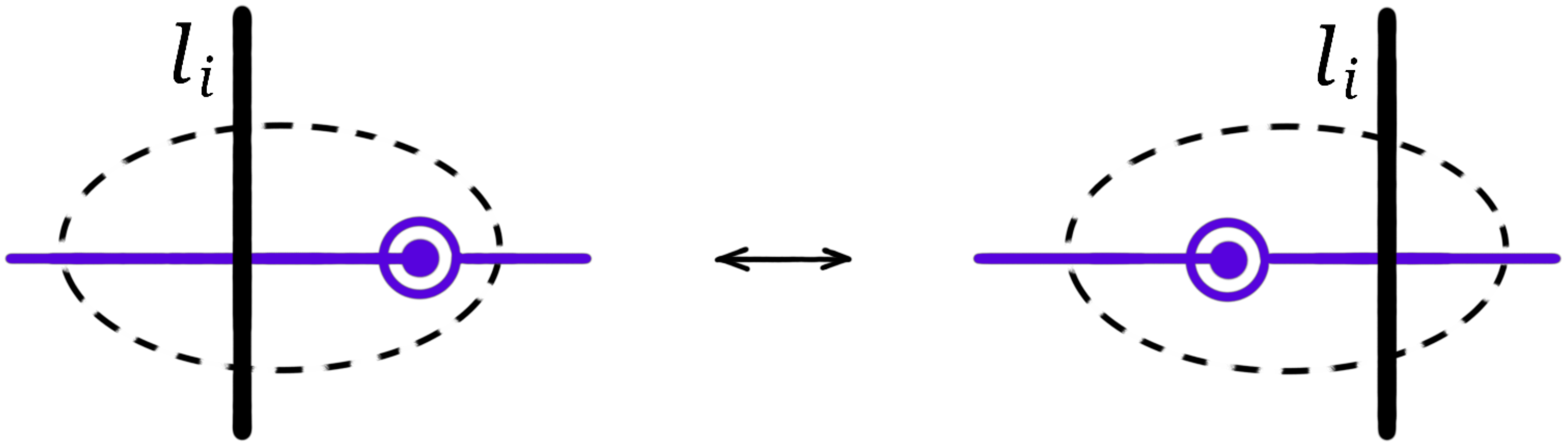}
        \caption{The $R_9$ move}
        \label{fig:Reid_move_IX}
        \vspace{0.375cm}
    \end{subfigure}
        \caption{The six new moves, necessary for 3-periodic tangles. Here $l_i$ and $l_j$, where $i,j = 1,2$, are the generating circles of the $2$-torus, that are also the edges of the square representing the torus.}
    \label{fig:new_Reid_moves}
\end{figure}

\begin{definition}
    The \textit{$R_4$ move} depicts a strand of the link which slides between the front part and back part of another strand that goes through the $N$-face as in shown Figure \ref{fig:Reid_move_IV}. The \textit{$R_5$ move} allows us to pass a part of a curve through the $N$-face as seen in Figure \ref{fig:Reid_move_V}.
    The \textit{$R_6$ move} depicts a curve that goes through one edge of the square representing the $2$-torus. The \textit{$R_7$ move} represents a crossing that goes through one edge of the square. The \textit{$R_8$ move} corresponds to a strand passing through the intersection of the two edges of the square. Finally, the \textit{$R_9$ move} depicts an $N$-point going through one edge of the square.
    We call \textit{$R$-moves} the set of nine moves given by the three Reidemeister moves and the $R_4$, $R_5$, $R_6$, $R_7$, $R_8$ and $R_9$ moves.
\end{definition}

\begin{remark}
    The $R_4$ and $R_5$ moves respectively mimic the isotopies that transform the unit cell of Figure \ref{fig:R_move_4_a} to that of Figure \ref{fig:R_move_4_b} (or vice-versa) and the unit cell of Figure \ref{fig:R_move_5_a} to that of Figure \ref{fig:R_move_5_b} (or vice-versa). We technically lose the information of the parity of the crossing with the $R_4$ move. This poses a challenge for defining invariants through the computation of signed crossings.

\begin{figure}[ht]
    \centering
    \begin{subfigure}[b]{4cm}
        \includegraphics[width=\textwidth]{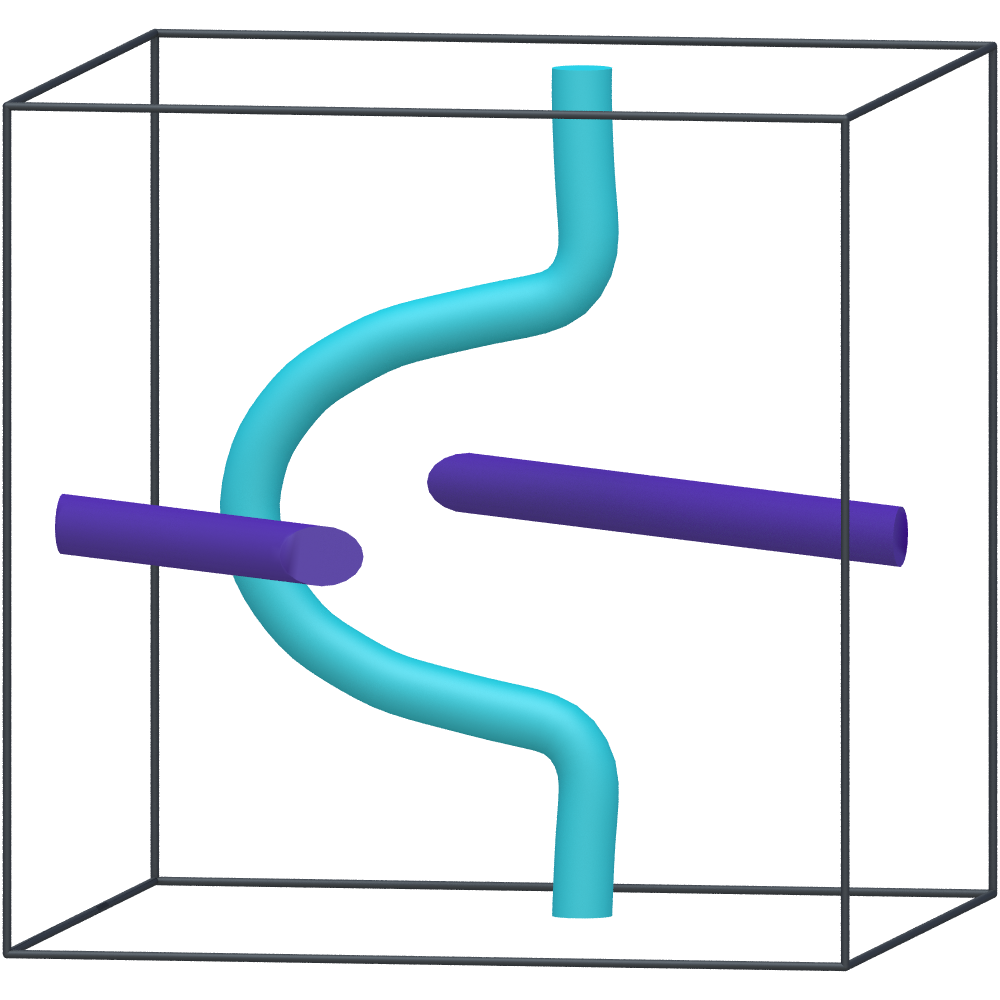}
        \caption{}
        \label{fig:R_move_4_a}
    \end{subfigure}
    \hspace{0.5cm}
    \begin{subfigure}[b]{4cm}        
        \includegraphics[width=\textwidth]{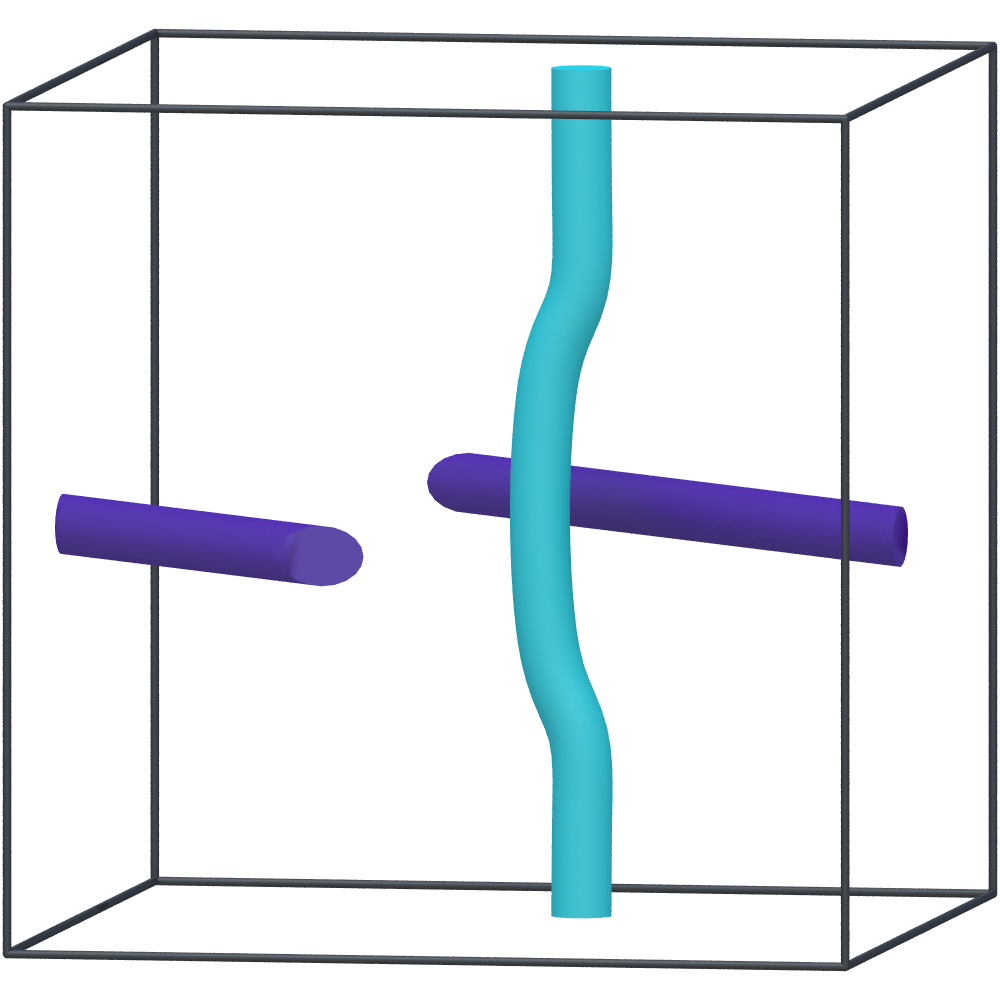}
        \caption{}
        \label{fig:R_move_4_b}
    \end{subfigure}
    \vskip\baselineskip
    \begin{subfigure}[b]{4cm}
        \includegraphics[width=\textwidth]{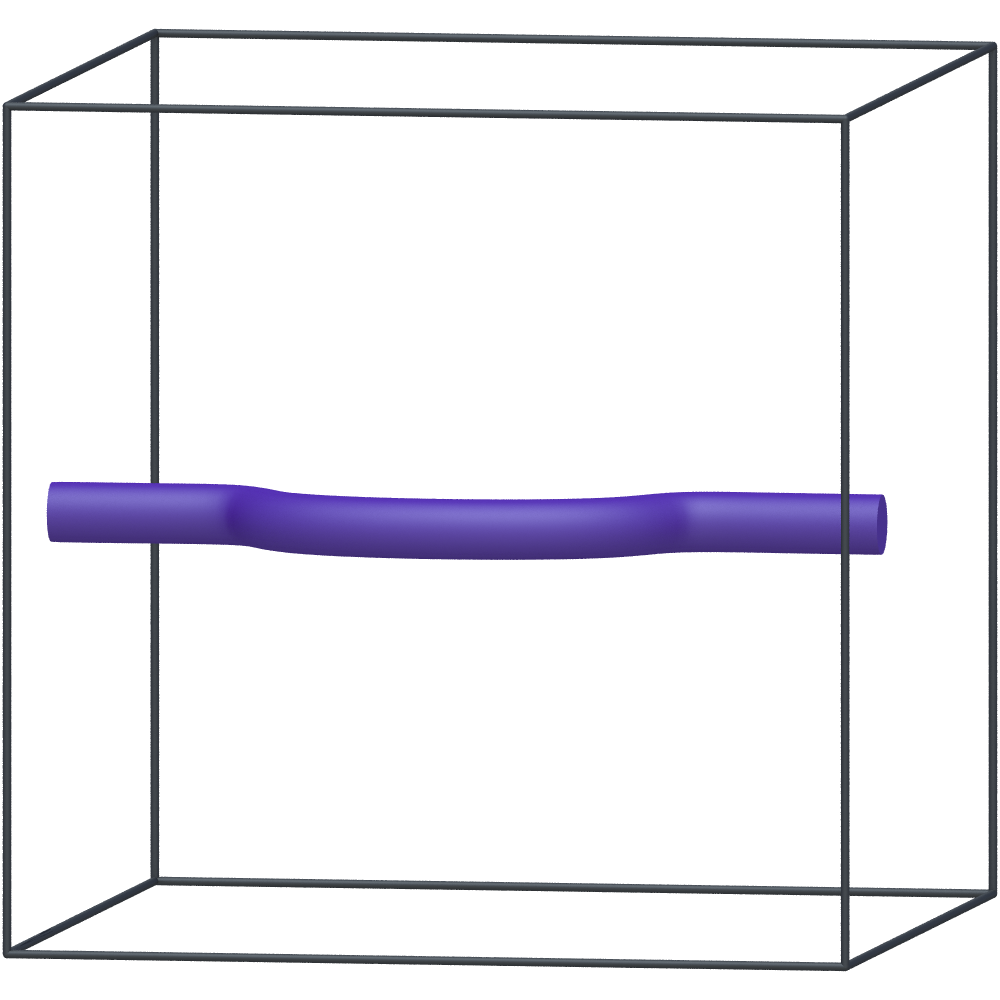}
        \caption{}
        \label{fig:R_move_5_a}
    \end{subfigure}
    \hspace{0.5cm}
    \begin{subfigure}[b]{4cm}
        \includegraphics[width=\textwidth]{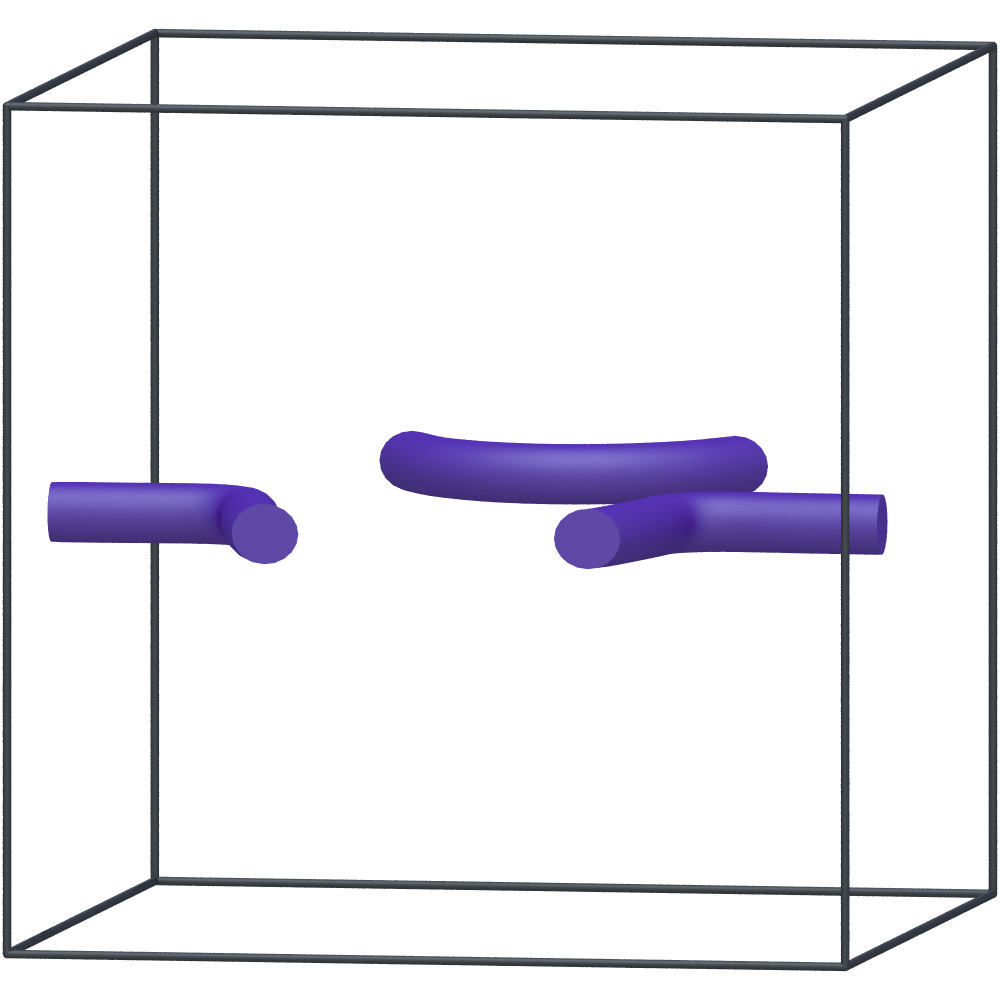}
        \caption{}
        \label{fig:R_move_5_b}
    \end{subfigure}
    \caption{Ambient isotopies that cannot be realised with the usual Reidemeister moves in a diagram: On the top, an isotopy from (a) to (b) or vice-versa that would correspond to the $R_4$ move in a diagram. On the bottom, an isotopy from (c) to (d) or vice-versa that would correspond to the $R_5$ move in a diagram. }
    \label{fig:R_moves}
\end{figure}

The $R_6$ move is similar to the $R_2$ move, The $R_7$ and $R_8$ moves are similar to the $R_3$ move, and the $R_9$ move is similar to the $R_4$ move. The difference is that there is no crossing involved with the edges $l_1$ and $l_2$.

The $R_6$ move also mimics the isotopy that transforms the unit cell of Figure \ref{fig:R_move_5_a} to that of Figure \ref{fig:R_move_5_b} (or vice-versa), but in a diagram projected from the top of the cube. This means that the $R_5$ and $R_6$ moves are two representations of the same isotopy regarded along different directions.
\end{remark}

\begin{remark}
    In Figure \ref{fig:Reid_move_I}, on the $R_1$ move, the second alternative can be obtained by a finite sequence of the first alternative and $R_2$ and $R_3$ moves \cite{Murasugi1996_chap4}. We list the two variants only because both are needed for the proof of the generalised Reidemeister theorem that will be given later in the paper, where an isotopy is locally translated to one variant of each $R$-move and not a sequence of variants. On the other hand, the two alternatives for an $R_5$ move in Figure \ref{fig:Reid_move_V}, or an $R_7$ move if Figure \ref{fig:Reid_move_VII}, are not connected by any sequence of other $R$-moves.
\end{remark}

\begin{remark}
    The $R_4$ and $R_5$ moves are similar to the $\Omega_4$ and $\Omega_5$ moves introduced in \cite{MROCZKOWSKI20091831} in the more general context of $F \times \mathbb{S}^1$ manifolds where $F$ is an orientable surface. The $R_6$, $R_7$ and $R_8$ moves are similar to the moves of mixed link diagrams \cite{LAMBROPOULOU199795}, and the moves of DP tangle diagrams \cite{diamantis2023equivalence}.
\end{remark}

\begin{definition}\label{def:r_equi}
    Two diagrams $D$ and $D'$ are \textit{$R$-equivalent} if they are connected by a finite sequence of $R$-moves. Furthermore, two tridiagrams $\left\lbrace D_i \right\rbrace_{i=1,2,3}$ and $\left\lbrace D_i' \right\rbrace_{i=1,2,3}$ are \textit{$R$-equivalent} if there is a finite sequence of tridiagrams $\left( \lbrace D_{i,k} \rbrace_{i=1,2,3} \right)_{k=0,\dots,n}$, for which we have 
    \[  D_{i,0}  =   D_{i} , \quad  D_{i,n}  =  D_{i}^{\prime}, \quad i=1,2,3\]
    and 
    \[\forall k = 1,\dots,n, \exists i,j = 1,2,3: \quad \text{ $D_{i,k-1}$ and $D_{j,k}$ are $R$-equivalent}.\]
    The equivalence is denoted by $\sim_R$ for both cases.
\end{definition}

\begin{remark}
    Strictly speaking, $R$-moves are $2$-dimensional moves, and thus performed on a diagram and not on a tridiagram. The example of Figure \ref{fig:r_moves_tridia} gives further motivation for that consideration. In Figure \ref{fig:r_moves_tridia_1}, we display a tridiagram that represents a unit cell, with each diagram projected respectively from the front, the top and the right faces. In Figure \ref{fig:r_moves_tridia_2}, an $R_2$ move is performed on the front diagram. From the resulting diagram, one projects the tridiagram displayed in Figure \ref{fig:r_moves_tridia_3}. One notices that the $R_2$ move performed on the front diagram induces another $R_2$ move in the right diagram. This means that $R$-moves in one diagram are not independent from $R$-moves in another diagram. As it is difficult to predict what sequence of $R$-moves are induced, it is better to always perform $R$-moves on a single diagram first, and then project the corresponding tridiagram afterwards, rather than attempting to perform $R$-moves directly on the tridiagrams. This also means that changing the number of crossings in one diagram with an $R$-move may or may not change the number of crossings in other diagrams.

The example also shows that not any ordered 3-tuple of diagrams can constitute a tridiagram. Indeed, although the three diagrams of Figure \ref{fig:r_moves_tridia_4} represent the same equivalence class of unit cells, there is no unit cell that projects into those three diagrams at the same time.\\

\begin{figure}[hbtp]
    \centering
    \begin{subfigure}[b]{0.65\textwidth}
        \includegraphics[width=\textwidth]{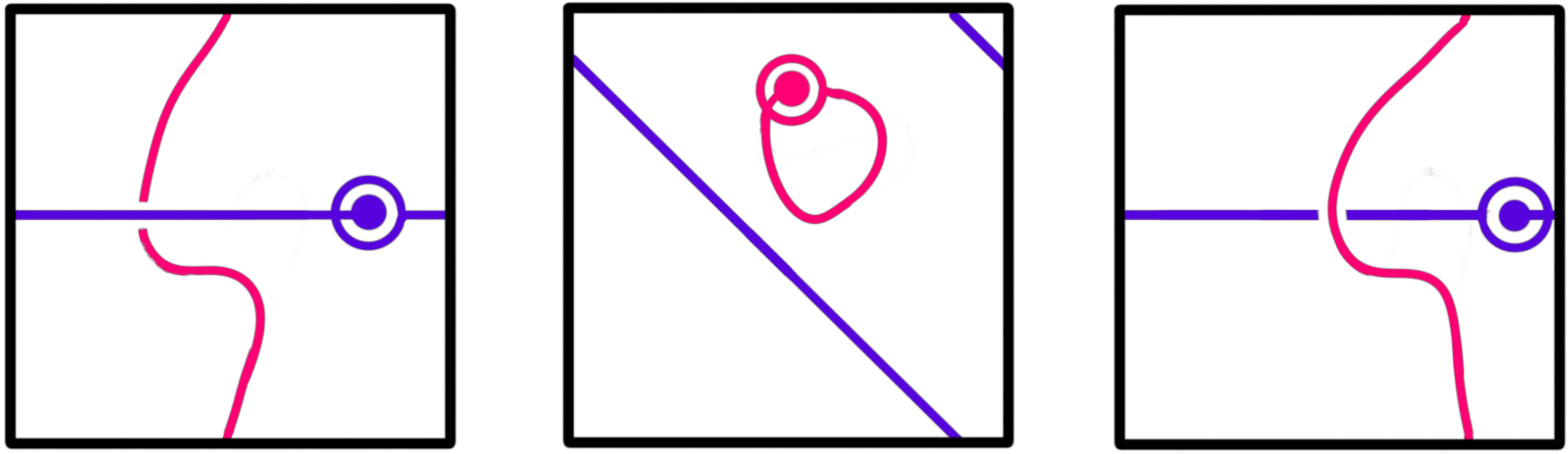}
        \caption{}
        \label{fig:r_moves_tridia_1}
    \end{subfigure}
    \vskip\baselineskip
    \begin{subfigure}[b]{0.65\textwidth}
        \includegraphics[width=\textwidth]{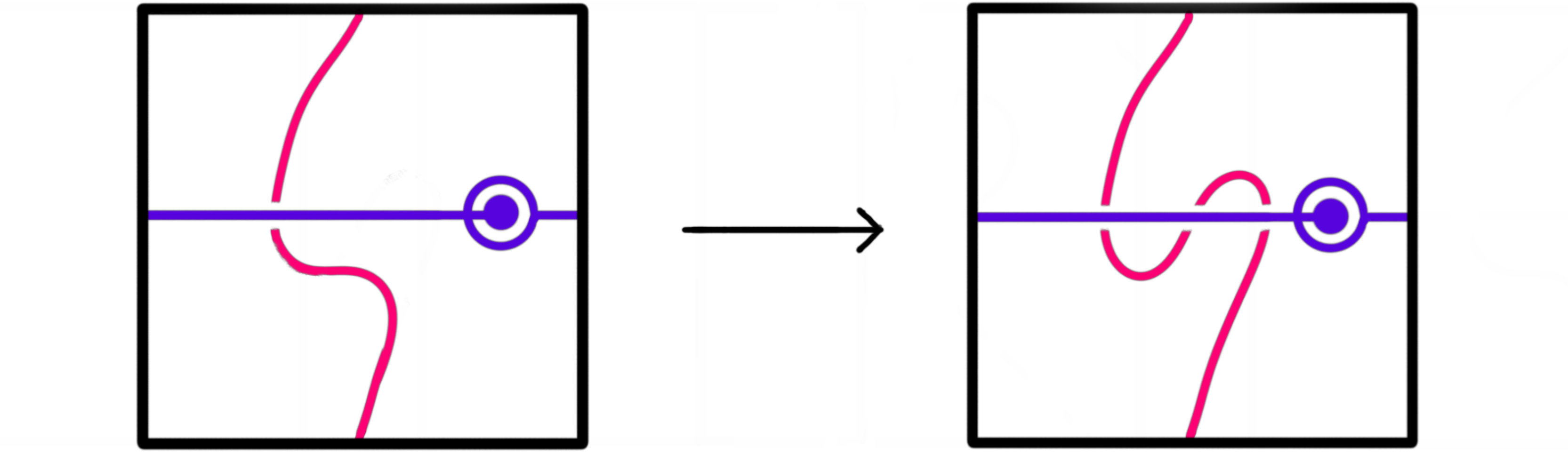}
        \caption{}
        \label{fig:r_moves_tridia_2}
    \end{subfigure}
    \vskip\baselineskip
    \begin{subfigure}[b]{0.65\textwidth}
        \includegraphics[width=\textwidth]{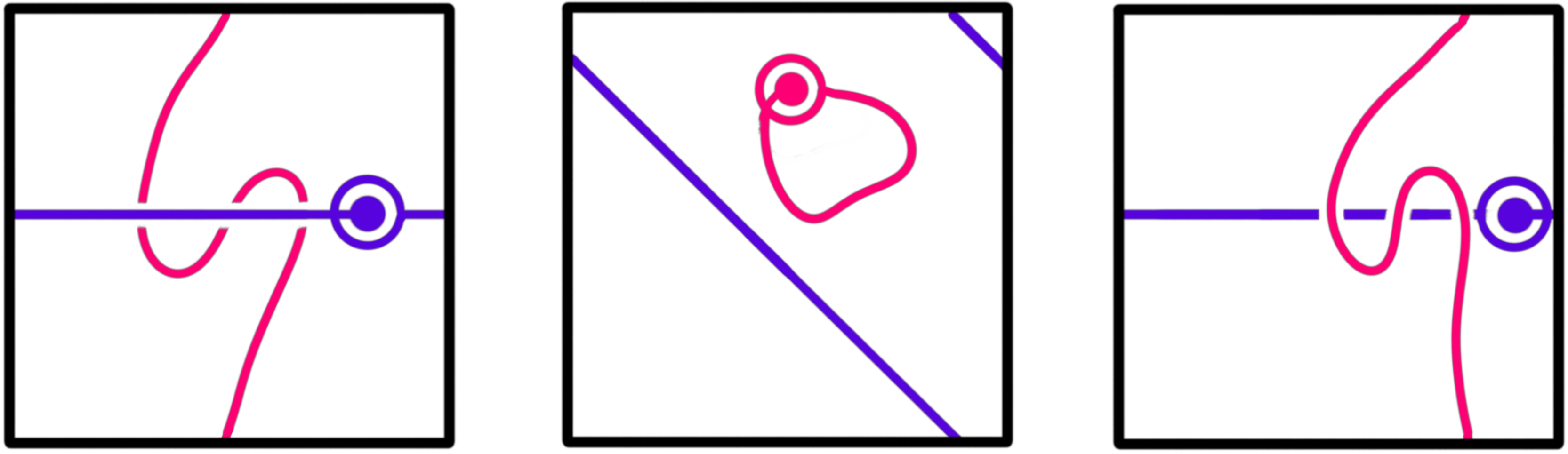}
        \caption{}
        \label{fig:r_moves_tridia_3}
    \end{subfigure}
    \vskip\baselineskip
    \begin{subfigure}[b]{0.65\textwidth}
        \includegraphics[width=\textwidth]{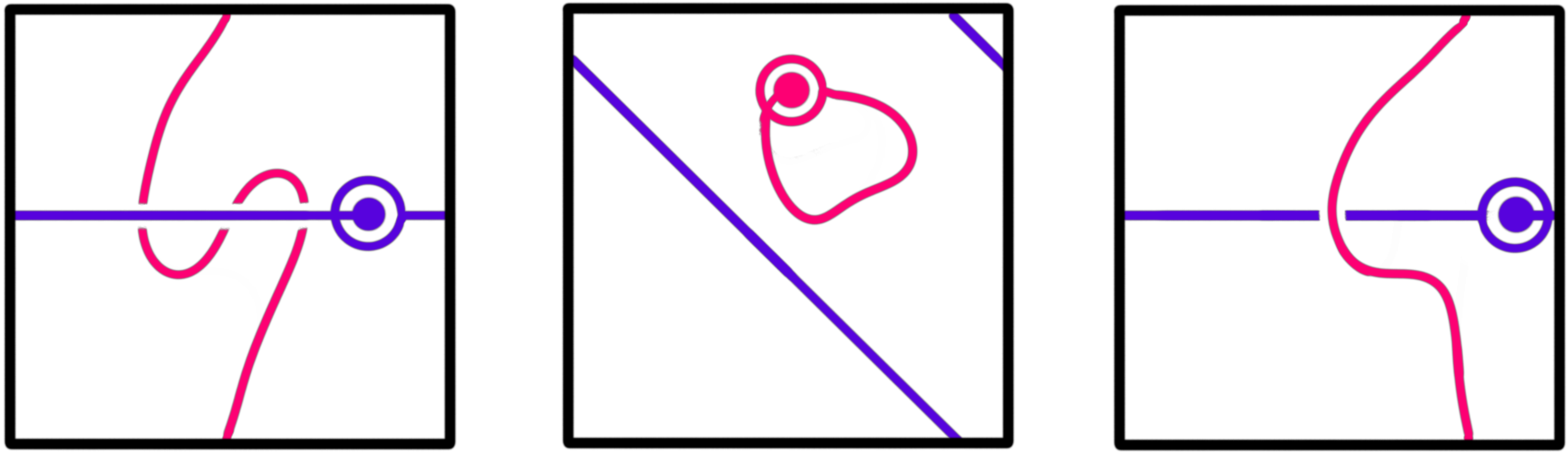}
        \caption{}
        \label{fig:r_moves_tridia_4}
    \end{subfigure}
    \caption{The behaviour of tridiagrams and $R$-moves: (a) A tridiagram where the diagrams represent the front, top and right faces of the associated unit cell. (b) One $R_2$ move is performed on the first diagram of the tridiagram of (a). (c) One possibility of a tridiagram resulting from the $R_2$ move. (d) Three diagrams that do not constitute a tridiagram.}
    \label{fig:r_moves_tridia}
\end{figure}
\end{remark}

Aside from ambient isotopies of the 3-torus, the other notions needed for the equivalence of TP tangles, such as the torus twists or the $CM^3$-equivalence for example, also need to be described at the diagrammatic level. The following definitions are given to account for those notions of equivalence.

\begin{definition}\label{def:2_torus_twists}
    A \textit{$2$-torus twist}, or \textit{torus twist} when the context leaves no ambiguity, is one of the following two homeomorphisms:
    \[
\begin{array}{cccc}
    \psi_{1} : & \mathbb{S}^1 \times \mathbb{S}^1  & \longrightarrow & \mathbb{S}^1 \times \mathbb{S}^1 \\
    & \left(z_1,z_2\right) & \longmapsto & \left(z_1z_2,z_2\right) \\
\end{array},
\]
 \[
\begin{array}{cccc}
    \psi_{2} : & \mathbb{S}^1 \times \mathbb{S}^1  & \longrightarrow & \mathbb{S}^1 \times \mathbb{S}^1 \\
    & \left(z_1,z_2\right) & \longmapsto & \left(z_1,z_1z_2\right) \\
\end{array}.
\]
\end{definition}

\begin{definition}\label{def:t2_equivalence}
    Two diagrams $D$ and $D'$ are \textit{$T^2$-equivalent} if they are connected by a finite composition of $2$-torus twists (and their inverses as well as the identity map). Furthermore, two tridiagrams $\left\lbrace D_i \right\rbrace_{i=1,2,3}$ and $\left\lbrace D_i' \right\rbrace_{i=1,2,3}$ are said to be \textit{$T^2$-equivalent} if there is a finite sequence of tridiagrams $\left( \lbrace D_{i,k} \rbrace_{i=1,2,3} \right)_{k=0,\dots,n}$, for which we have 
    \[  D_{i,0}  =   D_{i} , \quad  D_{i,n}  =  D_{i}^{\prime}, \quad i=1,2,3\]
    and 
    \[\forall k = 1,\dots,n, \exists i,j = 1,2,3: \quad \text{ $D_{i,k-1}$ and $D_{j,k}$ are $T^2$-equivalent}.\]
    The equivalence is denoted by $\sim_T$ for both cases.
\end{definition}

\begin{remark}
As suggested by Definition \ref{def:t2_equivalence}, when two tridiagrams are connected by a torus twist, it does not imply that all three diagrams of the tridiagrams are connected by a torus twist. See the example of Figure \ref{fig:pi_star_torus_twist_tridiagram}. A tridiagram is displayed in Figure \ref{fig:pi_star_tridiagram_2nd}, where the diagrams are respectively projected from the front, the top and the right faces of the corresponding unit cell. On the front diagram, a torus twist is performed, as shown in Figure \ref{fig:pi_star_2_torus_twist}. The resulting tridiagram is displayed in Figure \ref{fig:pi_star_2_torus_twist_tridia}, where the last two diagrams are not connected by a torus twist to the diagrams of Figure \ref{fig:pi_star_tridiagram_2nd}.
    \begin{figure}[htbp]
    \centering
    \begin{subfigure}[b]{\textwidth}
    \centering
        \includegraphics[width=0.75\textwidth]{pi_star_tridia.pdf}
        \caption{}
        \label{fig:pi_star_tridiagram_2nd}
    \end{subfigure}
    \vskip\baselineskip
    \begin{subfigure}[b]{\textwidth}
    \centering
        \includegraphics[width=0.75\textwidth]{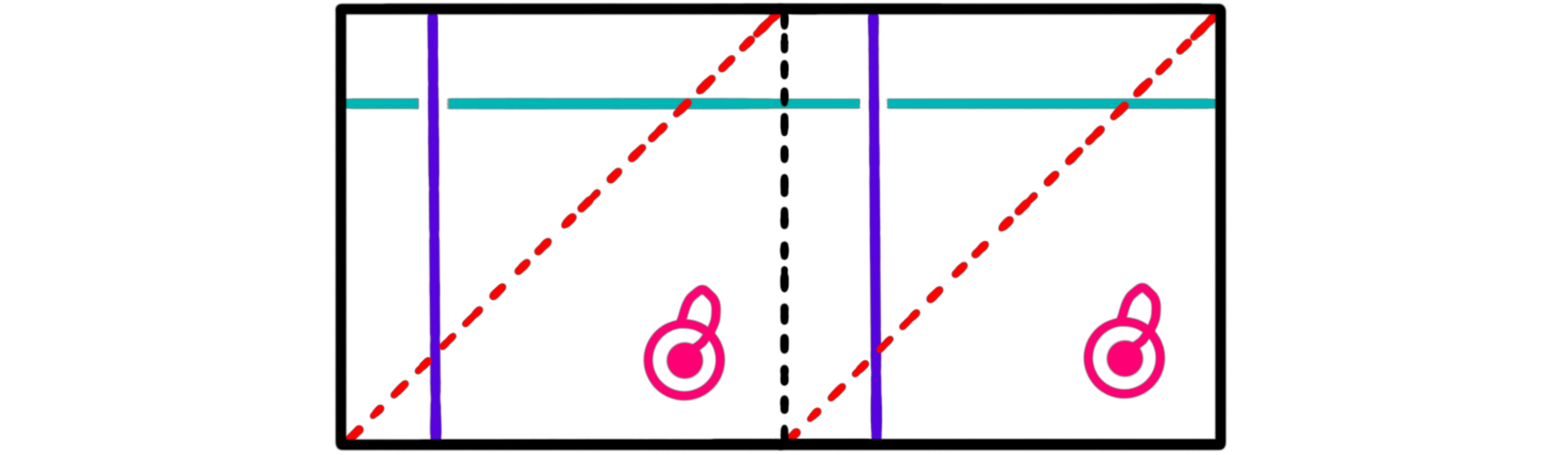}
        \caption{}
        \label{fig:pi_star_2_torus_twist}
    \end{subfigure}
    \vskip\baselineskip
    \begin{subfigure}[b]{\textwidth}
    \centering
        \includegraphics[width=0.75\textwidth]{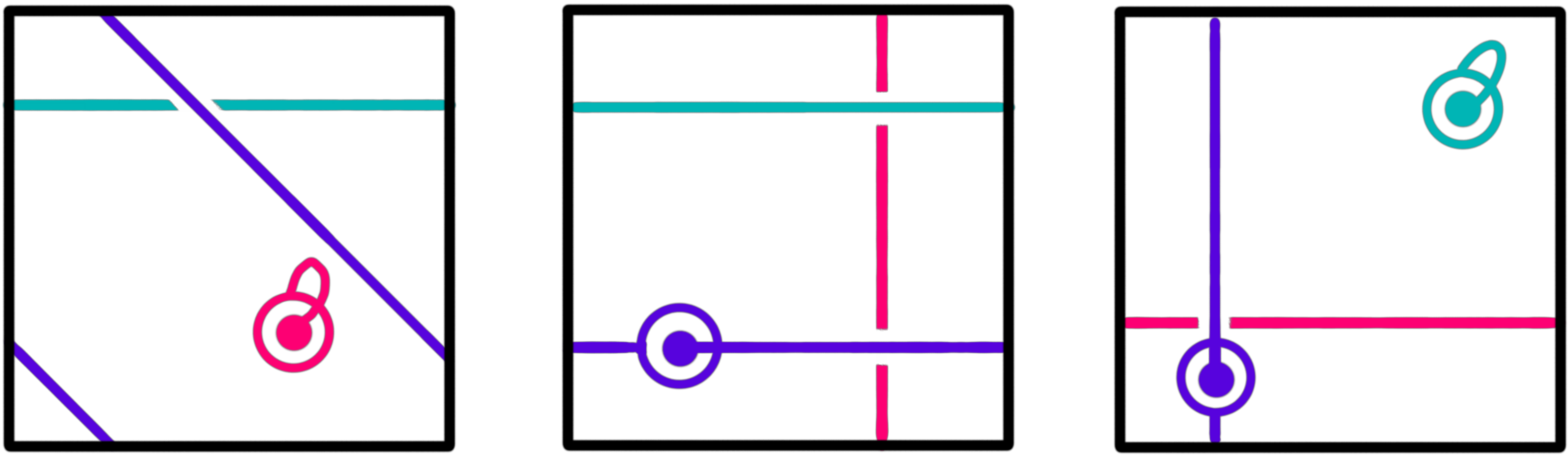}
        \caption{}
        \label{fig:pi_star_2_torus_twist_tridia}
    \end{subfigure}
    \caption{Visualisation of a $2$-torus twist in a tridiagram: (a) A tridiagram obtained from projections from the front, top and right faces of a unit cell. (b) A $2$-torus twist performed on the front diagram of (a). (c) The resulting tridiagram, where the last two diagrams are not connected by a torus twist to the last two diagrams of (a).}
    \label{fig:pi_star_torus_twist_tridiagram}
    \end{figure}

    It is for this reason that we actually defined tridiagrams. Indeed, in the usual theory of knots, one wants to prove that a single diagram is enough to encompass all information and transformations of a knot. However, doing so with 3-periodic tangles is very inconvenient. Take the example of the last diagrams of Figure \ref{fig:pi_star_tridiagram_2nd} and Figure \ref{fig:pi_star_2_torus_twist_tridia}, and suppose that one wants to work only with those diagrams. One can understand that the unit cells they represent are two unit cells belonging to the same lattice that differ only by a $3$-torus twist. However, the two diagrams differ by the appearance of an $N$-point in the diagram of Figure \ref{fig:pi_star_2_torus_twist_tridia}. The question, therefore, is to know whether there is a general rule about transforming the diagram of Figure \ref{fig:pi_star_tridiagram_2nd} to the one of Figure \ref{fig:pi_star_2_torus_twist_tridia}, given by some `generalised Reidemeister' move. But a priori, there is no such rule, or at least, it is very incovenient to work with such a rule.
\end{remark}

\begin{definition}\label{def:cover_2}
    Two diagrams $D$ and $D'$ are \textit{$CM^2$-equivalent} if there is a sequence of diagrams $(D_k)_{k=0,\dots,n}$ where $D_0 = D$, $D_n = D'$, and for all $k=1,\dots,n$, either $D_{k-1}$ is a cover $D_{k}$ or $D_{k}$ is a cover $D_{k-1}$, where the covering map is of type 
    \[
\begin{array}{ccc}
     \mathbb{S}^1 \times \mathbb{S}^1 & \longrightarrow & \mathbb{S}^1 \times \mathbb{S}^1 \\
      (z_1,z_2) & \longmapsto & \left({z_1}^{h_{11}}{z_2}^{h_{21}},{z_2}^{h_{22}}\right)
\end{array}, z_i \in \mathbb{C}, \|z_i\| = 1,
\]
where $0 \leqslant h_{ij} < h_{ii}$ for $j < i$. This is a covering map whith $h_{11}h_{22}$ sheets for each point of $\mathbb{T}^2$.

Furthermore, two tridiagrams $\left\lbrace D_i \right\rbrace_{i=1,2,3}$ and $\left\lbrace D_i' \right\rbrace_{i=1,2,3}$ are \textit{$CM^2$-equivalent} if there is a finite sequence of tridiagrams $\left( \lbrace D_{i,k} \rbrace_{i=1,2,3} \right)_{k=0,\dots,n}$, for which we have 
    \[  D_{i,0}  =   D_{i} , \quad  D_{i,n}  =  D_{i}^{\prime}, \quad i=1,2,3\]
    and 
    \[\forall k = 1,\dots,n, \exists i,j = 1,2,3: \quad \text{ $D_{i,k-1}$ and $D_{j,k}$ are $CM^2$-equivalent}.\]
    The equivalence is denoted by $\sim_{CM^2}$ for both cases.
\end{definition}

\begin{remark}
The fact that two tridiagrams are $CM^2$-equivalent does not imply that all the diagrams of the tridiagrams are covers of one another. See the example of Figure \ref{fig:pi_star_torus_m2_equi}. In Figure \ref{fig:pi_star_m2_1}, we display a cover of the first diagram of the tridiagram of Figure \ref{fig:pi_star_m2_0}. Clearly, the last diagram of Figure \ref{fig:pi_star_m2_2} is not a cover of the last diagram of Figure \ref{fig:pi_star_m2_0}.
    \begin{figure}[htbp]
    \centering
    \begin{subfigure}[b]{\textwidth}
    \centering
        \includegraphics[width=0.75\textwidth]{pi_star_2_torus_twist_tridia.pdf}
        \caption{}
        \label{fig:pi_star_m2_0}
    \end{subfigure}
    \vskip\baselineskip
    \begin{subfigure}[b]{\textwidth}
    \centering
        \includegraphics[width=0.75\textwidth]{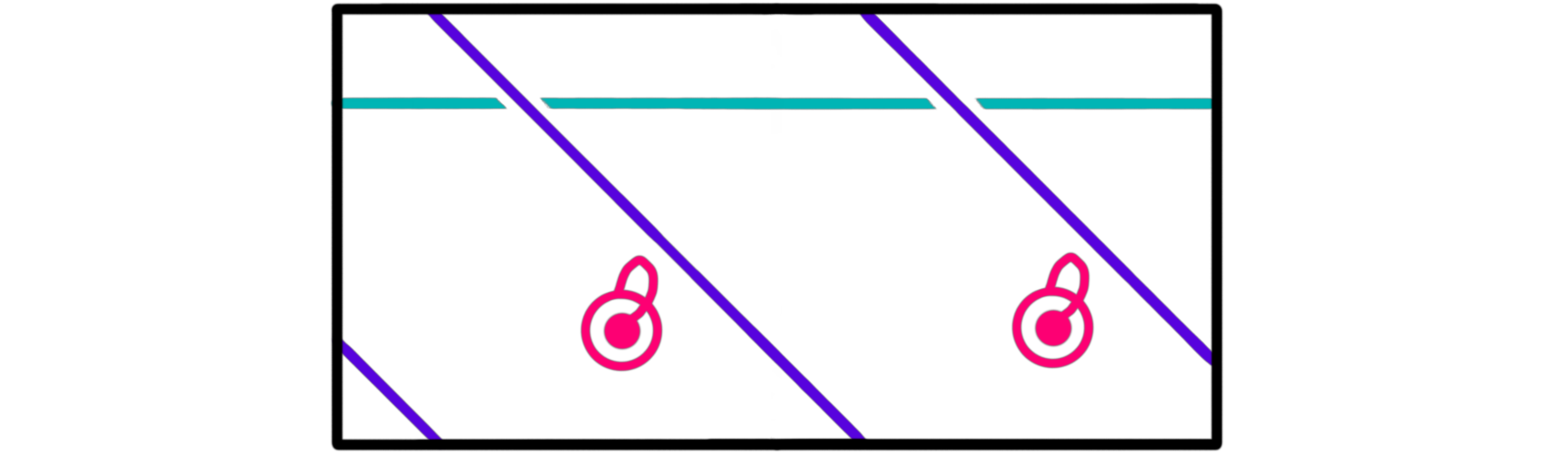}
        \caption{}
        \label{fig:pi_star_m2_1}
    \end{subfigure}
    \vskip\baselineskip
    \begin{subfigure}[b]{\textwidth}
    \centering
        \includegraphics[width=0.75\textwidth]{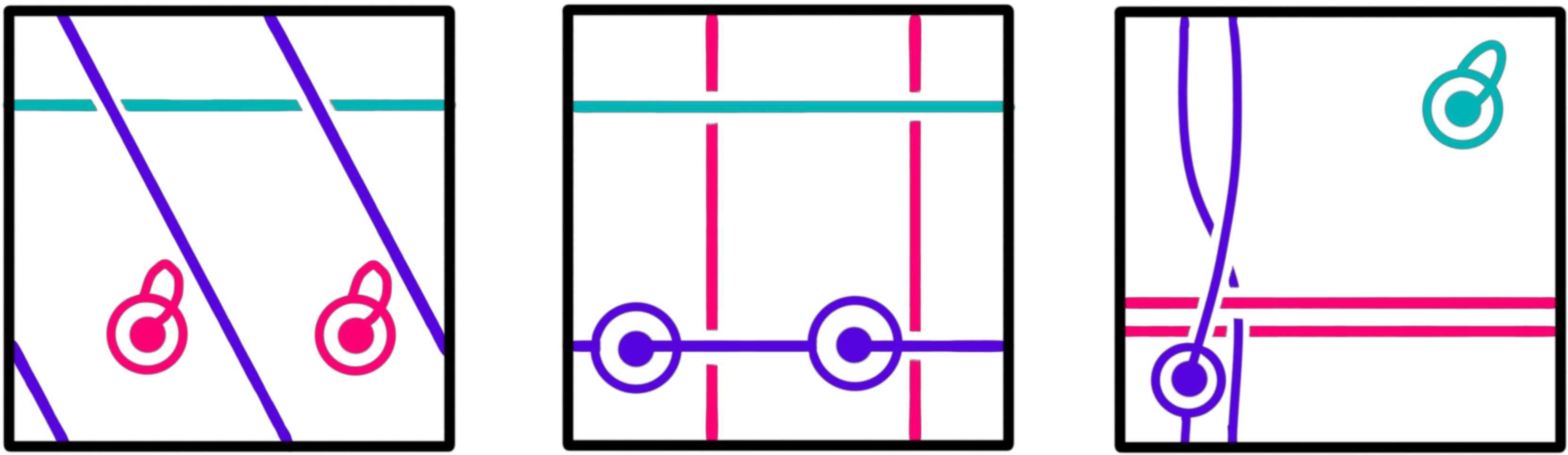}
        \caption{}
        \label{fig:pi_star_m2_2}
    \end{subfigure}
    \caption{Visualisation of the $CM^2$-equivalence of tridiagrams: (a) A tridiagram obtained from projections from the front, top and right faces of a unit cell. (b) A cover of the first diagram of (a). (c) A tridiagram corresponding to the cover diagram, where the last diagram is not a cover of the last diagram of (a).}
    \label{fig:pi_star_torus_m2_equi}
    \end{figure}

    \noindent Furthermore, a priori there is no way to connect the last diagrams of Figure \ref{fig:pi_star_m2_0} and Figure \ref{fig:pi_star_m2_2} by some `generalised Reidemeister' move.
\end{remark}

\begin{definition}
    Consider a diagram $D$ and an embedding $\Gamma$ of $D$ as a link in $\mathbb{T}^3 \cong \mathbb{S}^1_1 \times \mathbb{S}^1_2 \times  \mathbb{S}^1_3$. Suppose that $\pi : \mathbb{S}^1_1 \times \mathbb{S}^1_2 \times  (0,1)_N\longrightarrow \mathbb{S}^1_1 \times \mathbb{S}^1_2$ is the projection that yields $D$ from $\Gamma$. By projecting $\Gamma$ along $\mathbb{S}^1_3$ in the negative orientation, that is, from $1$ to $0$, one can obtain another diagram from the back face of the cube representing the $3$-torus. We say that two diagrams $D$ and $D'$ are \textit{$V$-equivalent} if $D = D'$ or if $D'$ is obtained from an embedding of $D$, where the projection is along the opposite direction of that of $D$. Furthermore, two tridiagrams $\left\lbrace D_i \right\rbrace_{i=1,2,3}$ and $\left\lbrace D_i' \right\rbrace_{i=1,2,3}$ are \textit{$V$-equivalent} if two of their respective diagrams are $V$-equivalent.

    \noindent The equivalence is denoted by $\sim_V$ for both cases.
\end{definition}

We wrap up all these notions of equivalence into one in the following.

\begin{definition}\label{def:equi_tridiagrams}
    Two tridiagrams $\left\lbrace D_i \right\rbrace_{i=1,2,3}$ and $\left\lbrace D_i' \right\rbrace_{i=1,2,3}$ are \textit{equivalent}, if there is a finite sequence of tridiagrams $\left( \lbrace D_{i,k} \rbrace_{i=1,2,3} \right)_{k=0,\dots,n}$, for which we have:
    \begin{itemize}
        \item[-]$ D_{i,0}  =   D_{i}$ and $ D_{i,n}  =  D_{i}^{\prime}, \quad i=1,2,3$,
        \item[-] $\forall k = 1,\dots,n$, $\lbrace D_{i,k-1} \rbrace_{i=1,2,3}$ and $\lbrace D_{i,k} \rbrace_{i=1,2,3}$ are $R$-equivalent, or $CM^2$-equivalent, or $T^2$-equivalent or $V$-equivalent. 
    \end{itemize}
     The equivalence is denoted by $\sim$.
\end{definition}

\subsection{Generalised Reidemeister theorem}

Before formulating our generalised Reidemeister theorem for 3-periodic tangles, we state several propositions relating notions of equivalence of links embedded in the $3$-torus to notions of equivalence of diagrams and tridiagrams.\\

We start by proving a proposition, which actually is one implication of the Reidemeister theorem for links in the 3-torus, using the diagrams as defined in this paper. Other generalisations of the Reidemeister theorem for (oriented) links in the 3-torus have been proven in the literature, such as in \cite{LAMBROPOULOU199795} where one makes use of mixed link diagrams, or in \cite{vuong2023fundamental} in a similar construction to the diagrams defined in this paper, however we give here a complete proof.

\begin{proposition}\label{prop:ambi_iso_links_imp_equi_dia}
    Two ambient isotopic links in the $3$-torus $\Gamma$ and $\Gamma'$ possess equivalent diagrams.
\end{proposition}

To prove this proposition, we extend to our diagrams the already known proof for usual knots in $\mathbb{R}^3$ using singularity theory \cite{grid_homology_appendix_B, Roseman2004ElementaryMF}. Without loss of generality, we can assume that $\Gamma$ and $\Gamma'$ are links in $\mathbb{T}^3$ with one component. The proof of the general case is similar.

We are going to prove that supposing the two 1-component links $\Gamma$ and $\Gamma'$ are ambient isotopic, the isotopy $ [0,1] \times \mathbb{S}^1 \longrightarrow \mathbb{T}^3$ between them gives rise to a 1-parameter family of diagrams which are connected by planar isotopies and $R$-moves.

Consider the following definition, in accordance to \cite{Hirsch1976_chap8}.
\begin{definition}
    Let $V$ and $M$ be manifolds, and $F: [0,1]  \times V \longrightarrow M$ an isotopy. We define the \textit{track} of $F$ as the map
    \[
    \begin{array}{cccc}
       \hat{F} : & [0,1]  \times V  & \longrightarrow & [0,1]  \times M \\
         & (t,v) & \longmapsto & \left(t,F(t,v)\right)
    \end{array}.
    \]
\end{definition}

One can compose the track of the isotopy between $\Gamma$ and $\Gamma'$ with the projection of $[0,1] \times \mathbb{T}^3$ onto $[0,1] \times \mathbb{T}^2$ to obtain a map from the 2-manifold $[0,1]\times \mathbb{S}^1$ to the 3-manifold $[0,1]\times \mathbb{T}^2$. By applying a small perturbation to this composite map, it can be put into general position for almost all $t \in [0,1]$. The remaining `non-generic diagrams' in $\lbrace t \rbrace \times \mathbb{T}^2$ for finitely many $t$'s will be dealt with using the nine $R$-moves.\\

Some cases of non-regularity of a projection in $\lbrace t \rbrace \times \mathbb{T}^2$ are due to the singularities of the map $[0,1]\times \mathbb{S}^1 \longrightarrow [0,1]\times \mathbb{T}^2$. To understand those singularities, a key local result we will use in our proof is the following theorem by Whitney.

\begin{theorem} \label{thm:whitney}
    (Whitney, \cite{Whitney1992}). Let $W$ be a smooth 2-manifold and $Y$ a smooth 3-manifold. Any smooth map $g_0 : W \longrightarrow Y$ can be approximated arbitrarily closely (in the $\mathcal{C}^2$ topology) by a smooth map $g : W \longrightarrow Y$ with the following property. Around each point $p \in W$, there are local coordinates $(x, y)$ so that $p$ corresponds to $(0, 0)$, and there are local coordinates $(u, v, z)$ around $g(p) \in Y$, so that $(0, 0, 0)$ corresponds to $g(p)$, with respect to which the function $g$ has the form $(x, y) \longmapsto (x, y, 0)$ or $(x, y) \longmapsto (x^2, xy, y)$.
\end{theorem}

\begin{figure}[hbtp]
\centering
\includegraphics[width=5cm]{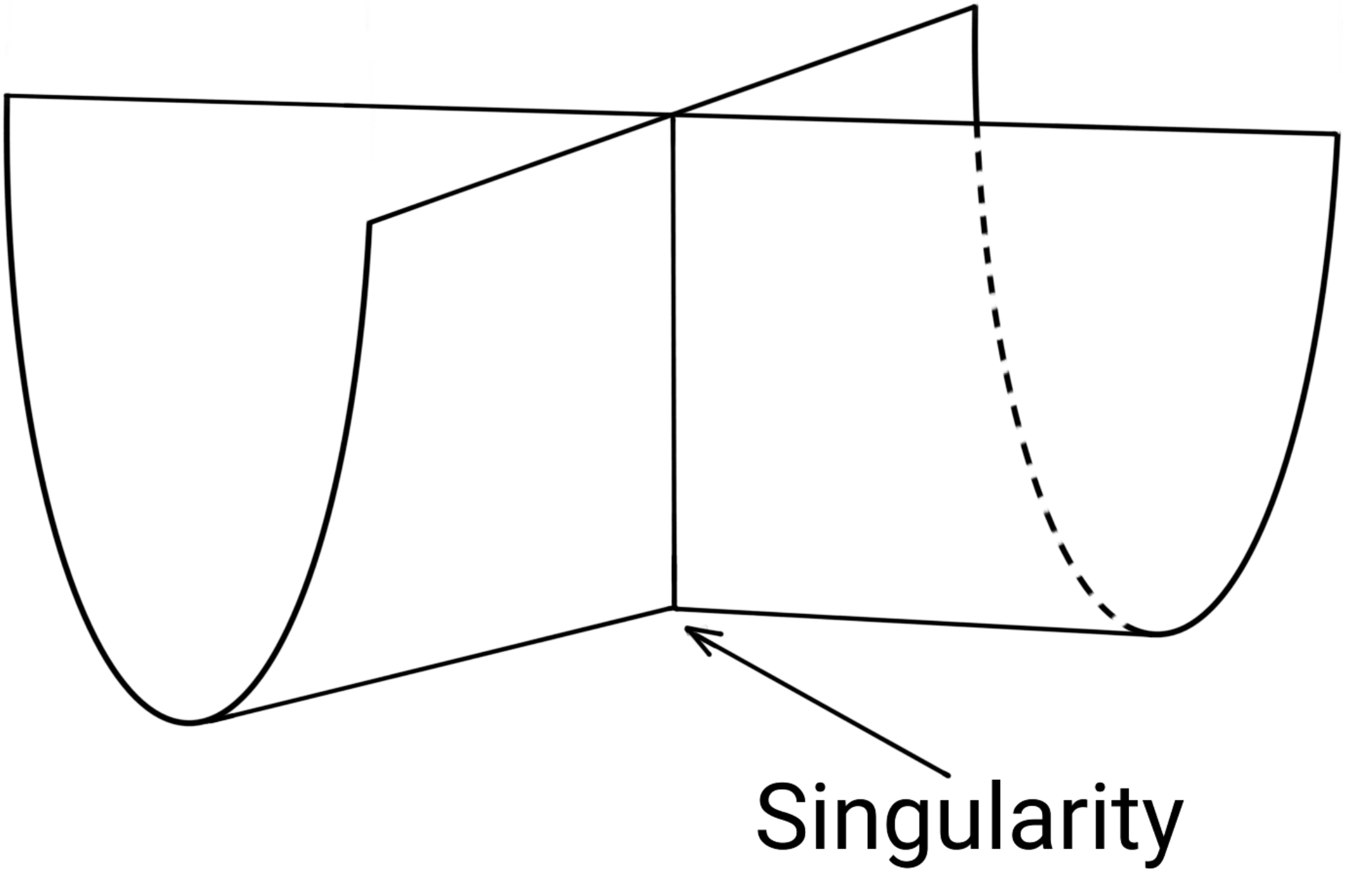}
\caption{The Whitney umbrella. It is the image of the map $(x, y) \longmapsto (x^2, xy, y)$.}
\label{fig:whitney}
\end{figure}

Whitney's theorem can be understood as follows. Consider the Jacobian $J_{g(p)}: T_pW \longrightarrow T_{g(p)}Y$ of the map $g$ at $p \in W$. On the one hand, it is clear that points of the first kind are those for which the Jacobian is injective, and so, at those points $g$ is an immersion. On the other hand, points of the second kind are the singular points of $g$, where the Jacobian has a 1-dimensional kernel. Note that for a generic choice of $g$, the rank of $J_{g(p)}$ is non-zero for all $p \in W$. Thus, Whitney’s theorem gives a canonical form for the neighbourhood of the singular points of $g$, depicted in Figure \ref{fig:whitney}. Such a singularity is called a \textit{Whitney umbrella}. A proof of Theorem \ref{thm:whitney} is given in \cite{grid_homology_appendix_B}.\\

\begin{proof}
\textbf{(Proof of Proposition \ref{prop:ambi_iso_links_imp_equi_dia}).}
Suppose that $D$ is the diagram of the smoothly embedded link with one component $\Gamma$. Assuming that $\Gamma'$ is ambient isotopic to $\Gamma$, the same projection map yielding $D$ from $\Gamma$, yields a diagram $D'$ from $\Gamma'$. Regard the track of the isotopy between the links $f: [0,1] \times \mathbb{S}^1 \longrightarrow [0,1] \times \mathbb{T}^3$ as a smooth embedding of the surface-with-boundary $[0,1] \times \mathbb{S}^1$ with the following properties:
\begin{enumerate}
    \item $f([0, 1] \times \mathbb{S}^1) \cap (\lbrace0\rbrace \times \mathbb{T}^3) = \lbrace0\rbrace \times \Gamma$,
    \item $f([0, 1] \times \mathbb{S}^1) \cap (\lbrace1\rbrace \times \mathbb{T}^3) = \lbrace1\rbrace \times \Gamma'$ and
    \item \label{item:property_3} the intersection $f([0, 1] \times \mathbb{S}^1) \cap (\lbrace t\rbrace \times \mathbb{T}^3)$ is transverse for all $t \in [0, 1]$.
\end{enumerate}

Consider the projection map $pr_1 : [0, 1] \times \mathbb{T}^3 \longrightarrow [0, 1]$. Property (\ref{item:property_3}) is actually a reformulation of the fact that $pr_1 \circ f : [0, 1] \times \mathbb{S}^1 \longrightarrow [0, 1]$ has no critical points. Indeed, if $y \in f([0, 1] \times \mathbb{S}^1) \cap (\lbrace t\rbrace \times \mathbb{T}^3)$, as $T_y \left(\lbrace t\rbrace \times \mathbb{T}^3 \right)$ has dimension 3, Property (\ref{item:property_3}) is true if and only if $T_y f([0, 1] \times \mathbb{S}^1)$ is of dimension greater than or equal to 1 and is not a subspace of $T_y \left(\lbrace t\rbrace \times \mathbb{T}^3\right)$. This is equivalent to $T_y f([0, 1] \times \mathbb{S}^1) \cap T_{pr_1(y)}[0,1] \neq \lbrace 0 \rbrace$. This occurs if and only if $T_{pr_1(y)} ((pr_1 \circ f)([0, 1] \times \mathbb{S}^1)) \neq \lbrace 0 \rbrace$, which is the characterisation of $pr_1 \circ f$ having no critical points, meaning that $pr_1 \circ f$ is a submersion. Now recall that being a submersion is an open property \cite{guillemin_differential_topology_chap1} which means that a small perturbation of $f$ does not have a critical point.

Furthermore, pictorially, consider $t \in [0, 1]$ of $[0, 1] \times \mathbb{T}^3$ as the `variation of time'. By doing so, $\lbrace t \rbrace \times \mathbb{T}^3$ has a front face indexed by $t$ which we call \textit{$t$-front face} $ \lbrace t \rbrace \times \mathbb{S}^1_1 \times \mathbb{S}^1_2 \times \lbrace 0+\varepsilon \rbrace$, a back face indexed by $t$ which we call \textit{$t$-back face} $ \lbrace t \rbrace \times \mathbb{S}^1_1 \times \mathbb{S}^1_2 \times \lbrace 1-\varepsilon \rbrace$, $\varepsilon > 0$, and a \textit{$t$-$N$-face} $\lbrace t \rbrace \times \mathbb{S}^1_1 \times \mathbb{S}^1_2 \times \lbrace N \rbrace$ extending Definition \ref{def:front_back_faces} for all $t \in [0,1]$.

We now extend the projection map $\pi : \mathbb{T}^3 \longrightarrow  \mathbb{T}^2 $ defining the diagram to the map $P = \mathrm{id} \times \pi : [0, 1] \times \mathbb{T}^3  \longrightarrow [0, 1] \times \mathbb{T}^2$, and compose the embedding $f$ of the annulus $[0, 1] \times \mathbb{S}^1$ with $P$ to get a map $\phi = P \circ f : [0, 1] \times \mathbb{S}^1 \longrightarrow [0, 1] \times \mathbb{T}^2$. Proposition \ref{prop:ambi_iso_links_imp_equi_dia} is proved by applying Theorem \ref{thm:whitney} to this map $\phi$, as follows.

Since $P$ is a submersion, $\phi$ can be put into general position by slightly perturbing the map $f$, so that it remains an isotopy between $\Gamma$ and $\Gamma'$. This is possible thanks to the three properties of $f$ listed above as explained in \cite{Roseman2004ElementaryMF}. Theorem \ref{thm:whitney} shows that there are finitely many points $\mathcal{W} \subset [0, 1] \times \mathbb{S}^1$ with a Whitney umbrella singularity, away from which the map $\phi$ is an immersion. By a further general position argument, we can assume that $\phi$ has only finitely many triple points (that is, points in $[0, 1] \times \mathbb{T}^2$ with three preimages), and a union of 1-dimensional submanifolds $\mathcal{D} \subset [0, 1] \times \mathbb{T}^2$ of double points. The closure of the set of double points includes the set of triple points and the set of Whitney umbrella singularities. Its boundary also includes the double points $\mathcal{D} \cap ( \lbrace 0 \rbrace \times \mathbb{T}^2)$, $\mathcal{D} \cap ( \lbrace 1 \rbrace \times \mathbb{T}^2)$ of the two original diagrams $D$ and $D'$.

Regard the intersections $\phi([0, 1] \times \mathbb{S}^1) \cap (\lbrace t\rbrace \times \mathbb{T}^2)$ as a 1-parameter family of diagrams. Again, by general position arguments, there are finitely many special $t \in [0, 1]$ where these diagrams are not generic, and where exactly one of the following occurs:
\begin{enumerate}
    \item \label{item:special_t_1} $\phi^{-1}(\lbrace t \rbrace \times \mathbb{T}^2)$ contains a Whitney umbrella singularity,
    \item \label{item:special_t_2} $\phi([0, 1] \times \mathbb{S}^1) \cap (\lbrace t\rbrace \times \mathbb{T}^2)$ contains a triple point,
    \item \label{item:special_t_3} $\lbrace t \rbrace \times \mathbb{T}^2$ is tangent to $\mathcal{D}$,
    \item \label{item:special_t_4} $\phi([0, 1] \times \mathbb{S}^1) \cap (\lbrace t\rbrace \times \mathbb{T}^2)$ contains a double point for which the embedding in $\lbrace t \rbrace \times \mathbb{T}^3$ of one of the two elements of its preimage is lying in the $t$-$N$-face,
    \item \label{item:special_t_5} $\phi([0, 1] \times \mathbb{S}^1) \cap (\lbrace t\rbrace \times \mathbb{T}^2)$ contains a point such that the $t$-$N$-face is tangent to the embedding of the preimage of that point in $\lbrace t \rbrace \times \mathbb{T}^3$,
    \item \label{item:special_t_6} $\lbrace t \rbrace \times l_i$, $i=1,2$ is tangent to $\phi([0, 1] \times \mathbb{S}^1) \cap (\lbrace t \rbrace \times \mathbb{T}^2)$.
    \item \label{item:special_t_7} $\mathcal{D} \cap (\lbrace t \rbrace \times l_i)$, $i=1,2$ is not empty,
    \item \label{item:special_t_8} $\phi([0, 1] \times \mathbb{S}^1) \cap \left(\lbrace t \rbrace \times l_1 \right) \cap  \left(\lbrace t \rbrace \times l_2 \right) $ is not empty,
    \item \label{item:special_t_9} a point of the embedded curve in $\lbrace t \rbrace \times \mathbb{T}^3$ belongs to the $t$-$N$-face, and is projected onto $\lbrace t \rbrace \times l_i$, $i=1,2$.
\end{enumerate}

Consider $0 \leqslant t_1 < t_2 \leqslant 1$, and suppose that there are no special values of $t \in [t_1, t_2]$. In this case, as $\phi$ is an immersion and thus locally an embedding, the projections at $t_1$ and $t_2$ are connected by planar isotopies.

Assume next that the interval $[t_1, t_2]$ contains a single special value. Suppose that the special value corresponds to a Whitney umbrella singularity as in (\ref{item:special_t_1}). Furthermore, assume that at that value $t$, $\lbrace t \rbrace \times \mathbb{T}^2$ is transverse to the 1-dimensional image of the Jacobian. Then the diagrams $\phi([0, 1] \times \mathbb{S}^1) \cap \lbrace t_1 \rbrace \times \mathbb{T}^2$ and $\phi([0, 1] \times \mathbb{S}^1) \cap \lbrace t_2\rbrace \times \mathbb{T}^2$ differ by a single $R_1$ move as seen in Figure \ref{fig:proof_R1}. Indeed the Whitney umbrella singularity occurs when a crossing is being taken out or added to the family of diagrams.

Next, as $\lbrace t \rbrace \times \mathbb{T}^2$ passes through a point of tangency with $\mathcal{D}$ as given by (\ref{item:special_t_2}), the projection undergoes an $R_2$ move as shown in Figure \ref{fig:proof_R2}. $\lbrace t \rbrace \times \mathbb{T}^2$ is tangent to $\mathcal{D}$ exactly when two crossings are deleted or added to the family of diagrams. 

In the case (\ref{item:special_t_3}), a triple point is locally modelled on three intersecting planes as in Figure \ref{fig:proof_R3}. The diagrams $\phi([0, 1] \times \mathbb{S}^1) \cap \lbrace t_1 \rbrace \times \mathbb{T}^2$ and $\phi([0, 1] \times \mathbb{S}^1) \cap \lbrace t_2\rbrace \times \mathbb{T}^2$ differ by an $R_3$ move. A triple point occurs when one strand passes through the crossing of other two strands in the family of diagrams.

\begin{figure}[hbtp]
    \centering
    \begin{subfigure}[b]{4.5cm}
        \includegraphics[width=\textwidth]{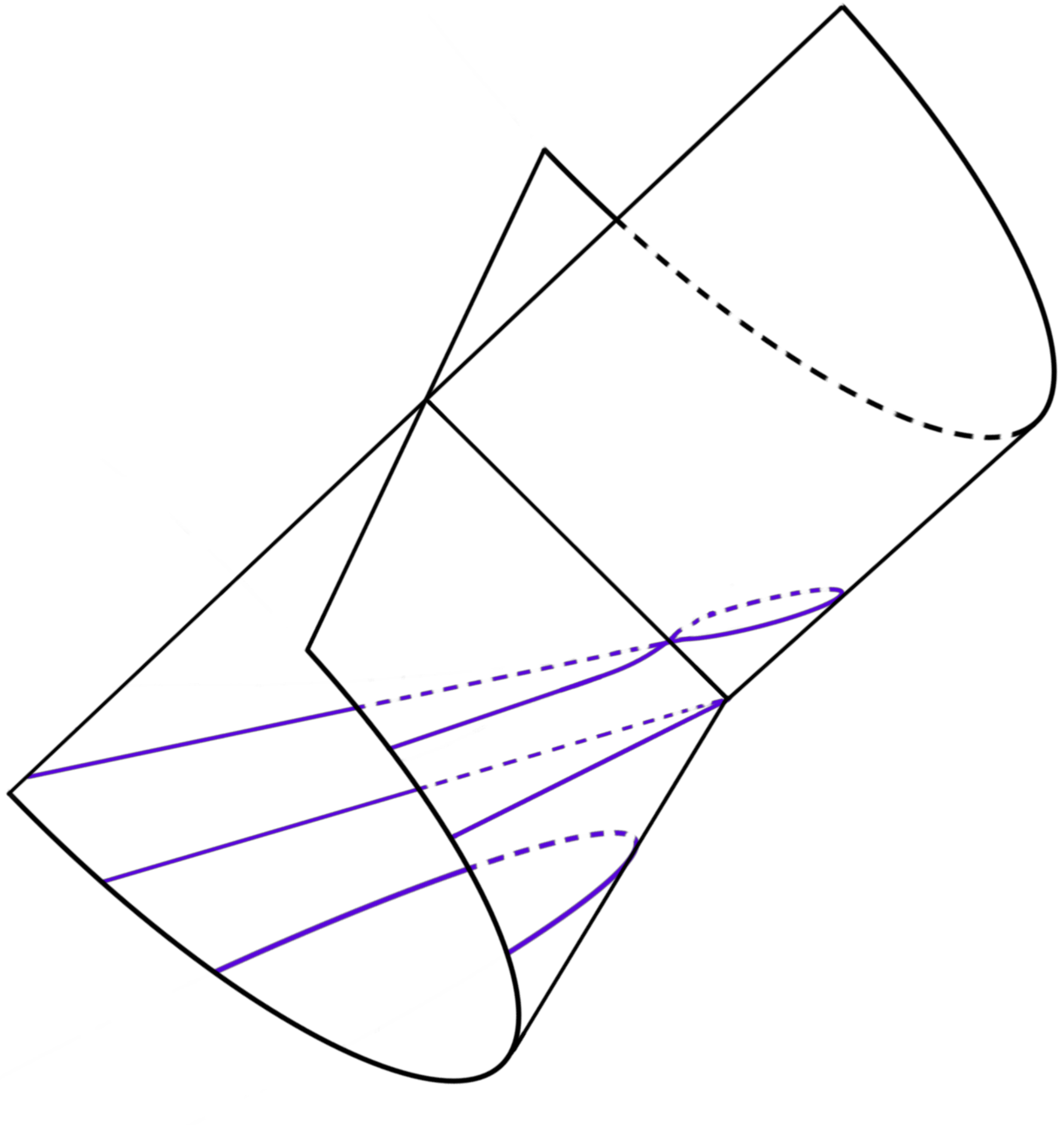}
        \caption{}
        \label{fig:proof_R1}
    \end{subfigure}
    \hspace{0.5cm}
    \begin{subfigure}[b]{4cm}        
        \includegraphics[width=\textwidth]{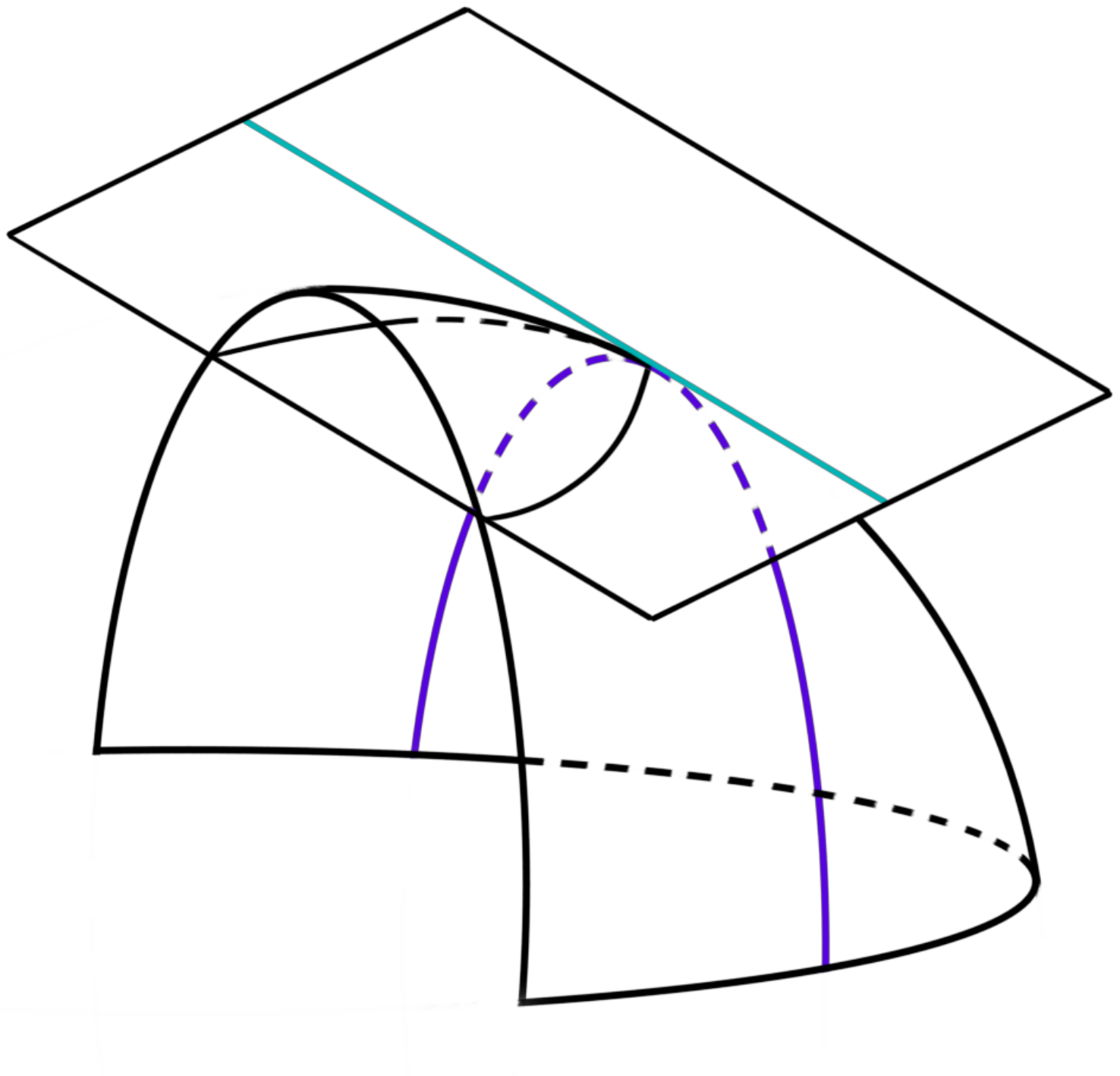}
        \caption{}
        \label{fig:proof_R2}
    \end{subfigure}
    \hspace{0.5cm}
    \begin{subfigure}[b]{3.25cm}
        \includegraphics[width=\textwidth]{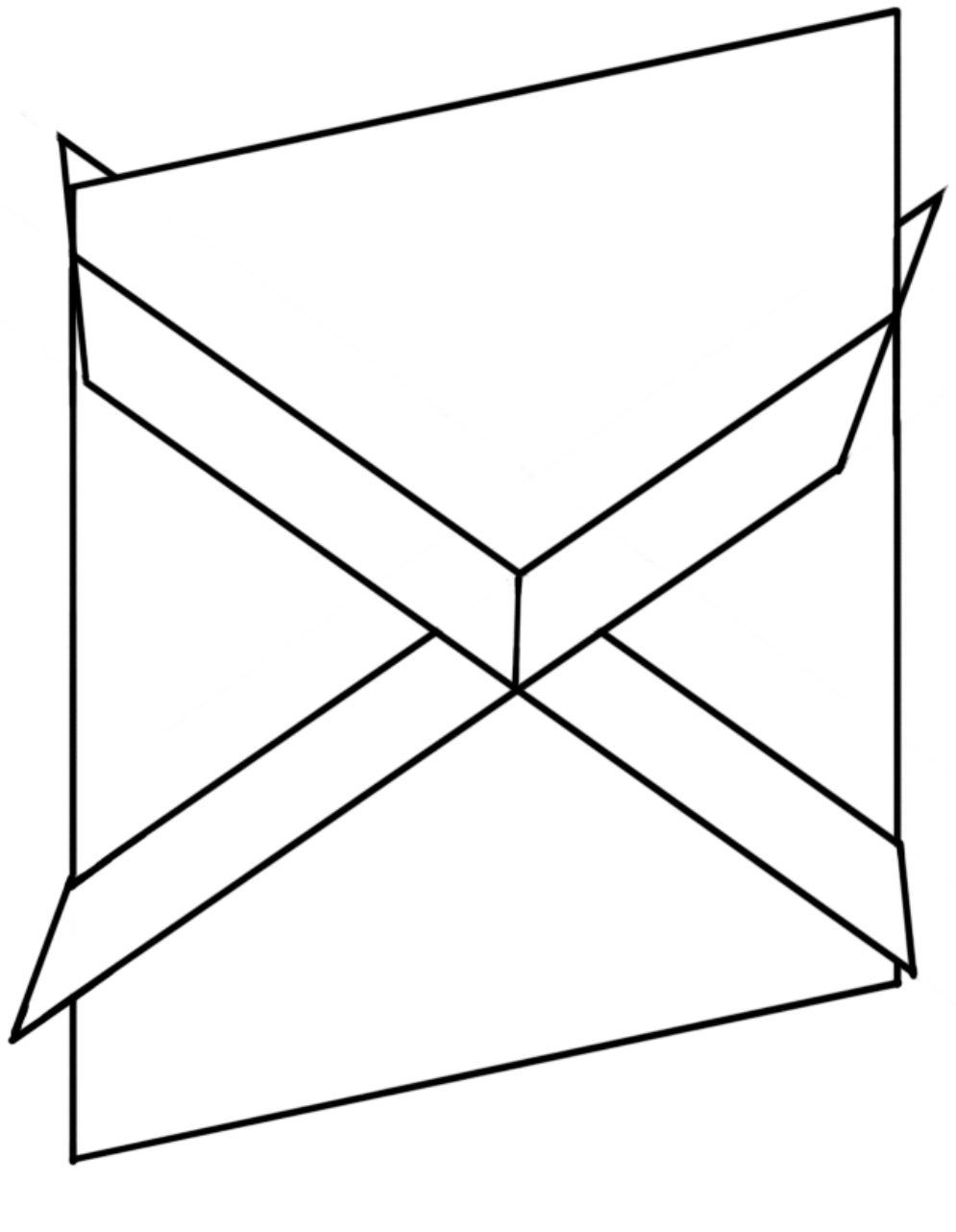}
        \caption{}
        \label{fig:proof_R3}
    \end{subfigure}
    \caption{Non-generic occurences: (a) The Whitney umbrella and the $R_1$ move. We can assume that the singularity is transverse to the projection. (b) The double point set of the projection. When crossing the double point set with a plane (locally), we get crossings in the projection. When $\lbrace t \rbrace \times \mathbb{T}^2 $ is tangent to the double point set, we get an $R_2$ move. (c) A triple point locally modelled by three intersecting planes. This corresponds to the $R_3$ move.}
    \label{fig:proof_R1_R3}
\end{figure}

\begin{figure}[hbtp]
    \centering
    \begin{subfigure}[b]{4.5cm}
        \includegraphics[width=\textwidth]{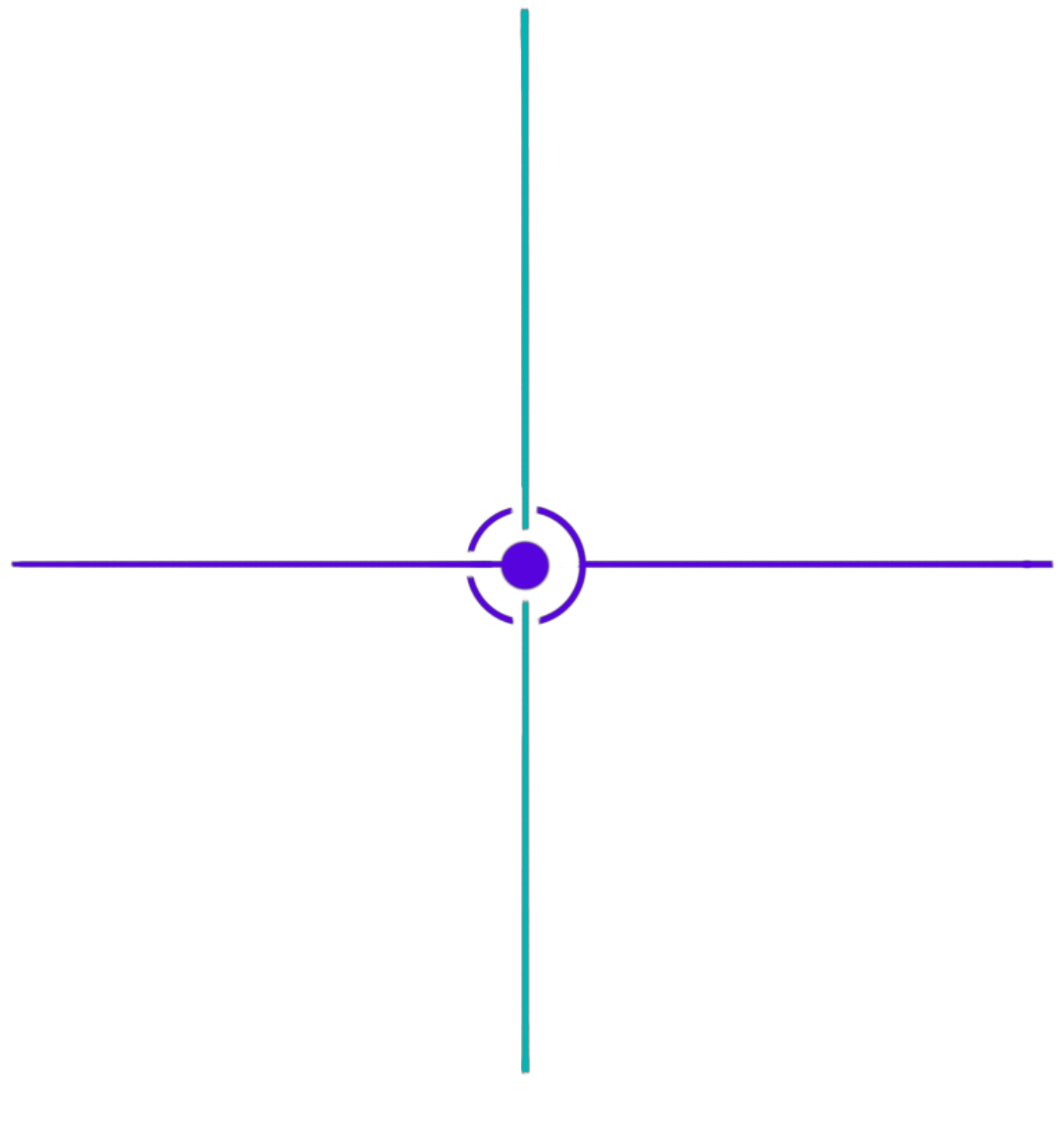}
        \caption{}
        \label{fig:proof_R4}
    \end{subfigure}
    \hspace{1.5cm}
    \begin{subfigure}[b]{3cm}        
        \includegraphics[width=\textwidth]{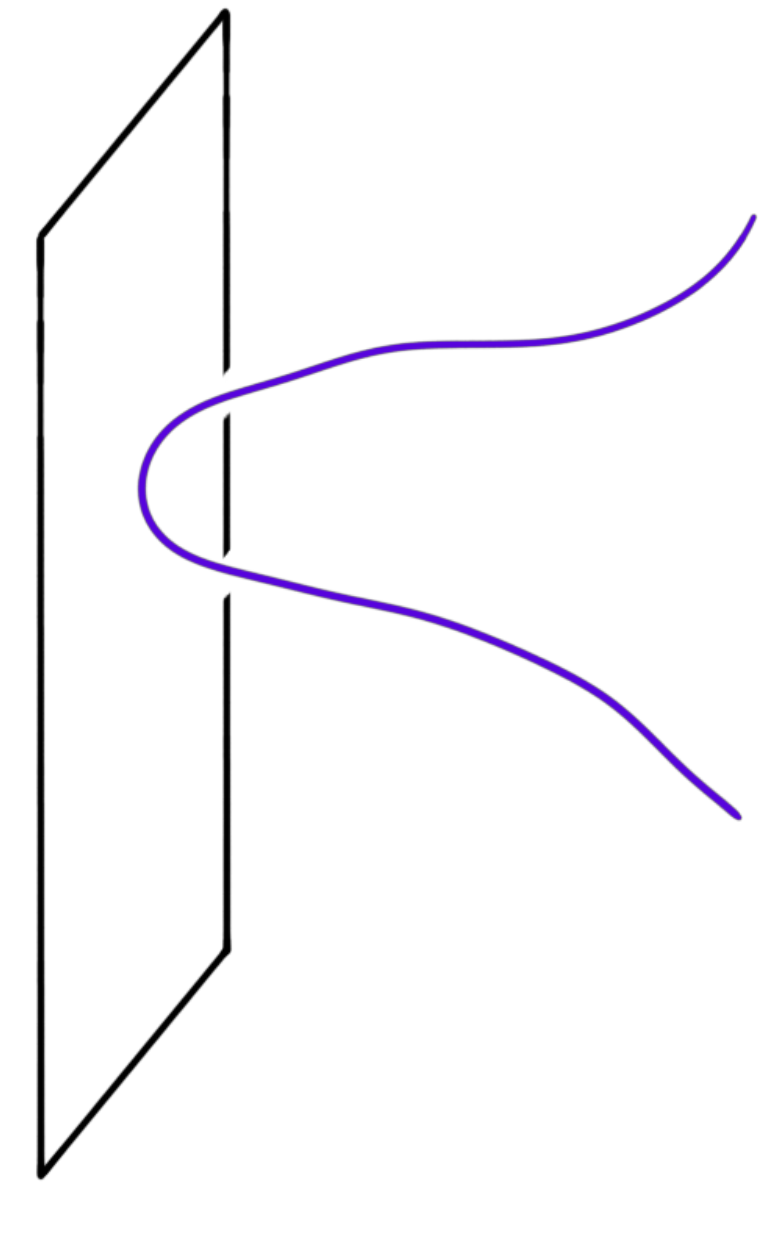}
        \caption{}
        \label{fig:proof_R5}
    \end{subfigure}
    \caption{Non-generic occurences: (a) A strand goes between the front part and back part of another that intersects the $t$-$N$-face. This requires the $R_4$ move. (b) The $t$-$N$-face is tangent to a strand. This corresponds to the $R_5$ move.}
    \label{fig:proof_R4_R5}
\end{figure}

Assume now that the special value corresponds to the case (\ref{item:special_t_4}). The projections in $\lbrace t_1 \rbrace \times \mathbb{T}^2$ and $\lbrace t_2\rbrace \times \mathbb{T}^2$ differ by an $R_4$ move. Indeed, such a double point appears when a strand passes between the front part and back part of another strand as shown in Figure \ref{fig:proof_R4}.

As the $t$-$N$-face passes through a point of tangency with the component of the link, as given in (\ref{item:special_t_5}), the projection undergoes an $R_5$ move. This happens exactly when a strand goes through the $N$-face as seen in Figure \ref{fig:proof_R5}.

The case (\ref{item:special_t_6}) occurs precisely when the projected curve goes through the edges of the square delimiting the $2$-torus. This corresponds to the $R_6$ move.

Case (\ref{item:special_t_7}) occurs when a crossing goes through the edges of the square delimiting the $2$-torus. This is the $R_7$ move.

 The special value corresponds to the case (\ref{item:special_t_8}) when the curve goes through the intersection of the edges $l_1$ and $l_2$ of the square representing the $2$-torus. This is the $R_8$ move.

Finally, the case (\ref{item:special_t_9}) occurs when an $N$-point goes through the edges of the square representing the $2$-torus. This corresponds to the $R_9$ move.
\end{proof}

\begin{remark}\label{rmk:dia_has_tridia_R_equivalent}
    Proposition \ref{prop:ambi_iso_links_imp_equi_dia} implies that in Lemma \ref{lem:dia_has_tridia}, although there are infinitely many ways of embedding a diagram, and thus, obtaining a tridiagram from it, any two of those tridiagrams are $R$-equivalent.
\end{remark}

\begin{proposition}\label{prop:t3_equi_equi_t2_equi}
    Two $T^3$-equivalent links in the $3$-torus $\Gamma$ and $\Gamma'$ possess $T^2$-equivalent tridiagrams, up to $R$-equivalence.
\end{proposition}

\begin{proof}
It suffices to prove that a $3$-torus twist corresponds to a $2$-torus twist performed on a `well-chosen' diagram. Suppose therefore, without loss of generality, that $\Gamma$ and $\Gamma'$ are connected by a single $3$-torus twist. The general case is proven by iterating the method used for this particular case.

The homeomorphisms $\phi_j$'s of Definition \ref{def:3_torus_twists} are the identity on one of the circles of the $3$-torus, and thus, correspond to the homeomorphisms $\psi_k$'s of Definition \ref{def:2_torus_twists}. More precisely, $\phi_1$, $\phi_3$ and $\phi_5$ correspond to $\psi_1$, and $\phi_2$, $\phi_4$ and $\phi_6$ correspond to $\psi_2$.

Suppose that $\left\lbrace D_x, D_y, D_z \right\rbrace$ is a tridiagram of $\Gamma$, where $x$, $y$, $z$ indicates the axis of projection. Suppose that the $3$-torus twist $\phi_{j_0}$, with $j_0 \in \left\lbrace 1,\dots,6 \right\rbrace$ that connects $\Gamma$ and $\Gamma'$, is the identity when restricted along $t_0 \in \left\lbrace x,y,z\right\rbrace$ and corresponds to $\psi_{k_0}$, $k_0 \in \lbrace 1,2\rbrace$. The action of $\phi_{j_0}$ on $\Gamma$ is thus equivalent to the action of $\psi_{k_0}$ on $D_{t_0}$. The resulting diagram is called $D_{t_0,1}$. Thanks to Lemma \ref{lem:dia_has_tridia}, from $D_{t_0,1}$ one can obtain a tridiagram $T_1$. Remark \ref{rmk:dia_has_tridia_R_equivalent} ensures that although there are infinitely many tridiagrams that can be obtained from $D_{t_0,1}$, they are all $R$-equivalent. The tridiagram $T_1$ is actually $R$-equivalent to a tridiagram of $\Gamma'$. Thus, we obtained a sequence of tridiagrams that satisfy the requirements of Definition \ref{def:t2_equivalence}, up to $R$-equivalence.
\end{proof}

\begin{proposition}\label{prop:m3_equi_equi_m2_equi}
    Two $CM^3$-equivalent links in the $3$-torus $\Gamma$ and $\Gamma'$ possess $CM^2$-equivalent tridiagrams, up to $R$-equivalence.
\end{proposition}

\begin{proof}
    It suffices to show that a covering map of a link in the $3$-torus can be interpreted as a finite sequence of covering maps of diagrams. The general case is just an iteration of the proof of this particular case. Therefore, without loss of generality, assume that $\Gamma$ is a cover of $\Gamma'$ where the covering map is
    \[
\begin{array}{cccc}
  p: &  \mathbb{S}^1_1 \times \mathbb{S}^1_2 \times \mathbb{S}^1_3 & \longrightarrow & \mathbb{S}^1_1 \times \mathbb{S}^1_2 \times \mathbb{S}^1_3 \\
    &  (z_1,z_2,z_3) & \longmapsto & \left({z_1}^{h_{11}}{z_2}^{h_{21}}{z_3}^{h_{31}},{z_2}^{h_{22}}{z_3}^{h_{32}},{z_3}^{h_{33}}\right)
\end{array}, z_i \in \mathbb{C}, \|z_i\| = 1,
\]
where $0 \leqslant h_{ij} < h_{ii}$ for $j < i$.

From $p$, define $p_1$, $p_2$ and $p_3$ in the following way:
\[
\begin{array}{cccc}
  p_1: &  \mathbb{S}^1_1 \times \mathbb{S}^1_2 \times \mathbb{S}^1_3 & \longrightarrow & \mathbb{S}^1_1 \times \mathbb{S}^1_2 \times \mathbb{S}^1_3 \\
    &  (z_1,z_2,z_3) & \longmapsto & \left({z_1}^{h_{11}}{z_2}^{h_{21}},{z_2}^{h_{22}},z_3\right)
\end{array}, z_i \in \mathbb{C}, \|z_i\| = 1,
\]

\[
\begin{array}{cccc}
  p_2: &  \mathbb{S}^1_1 \times \mathbb{S}^1_2 \times \mathbb{S}^1_3 & \longrightarrow & \mathbb{S}^1_1 \times \mathbb{S}^1_2 \times \mathbb{S}^1_3 \\
    &  (z_1,z_2,z_3) & \longmapsto & \left(z_1{z_3}^{h_{31}},z_2,{z_3}^{h_{33}}\right)
\end{array}, z_i \in \mathbb{C}, \|z_i\| = 1,
\]

\[
\begin{array}{cccc}
  p_3: &  \mathbb{S}^1_1 \times \mathbb{S}^1_2 \times \mathbb{S}^1_3 & \longrightarrow & \mathbb{S}^1_1 \times \mathbb{S}^1_2 \times \mathbb{S}^1_3 \\
    &  (z_1,z_2,z_3) & \longmapsto & \left(z_1,z_2{z_3}^{h_{32}{h_{33}}^{-1}},z_3\right)
\end{array}, z_i \in \mathbb{C}, \|z_i\| = 1.
\]

One may notice that $p_3 \circ p_2 \circ p_1 = p$, and that each of the $p_i$'s corresponds to the identity when restricted to one of the circles $\mathbb{S}^1_j$'s. Therefore, from $p_1$, $p_2$ and $p_3$ define $\tilde{p}_1$, $\tilde{p}_2$ and $\tilde{p}_3$ in the following way:

\[
\begin{array}{cccc}
  \tilde{p}_1: &  \mathbb{S}^1_1 \times \mathbb{S}^1_2 & \longrightarrow & \mathbb{S}^1_1 \times \mathbb{S}^1_2 \\
    &  (z_1,z_2) & \longmapsto & \left({z_1}^{h_{11}}{z_2}^{h_{21}},{z_2}^{h_{22}}\right)
\end{array}, z_i \in \mathbb{C}, \|z_i\| = 1,
\]

\[
\begin{array}{cccc}
  \tilde{p}_2: &  \mathbb{S}^1_1 \times \mathbb{S}^1_3 & \longrightarrow & \mathbb{S}^1_1 \times \mathbb{S}^1_3 \\
    &  (z_1,z_3) & \longmapsto & \left(z_1{z_3}^{h_{31}},{z_3}^{h_{33}}\right)
\end{array}, z_i \in \mathbb{C}, \|z_i\| = 1,
\]

\[
\begin{array}{cccc}
  \tilde{p}_3: &  \mathbb{S}^1_2 \times \mathbb{S}^1_3 & \longrightarrow &  \mathbb{S}^1_2 \times \mathbb{S}^1_3 \\
    &  (z_2,z_3) & \longmapsto & \left(z_2{z_3}^{h_{32}{h_{33}}^{-1}},z_3\right)
\end{array}, z_i \in \mathbb{C}, \|z_i\| = 1.
\]

The maps $\tilde{p}_1$ and $\tilde{p}_2$ are covering maps of the $2$-torus. However, it is not the case for $\tilde{p}_3$. Nevertheless, one can compose $\tilde{p}_3$ with the map 
\[
\begin{array}{cccc}
    h: & \mathbb{S}^1_2 \times \mathbb{S}^1_3 & \longrightarrow &  \mathbb{S}^1_2 \times \mathbb{S}^1_3\\
    & (z_2,z_3) & \longmapsto & \left(z_2,{z_3}^{h_{33}}\right)
\end{array},
\]
where $h$ is regarded as a map that reparametrises the circle $\mathbb{S}^1_3$. Therefore, by setting $\hat{p}_3 = \tilde{p}_3 \circ h$, we have the covering map

\[
\begin{array}{cccc}
  \hat{p}_3: &  \mathbb{S}^1_{2'} \times \mathbb{S}^1_{3'} & \longrightarrow &  \mathbb{S}^1_{2'} \times \mathbb{S}^1_{3'} \\
    &  (z_2,z_3) & \longmapsto & \left(z_2{z_3}^{h_{32}},{z_3}^{h_{33}}\right)
\end{array}, z_i \in \mathbb{C}, \|z_i\| = 1,
\]
where $\mathbb{S}^1_{2'} \times \mathbb{S}^1_{3'} = h^{-1}\left(\mathbb{S}^1_2 \times \mathbb{S}^1_3\right)$.

Since we have $0 \leqslant h_{ij} < h_{ii}$ for $j < i$, $\tilde{p}_1$, $\tilde{p}_2$ and $\hat{p}_3$ satisfy the conditions on the covering maps of Definition \ref{def:cover_2}.

 Assume now that $\left\lbrace D_x, D_y, D_z \right\rbrace$ is a tridiagram of $\Gamma$, where $x$, $y$, $z$ indicates the axis of projection. The cover $p$ can be divided in the following way. Apply $\tilde{p}_1$ on $D_z$ to get a diagram $D_{z,1}$. Lemma \ref{lem:dia_has_tridia} ensures that from $D_{z,1}$, one can get a tridiagram $\left\lbrace D_{x,1}, D_{y,1}, D_{z,1}\right\rbrace$. There are infinitely many ways to obtain a tridiagram from $D_{z,1}$. However, one must make sure that the chosen tridiagram must be one onto which $\tilde{p}_2$ is indeed a covering map of a certain diagram. Once the tridiagram is well chosen, apply $\tilde{p}_2$ on $D_{y,1}$ to get a diagram $D_{y,2}$. Again, from $D_{y,2}$, one can get a well-chosen tridiagram $\left\lbrace D_{x,2}, D_{y,2}, D_{z,2}\right\rbrace$. Apply $\hat{p}_3$ on $D_{x,2}$ to get a diagram $D_{x,3}$. Finally, embed $D_{x,3}$ in the $3$-torus to obtain a tridiagram $\left\lbrace D_{x,3}, D_{y,3}, D_{z,3}\right\rbrace$. It is a tridiagram that is $R$-equivalent to a tridiagram of $\Gamma'$.
\end{proof}

We are now in a position to give a generalised Reidemeister theorem for TP tangles, with respect to the $U$-equivalence.

\begin{theorem} \textbf{(Generalised Reidemeister theorem)}\\
    Two 3-periodic tangles $K$ and $K'$ are equivalent if and only if, up to $a$-equivalence, any two of their respective tridiagrams $\lbrace D_i \rbrace_{i=1,2,3}$ and $\lbrace D_i^{\prime} \rbrace_{i=1,2,3}$ are equivalent.
\end{theorem}

\begin{proof}
    Suppose we have two equivalent tridiagrams $\lbrace D_i \rbrace_{i=1,2,3}$ and $\lbrace D_i^{\prime} \rbrace_{i=1,2,3}$. There exists a finite sequence of tridiagrams satisfying the conditions given in Definition \ref{def:equi_tridiagrams}. But the $R$-moves and $2$-torus isotopies, $2$-torus twists, $CM^2$-equivalence and $V$-equivalence respectively translate the $3$-torus isotopies, $3$-torus twists, $CM^3$-equivalence and change of orientation of projection. Thus, their associated embeddings $K_a$ and $K_a'$ are also equivalent. By Remark \ref{rmk:diagram_loss_of_geometry} $K_a$ and $K_a'$ are respectively $a$-equivalent to $K$ and $K'$.

   Conversely, consider $\Gamma$ and $\Gamma'$ two unit cells of $K$ and $K'$. By Definition \ref{def:u_equivalence}, there is a finite sequence of links in the $3$-torus $(\Gamma_k)_{k=0,\dots,n}$ that are connected by ambient isotopies, $T^3$-equivalence, $CM^3$-equivalence, such that $\Gamma_0$ is a unit cell of $K_a$ and $\Gamma_n$ is a unit cell of $K_a'$, where $K_a$ and $K_a'$ are $a$-equivalent to $K$ and $K'$. By Proposition \ref{prop:unit_cells_are_equi}, $\Gamma$ and $\Gamma'$ are respectively connected to $\Gamma_0$ and $\Gamma_n$ by a finite sequence of links in the $3$-torus that are connected by ambient isotopies, $T^3$-equivalence, $CM^3$-equivalence. By Proposition \ref{prop:ambi_iso_links_imp_equi_dia}, Proposition \ref{prop:t3_equi_equi_t2_equi} and Proposition \ref{prop:m3_equi_equi_m2_equi}, the $3$-torus isotopies, $T^3$-equivalence and $CM^3$-equivalence can be converted into $2$-torus isotopies, $T^2$-equivalence and $CM^2$-equivalence of a finite sequence of tridiagrams.
\end{proof}

\section{Classification problem}\label{sec:4}
In all diagrammatic theories of links in 3-manifolds, diagrams enable the definition of invariants. The primary use of these invariants is the classification of those links. With the notions of diagrams and tridiagrams described in Section \ref{sec:3}, we define here the \textit{crossing number} of TP tangles. Crossing number is the primary measure of complexity used to classify knots and links, and it is thus important to extend this notion to classify TP tangles. This can then serve as a basis for future work to construct more nuanced invariants.

As explained in the previous sections, the crossing number defined here is a $U$-equivalence invariant (a particular case of ambient isotopy in $3$-space), and is not necessarily an ambient isotopy invariant in its most general sense.

\subsection{Crossing number}
We now consider the crossing number, which is previously defined as the minimum number of crossings over all diagrams of a knot \cite{Murasugi1996_chap4}.
The crossing number is the main complexity ordering of finite knots and links.

Naturally, crossings of TP tangles depend on unit cells. We define first a crossing number with respect to a unit cell, or in other words, the crossing number of a link embedded in the $3$-torus. Then, by minimising this crossing number over all unit cells, one obtains an invariant of TP tangles.

\begin{definition}
Consider a 3-periodic tangle $K$, and a unit cell $\Gamma$ of $K$. To a tridiagram $T$ associated to $\Gamma$, one can associate an unordered triplet of numbers of crossings $(a,b,c)$. It is referred to as a \textit{triplet of crossings}. We define the \textit{crossing number of $\Gamma$}, denoted by $c(\Gamma)$, as the minimum of $c(T) = a^2 + b^2 +c^2$ among all $R$-equivalent (and $V$-equivalent) tridiagrams. A tridiagram realising $c(\Gamma)$ is a called a \textit{minimal tridiagram}.
\end{definition}

\begin{remark}
    We note that the circle symbol representing the curve passing through the front and back faces depicted in Figure \ref{fig:newsignsdiagram} is not considered a crossing.
\end{remark}

\begin{definition}
    For a given 3-periodic tangle $K$, we define the \textit{crossing number}, denoted by $c(K)$, as the minimum of all $c(\Gamma)$ among all unit cells $\Gamma$ of $K$. A triplet of crossings of a tridiagram realising $c(K)$ is called a \textit{minimum crossing number triplet}.
\end{definition}

\begin{figure}[ht]
\centering
\includegraphics[width=0.75\textwidth]{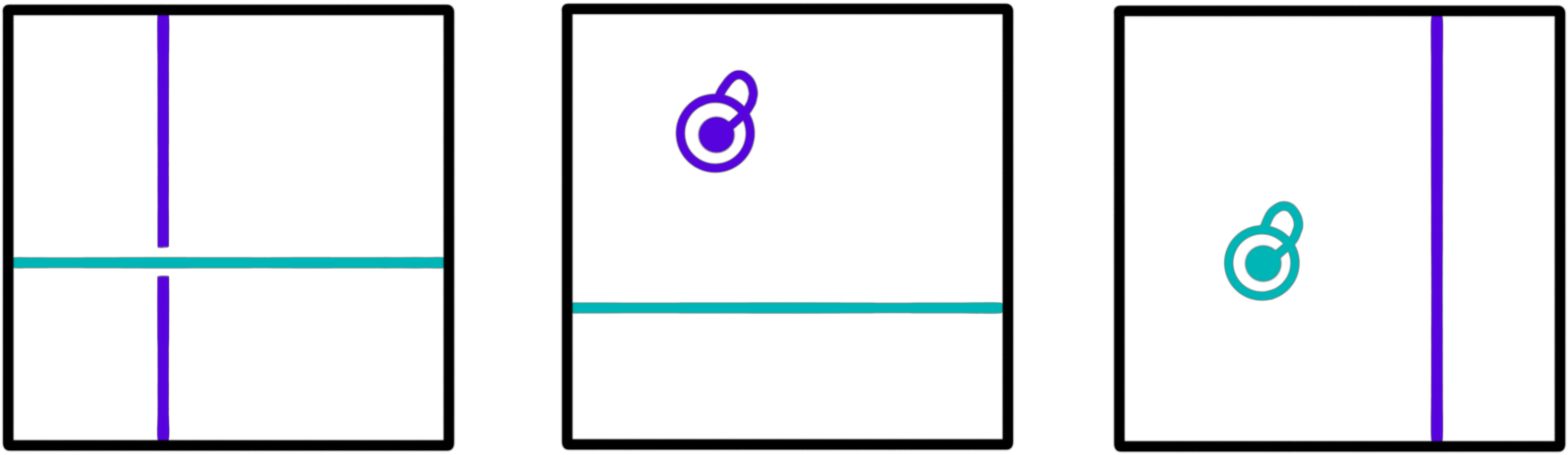}
\vskip\baselineskip
\includegraphics[width=0.75\textwidth]{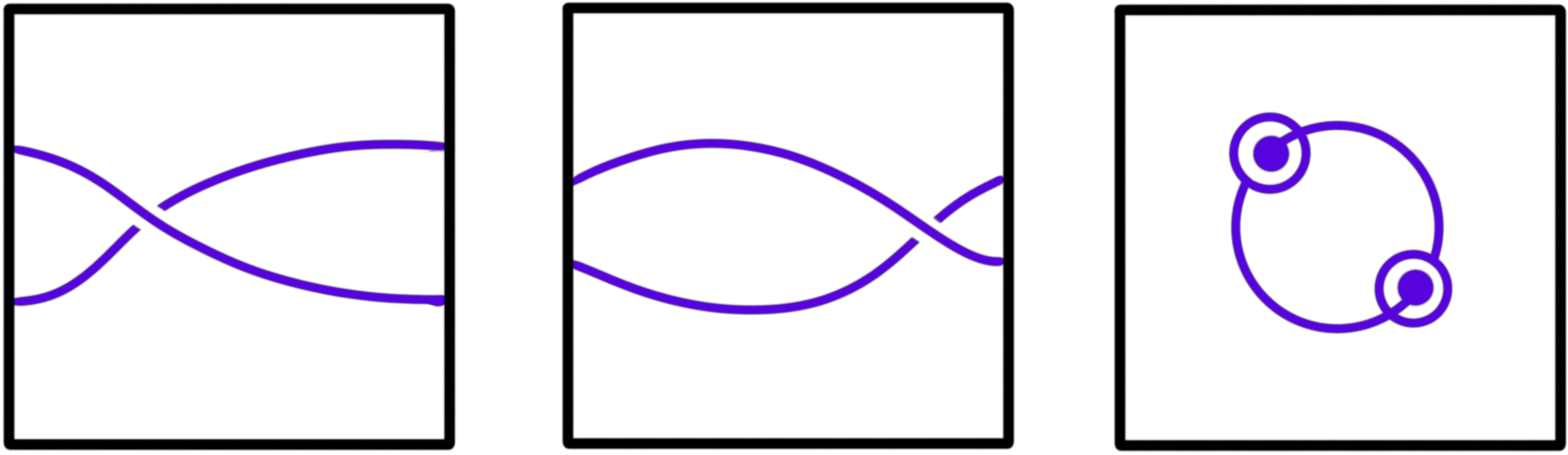}
\caption{Examples of 3-periodic tangles possessing one diagram having no crossings but with no tridiagrams with no crossings. On the top, a tridiagram of a bidirectional layer packing of straight lines. On the bottom, a tridiagram of a unidirectional packing of double helices.}
\label{fig:crossing_num_1}
\end{figure}

We now motivate our definition of crossing number of a TP tangle. Firstly, considering our characterisation of TP tangles via tridiagrams, three numbers of crossings are necessary to describe the entanglement. Indeed, some TP tangles, such as the ones shown in Figure \ref{fig:crossing_num_1}, have particular diagrams with no crossings, but cannot have a tridiagram with a $(0,0,0)$ crossing triplet where all three diagrams have no crossings simultaneously.

Secondly, there can be many ways of minimising the number of crossings of tridiagrams. However, we choose the sum of squares to bias the triplet towards a type $(a,a,a)$ where possible, favouring a more isotropic or symmetric configuration.

\begin{figure}[htbp]
\centering
    \begin{subfigure}[b]{\textwidth}
    \centering
        \includegraphics[width=0.75\textwidth]{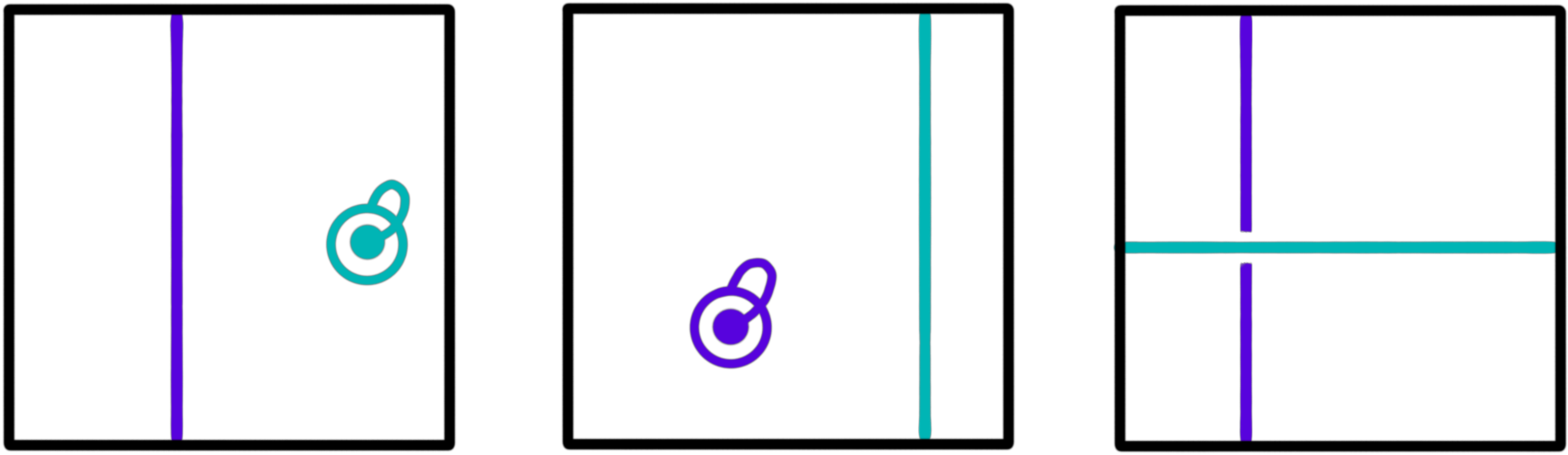} 
        \caption{}
        \label{fig:t_twi_compa_0}
    \end{subfigure}
    \vskip\baselineskip
    \begin{subfigure}[b]{6cm}
    \centering
        \includegraphics[width=\textwidth]{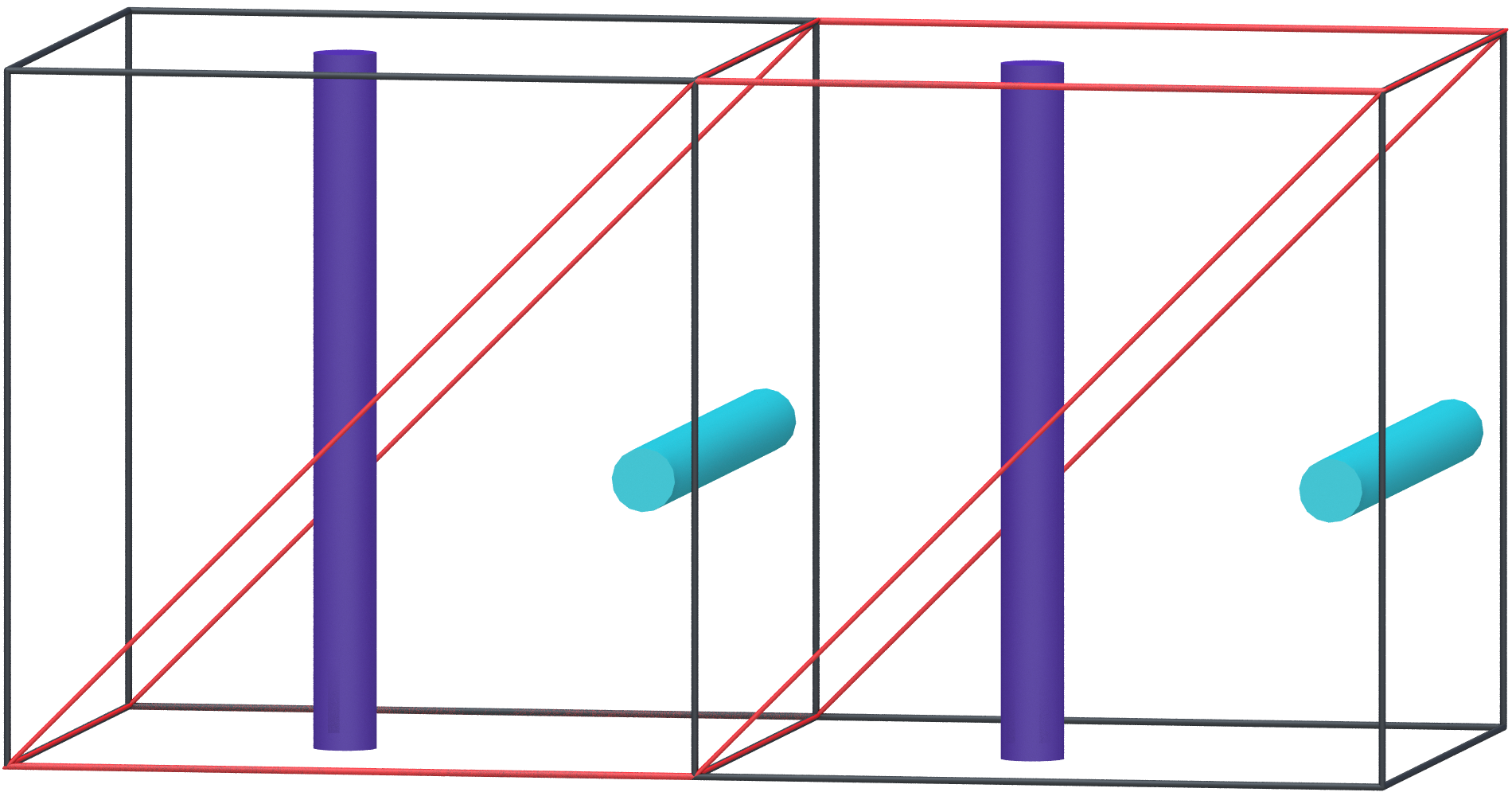}
        \caption{}
        \label{fig:t_twi_compa_1}
    \end{subfigure}
    \hspace{0.5cm}
    \begin{subfigure}[b]{3.25cm}
        \includegraphics[width=\textwidth]{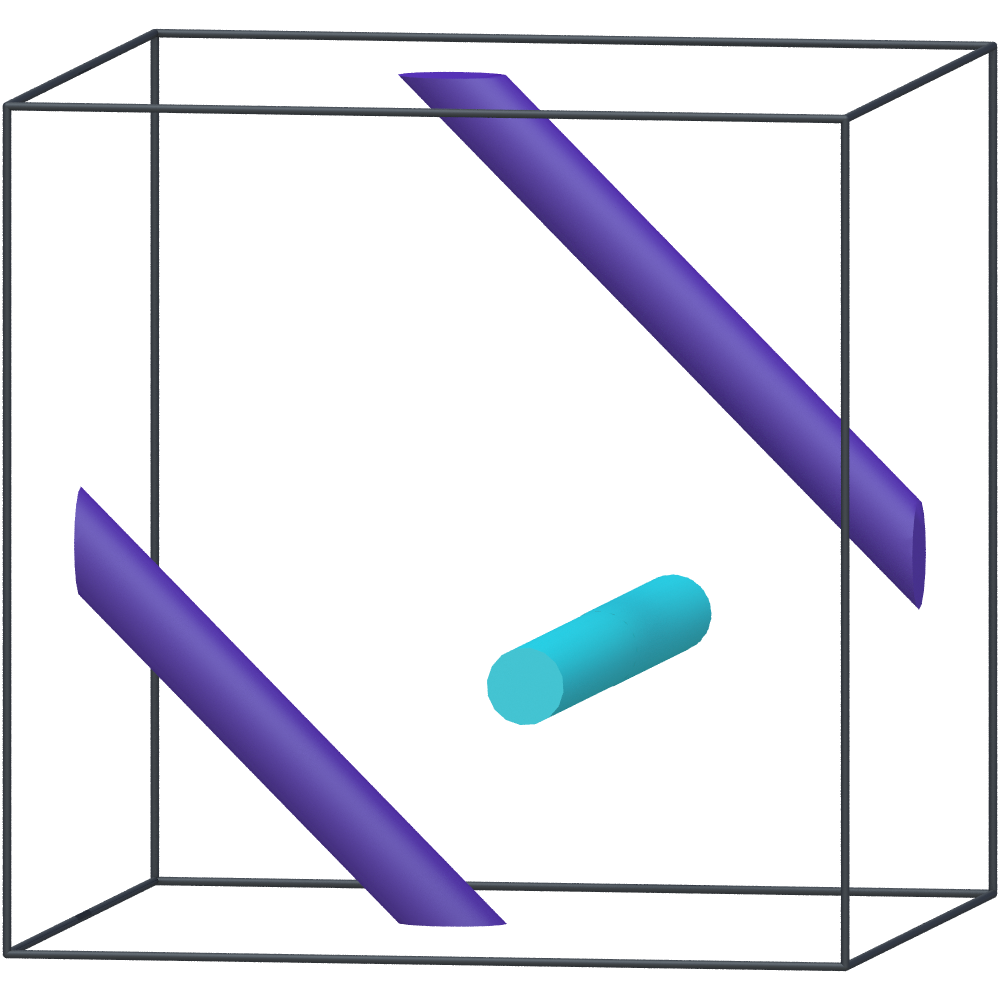}
        \caption{}
        \label{fig:t_twi_compa_2}
    \end{subfigure}
    \vskip\baselineskip
    \begin{subfigure}[b]{\textwidth}
    \centering
        \includegraphics[width=0.75\textwidth]{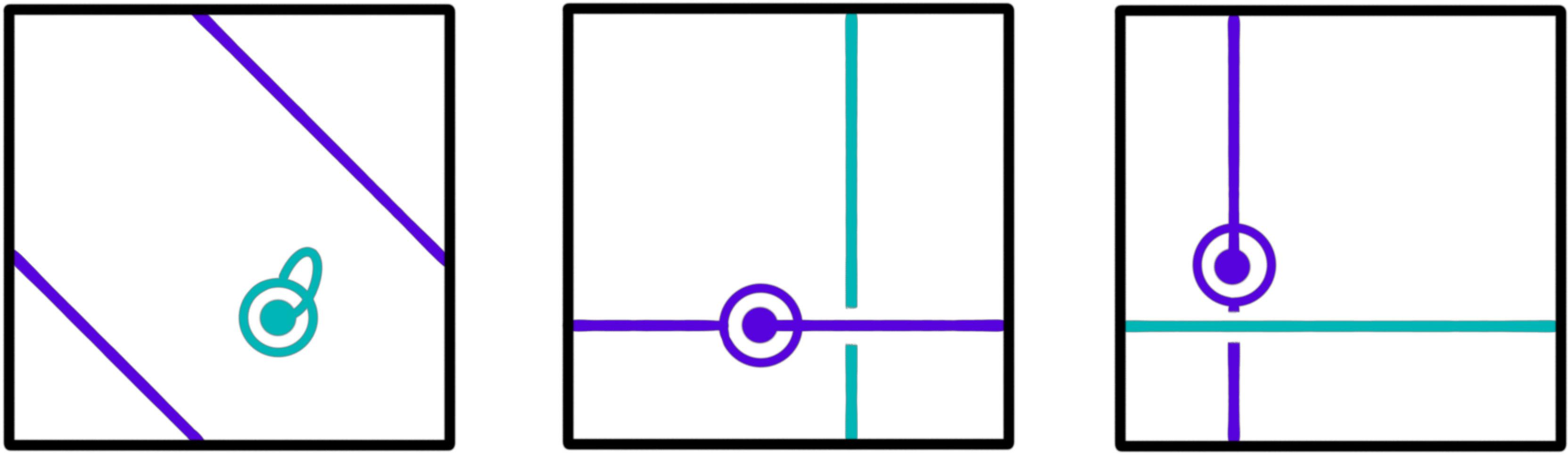}
        \caption{}
        \label{fig:t_twi_compa_3}
    \end{subfigure}
\caption{(a) A tridiagram of a bidirectional layer packing of straight lines, which has triplet of crossings $(0,0,1)$. (b) A change of basis of the lattice which corresponds to a torus twist in the unit cell. (c) The resulting unit cell after torus twist. (d) A tridiagram of the new unit cell, which now has triplet of crossings $(0,1,1)$.}
\label{fig:t_twi_compa}
\end{figure}

In general, symmetric embeddings of the structure tend to decrease the number of crossings. As an example, consider the TP tangle of Figure \ref{fig:t_twi_compa}, where we change the basis of its associated lattice. With the canonical basis, the axes of the curves of the embedding are parallel to the vectors delimiting the unit cell. The triplet of crossings is $(0,0,1)$ as given by the tridiagram in Figure \ref{fig:t_twi_compa_0}. By changing the basis, we `break' the symmetry of the TP tangle, and a crossing is added as seen in Figure \ref{fig:t_twi_compa_3}. The resulting triplet of crossings is $(0,1,1)$.

\begin{remark}
    This phenomenon is proper to TP tangles and does not occur in links in the 3-torus or DP tangles. Indeed, regarded as links in the 3-torus, the two unit cells of Figure \ref{fig:t_twi_compa_1} and \ref{fig:t_twi_compa_2} are not isotopic and thus the comparison of number of crossings is irrelevant. In the case of DP tangles, torus twists do not change the number of crossings as explained in \cite{grishanovmeshkov2007}.
\end{remark}

\begin{example}
Figure \ref{fig:pi_plus_uc_and_tridia} and Figure \ref{fig:sigma_plus_uc_and_tridia} show tridiagrams of the $\Pi^{+}$ and $\Sigma^{+}$ rod packings that we introduced in Figure \ref{fig:TPT}. From the tridiagrams, one can understand that the $\Pi^{+}$ rod packing has minimum crossing number triplet $(4,4,4)$, and $\Sigma^{+}$ has $(8,8,8)$. The cubic symmetry of the rod packings is exemplified in the equal crossing numbers from different projections. The tridiagrams give a simple encoding of these 3-dimensional structures, with the crossing triplet giving a measure of complexity of the structure.

    \begin{figure}[hbtp]
    \centering
    \includegraphics[width=0.8\textwidth]{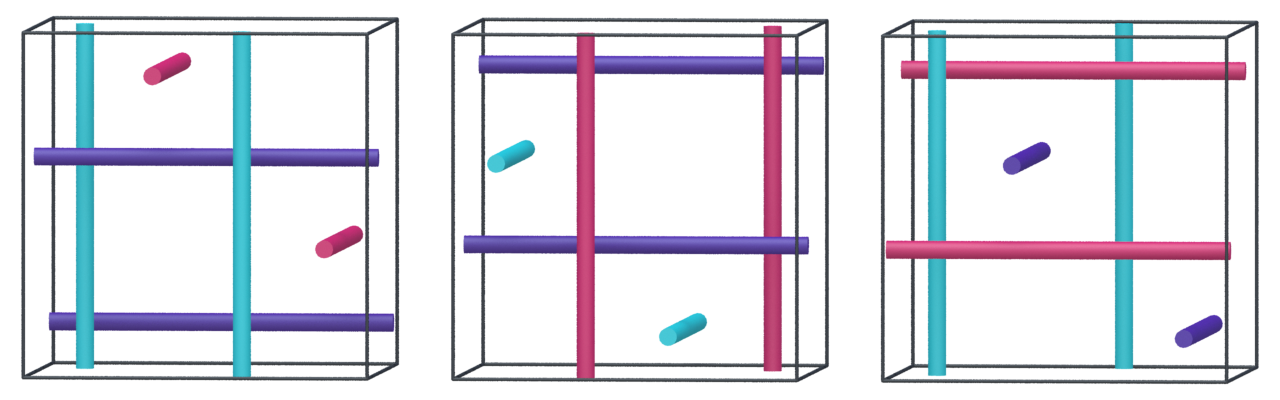}
    \vskip\baselineskip
    \includegraphics[width=0.75\textwidth]{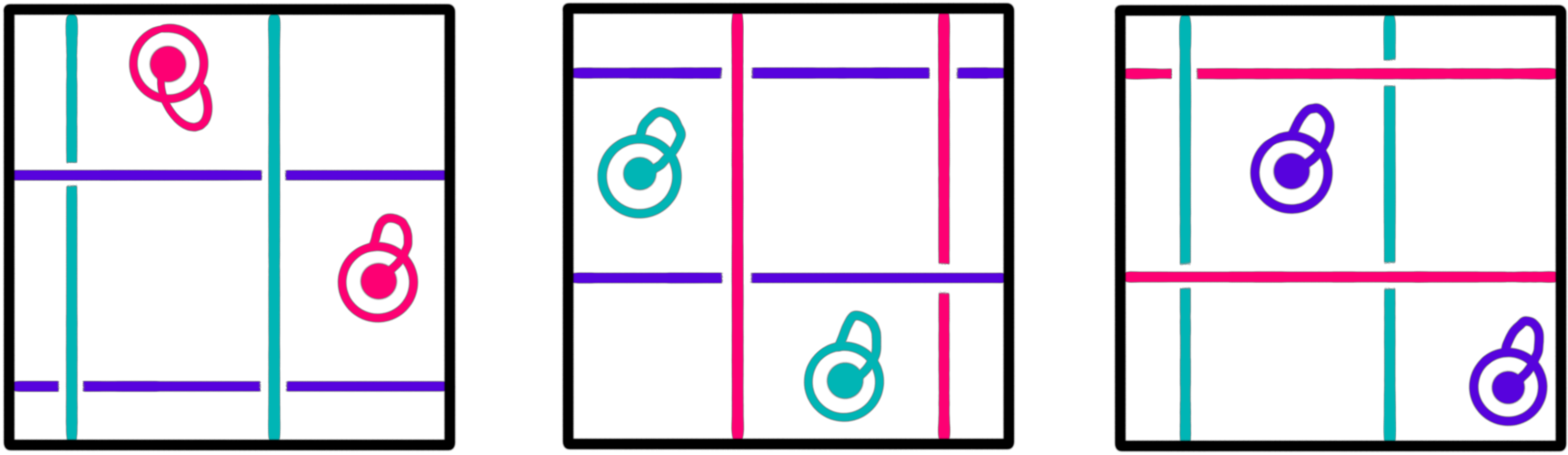}
    \caption{(Top) Unit cell of the $\Pi^{+}$ rod packing, presented in Figure \ref{fig:TPT} left, viewed along three non-coplanar axes. (Bottom) Tridiagram of the $\Pi^{+}$ rod packing. The minimum crossing number triplet is $(4,4,4)$.}
    \label{fig:pi_plus_uc_and_tridia}
\end{figure}

\begin{figure}[hbtp]
    \centering
    \includegraphics[width=0.8\textwidth]{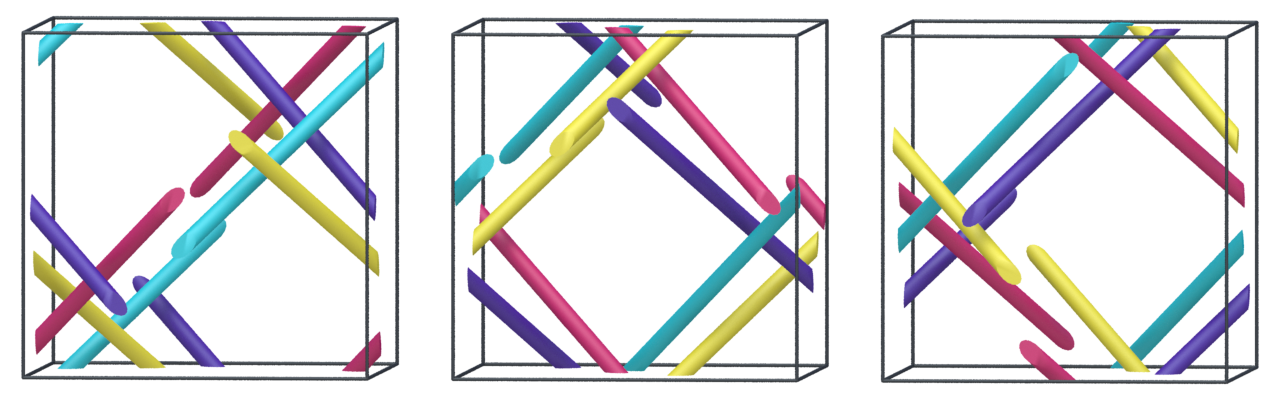}
    \vskip\baselineskip
    \includegraphics[width=0.75\textwidth]{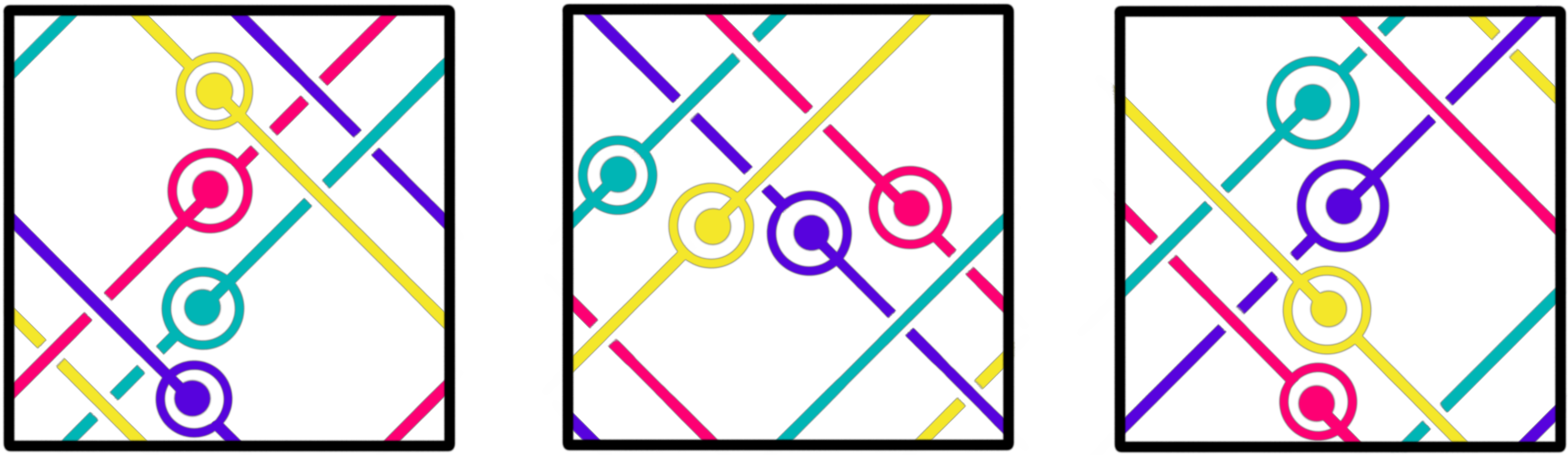}
    \caption{(Top) Unit cell of the $\Sigma^{+}$ rod packing, presented in Figure \ref{fig:TPT} right, viewed along three non-coplanar axes. (Bottom) Tridiagram of the $\Sigma^{+}$ rod packing. The minimum crossing number triplet is $(8,8,8)$.}
    \label{fig:sigma_plus_uc_and_tridia}
\end{figure}

\end{example}

\subsection*{Acknowledgements}
This work is funded by the Deutsche Forschungsgemeinschaft (DFG, German Research Foundation) - Project number 468308535, and partially supported by RIKEN iTHEMS.
We also acknowledge the software Houdini/SideFX.

%
%
%
\bibliographystyle{gtart}

\bibliography{biblio.bib}



\end{document}